\documentclass[a4paper,11pt,english]{article}
\usepackage{babel}
\usepackage{authblk} 
\usepackage[latin1]{inputenc}
\usepackage[pagewise]{lineno}
\usepackage[percent]{overpic} 

\author[1]{Juli\'an L\'opez-G\'omez\thanks{Corresponding author: jlopezgo@ucm.es}}
\author[1]{Alejandro Sahuquillo\thanks{alejsahu@ucm.es}}

\author[2]{Andrea Tellini\thanks{andrea.tellini@upm.es}}
\affil[1]{\small Universidad Complutense de Madrid\\ Departamento de An\'alisis Matem\'atico y Matem\'atica Aplicada\\ Plaza de las Ciencias 3\\ 28040 Madrid, Spain}
\affil[2]{\small Universidad Polit\'ecnica de Madrid\\ E.T.S.I.D.I.\\ Departamento de Matem\'atica Aplicada a la Ingenier\'ia Industrial\\ Ronda de Valencia 3\\ 28012 Madrid, Spain}

\title{\textbf{Large positive solutions for a class of 1-D diffusive logistic problems with general boundary conditions}}

\date{\today}

\usepackage{amsmath,amsfonts,amssymb}
\usepackage{yhmath,mathrsfs,yfonts,arcs}
\usepackage{amsthm}
\usepackage{multirow}
\usepackage{graphicx}
\usepackage[percent]{overpic} 
\usepackage{paralist}
\usepackage[titletoc,title]{appendix}
\usepackage{cite}
\usepackage{leftidx}
\usepackage{soul}
\usepackage{color}
\usepackage{bbm}
\usepackage{mathtools}
\mathtoolsset{showonlyrefs}
\usepackage{subcaption}
\usepackage{enumitem}
\usepackage{hyperref}

\theoremstyle{plain}
\newtheorem{theorem}{Theorem}[section]
\newtheorem{proposition}[theorem]{Proposition}
\newtheorem{lemma}[theorem]{Lemma}
\newtheorem{corollary}[theorem]{Corollary}

\theoremstyle{definition}

\theoremstyle{remark}
\numberwithin{equation}{section}

\usepackage[margin=0.80in]{geometry}

\newcommand\Item[1][]{%
	\ifx\relax#1\relax  \item \else \item[#1] \fi
	\abovedisplayskip=0pt\abovedisplayshortskip=0pt~\vspace*{-\baselineskip}}

\newcommand{\field}[1] {\mathbb{#1}}

\newcommand{\R}{\field{R}}

\def\b{\beta}
\def\e{\varepsilon}
\def\D{\Delta}
\def\d{\delta}

\def\l{\lambda}

\def\o{\omega}
\def\O{\Omega}
\def\p{\partial}
\def\r{\rho}

\def\t{\theta}

\def\v{\varphi}

\def\un{\underline}
\def\ua{\uparrow}
\def\da{\downarrow}

\newcommand{\mc}{\mathcal}

\makeatletter
\newcommand{\pushright}[1]{\ifmeasuring@#1\else\omit\hfill$\displaystyle#1$\fi\ignorespaces}
\newcommand{\pushleft}[1]{\ifmeasuring@#1\else\omit$\displaystyle#1$\hfill\fi\ignorespaces}
\makeatother

\begin{document}
\maketitle

\begin{abstract}
The first goal of this paper is to establish the existence  of a positive solution for the singular boundary value problem \eqref{eq:1.1}, where $\mc{B}$ is a general boundary operator of Dirichlet, Neumann or Robin type, either classical or non-classical; in the sense that, as soon as
$\mc{B}u(0)=-u'(0)+\b u(0)$, the coefficient $\b$ can take any real value, not necessarily $\b\geq 0$ as in the classical Sturm--Liouville theory.   Since the function $f(u):=au^p -\l u$, $u\geq 0$, is not increasing if $\l>0$, the uniqueness of the positive solution of \eqref{eq:1.1} is far from obvious, in general,  even for the simplest case when $a(x)$ is a positive constant. The second goal of this paper is to establish the uniqueness of the positive solution of
\eqref{eq:1.1} in that case. At a later stage, denoting by $L_\l$ the unique positive solution of \eqref{eq:1.1} when $a(x)$ is a positive constant, we will characterize  the  point-wise behavior of $L_\l$ as $\l\to \pm \infty$. It turns out that any positive solution of \eqref{eq:1.1} mimics the behavior of $L_\l$ as $\l \to \pm\infty$. Finally, we will establish  the uniqueness of the positive solution of \eqref{eq:1.1} when $a(x)$ is non-increasing in $[0,R]$, $\l\geq 0$, and $\b<0$ if $-u'(0)+\b u(0)=0$.
 \end{abstract}


\smallskip
\noindent \textbf{Acknowledgements:} This work has been supported by the Ministry of Science and Innovation of Spain under Research Grant PID2024-155890NB-I00. Andrea Tellini has been also supported by the Ram\'{o}n y Cajal program RYC2022-038091-I, funded by MCIN/AEI/10.13039/501100011033 and by the FSE+.

\section{Introduction}

\label{sec:1}

\noindent This paper studies the positive solutions of the one-dimensional singular boundary value problem
\begin{equation}
	\label{eq:1.1}
	\left \{ \begin{array}{l}
		-u''=\lambda u-a(x) |u|^{p-1}u\qquad \hbox{in} \;\; [0,R),\\[1ex]
		\mathcal{B}u(0)=0,\quad           u(R) = +\infty,          \end{array} \right.
\end{equation}
where $\l\in\R$ is regarded as a bifurcation parameter, $p >1$, $a\in\mc{C}([0,R];\R)$ satisfies $a(x)>0$ for all $x\in[0,R]$, and $\mathcal{B} \in \{ \mathcal{D}, \mathcal{N}, \mathcal{R}_{\beta} \}$ is a boundary operator of some of the following types: either $\mathcal{D}u(0) = u(0)$, or
$$
   \mathcal{N}u(0) =\frac{\p u}{\p n}(0)=-u'(0),
$$
where $\frac{\p}{\p n}$ denotes the exterior normal derivative at $x=0$, or
$$
   \mathcal{R}_{\beta}u(0) = \frac{\p u}{\p n}(0) + \beta u(0)=-u'(0)+\b u(0),
$$
for some $\beta \in \R$. Thus, $\mc{B}$ can be either the Dirichlet operator, or the Neumann operator, or a Robin operator for some value of $\b$ not necessarily positive as in the classical Sturm--Liouville theory. Note that $\mc{R}_0=\mc{N}$. The singular boundary condition at $x=R$ is understood as
$$
 \lim_{x\ua R}u(x)=+\infty.
$$
As we are focusing our attention only on positive solutions, the absolute value can be omitted in \eqref{eq:1.1},
though for some arguments in this paper  it is helpful to consider the nonlinearity of \eqref{eq:1.1}
as being odd.
\par
The most pioneering contributions to this field were given by Loewner and Nirenberg \cite{LoNi},
Kondratiev and Nikishin \cite{KoNi}, Bandle and Marcus \cite{BaMa90,BaMa91,BaMa98}, Marcus and V\'{e}ron
\cite{MaVe97,MaVe03,MaVe04}, and V\'{e}ron \cite{Ve92}, where  some multidimensional counterparts of \eqref{eq:1.1} of the type
\begin{equation}
\label{i.2}
  \left\{ \begin{array}{ll} \D u = a(x) f(u) & \quad \hbox{in}\;\;\O,\\
  u=\infty& \quad \hbox{on}\;\; \p\O,\end{array} \right.
\end{equation}
were analyzed. In \eqref{i.2}, $\O$ is a smooth bounded domain of $\R^N$, $N\geq 1$, and $f$ is a nonlinear positive increasing function satisfying the classical condition of Keller \cite{Ke} and Osserman \cite{Os}, to guarantee the existence of at least a positive solution for \eqref{i.2}.  These papers treat the  case when $a\in\mc{C}(\bar\O)$ is a positive function separated away from zero. The case when $a(x)$ vanishes somewhere on $\bar\O$ was treated some years later by G\'{o}mez-Re\~{n}asco \cite{Go99}, G\'{o}mez-Re\~{n}asco and L\'{o}pez-G\'{o}mez \cite{GoLo}, Du and Huang \cite{DuHu}, Garc\'{\i}a-Meli\'{a}n, Letelier-Albornoz and Sabina de Lis \cite{GLS}, and L\'{o}pez-G\'{o}mez \cite{LG00,LG03,LG06}. The monograph of L\'{o}pez-G\'{o}mez \cite{LG15} contains a rather exhaustive treatment of the theory developed until 2015. The most important open problem in this field is to ascertain whether, or not, the singular boundary value problem \eqref{i.2} has a unique positive solution. Up to the best of our knowledge, the sharpest uniqueness and multiplicity results available for \eqref{i.2} are those of L\'{o}pez-G\'{o}mez and Maire \cite{LoM16,LoM17,LoM19} and  L\'{o}pez-G\'{o}mez, Maire and V\'{e}ron \cite{LMV}, where the concept of super-additivity introduced by Marcus and V\'{e}ron  \cite{MaVe04} was strengthened and the technical device to get uniqueness
based on the maximum principle of L\'{o}pez-G\'{o}mez \cite{LG07} was invoked systematically even in the absence of radial symmetry.
\par
According to Corollary 3.3 of L\'{o}pez-G\'{o}mez and Maire \cite{LoM18}, when $f\in\mc{C}^1[0,\infty)$
is increasing, $f(0)=0$, and
\begin{equation*}
  \mc{T}(x)\equiv \frac{1}{\sqrt{2}}\int_0^{+\infty} \frac{d\tau }{\sqrt{\int_0^\tau f(x+t)\,dt}}<\infty
\end{equation*}
for some $x>0$, then, for every $R>0$, the singular boundary value problem
\begin{equation}
\label{i.3}
	\left \{ \begin{array}{l}
		u''=f(u) \qquad \hbox{in} \;\; [0,R),\\[1ex]
		u'(0)=0,\quad       u(R) = +\infty,          \end{array} \right.
\end{equation}
possesses a unique positive solution.  The classical condition of Keller \cite{Ke} and Ossermann \cite{Os} requires the much stronger assumption that $T(x)<\infty$ for all $x>0$.
\par
However, thanks to Theorem 3.5 of \cite{LoM18}, it is already known that destroying the monotonicity of $f(u)$ on a compact set with arbitrarily small measure can provoke the existence of an arbitrarily large number of positive solutions of \eqref{i.3}. Thus, since, for any given constant $a>0$, the function
$$
  f(u):=au^p -\l u,\qquad u\geq 0,
$$
is decreasing for every $u\in [0,\left(\frac{\l}{pa}\right)^\frac{1}{p-1}]$ if $\l>0$, it is far from obvious whether, or not, the problem \eqref{eq:1.1} admits a unique positive solution.
\par
In the context of superlinear indefinite problems the most astonishing multiplicity results are those of L\'{o}pez-G\'{o}mez, Tellini and Zanolin \cite{LTZ} and L\'{o}pez-G\'{o}mez and Tellini \cite{LoT},  where it was established that for a class of weight functions $a(x)$ changing sign in $[0,R]$ the problem
\begin{equation}
\label{i.7}
	\left \{ \begin{array}{l}
		-u''=\lambda u-a(x) |u|^{p-1}u\qquad \hbox{in} \;\; [0,R),\\[1ex]
		u'(0)=0,\quad u(R)=+\infty,     \end{array} \right.
\end{equation}
can admit an arbitrarily large number of positive solutions as soon as $\l<0$ is sufficiently large. These
multiplicity results are provoked by the fact that $a(x)$ changes of sign. The associated bifurcation
diagramas were computed by L\'{o}pez-G\'{o}mez, Molina-Meyer and Tellini \cite{LMT14,LMT15}.
\par
The first result of this paper establishes the existence and the uniqueness of a positive solution for  \eqref{eq:1.1}  for every $\l\in\R$ in the special case when $a(x)$ is a positive constant, and
provides us with the sharp point-wise behavior of these solutions as $\l\ua +\infty$ and as $\l\da -\infty$. Although we are exploiting the fact that $a(x)$ is constant to make use of standard phase plane techniques, the reader should be aware that, since $\b$ can take any value, the boundary condition at $x=0$ is far from being of Sturm--Liouville type. Moreover, the associated nonlinearity is far from monotone if $\l>0$. These circumstances
add a significant value to the analysis carried out in this paper.

\begin{theorem}
\label{th:1.1} Assume that $a(x)$ is constant, $a(x)= a>0$ for all $x\in[0,R]$. Then:
\begin{enumerate}
\item[{\rm (a)}] For every $\l \in \mathbb{R}$, the problem \eqref{eq:1.1} has a unique positive solution, $L_{\l}\in \mc{C}^2([0,R);\R)$.
\item[{\rm (b)}] The map
$$
    \begin{array}{rcl} \R & \longrightarrow & \mc{C}([0,R);\R) \\
      \lambda & \mapsto & L_{\l}\end{array}
$$
is differentiable and point-wise increasing in $\l$.
\item[{\rm (c)}] One has that
$$
   \lim_{\l \uparrow \infty}L_{\l} = +\infty\;\;\hbox{uniformly in compact subsets of}\;\; [0,R),
$$
unless
$\mathcal{B} = \mathcal{D}$ at $x=0$, in which case $L_{\l}(0) = 0$ for all $\l \in \mathbb{R}$.
\item[{\rm (d)}] One has that
\begin{equation}
\lim_{\l \downarrow -\infty} L_{\l} = 0 \;\; \text{uniformly  in compact subsets of}\;\; [0,R).
\end{equation}
\end{enumerate}
\end{theorem}

Based on Theorem \ref{th:1.1}, we can also derive the following result for general weight functions $a(x)$ continuous and positive, not necessarily constant.

\begin{theorem}
\label{th:1.2}
Assume that $a(x)$ is continuous and positive in $[0,R]$. Then:
	\begin{enumerate}
		\item[{\rm (a)}] For every $\l\in\R$, \eqref{eq:1.1} possesses a minimal and a maximal positive solution, respectively denoted by $L^{\rm min}_{\l}$ and $L^{\rm max}_{\l}$, in the sense that any other positive solution $L$ of \eqref{eq:1.1} satisfies
	\begin{equation}
		\label{i.4}
		L^{\rm min}_{\l}(x) \leq L(x) \leq L^{\rm max}_{\l}(x) \quad \hbox{for all } x \in [0,R).
	\end{equation}
	\item[{\rm (b)}] As in Theorem \ref{th:1.1},
	\begin{equation}
		\label{i.5}
		\lim_{\l \uparrow +\infty} L^{\rm min}_{\l} = +\infty \;\;\hbox{uniformly in compact subsets of}\;\; [0,R),
	\end{equation}
unless  $\mathcal{B} = \mathcal{D}$ at $x=0$, in which case $L^{\rm min}_{\l}(0) = 0$ for all $\l \in \mathbb{R}$.
	\item[{\rm (c)}] As in Theorem \ref{th:1.1},
\begin{equation}
		\label{i.6}
		\lim_{\l \downarrow -\infty} L^{\rm max}_{\l} = 0 \;\; \hbox{uniformly in compact subsets of}\;\; [0,R).
	\end{equation}
\item[{\rm (d)}] If, in addition, $a$ is non-increasing in $[0,R]$, $\l \geq 0$ and $\mc{B} \in \{ \mc{D}, \mc{N}, \mc{R}_\b : \b <0 \}$, we have that $L^{\rm min}_{\l} = L^{\rm max}_{\l}$. Thus, in such cases, \eqref{eq:1.1} has a unique positive solution.
\end{enumerate}
\end{theorem}

Therefore, except for the uniqueness of the positive solution, the point-wise behavior of the positive solutions of \eqref{eq:1.1} as $\l\ua +\infty$ and $\l\da -\infty$ mimic the patterns of  Theorem \ref{th:1.1} with  $a(x)$ constant, though for $a(x)$ non-constant \eqref{eq:1.1} might have several positive solutions.
Theorem \ref{th:1.2} (d) delivers a new uniqueness result based on the monotonicity of $a(x)$, which is based on a sophisticated use of the maximum principle rather than on the analysis of the exact blow-up rate of the positive solutions at $R$, as it is usual in most of the existing uniqueness results (see \cite{LG15} and the
references therein, for an exhaustive survey on this issue).
\par
The distribution of this paper is the following. Section \ref{sec:2} delivers some preliminaries;  among them, the comparison principles that we will use throughout this paper. Section  \ref{sec:3} shows the existence and uniqueness parts of Theorem \ref{th:1.1}, which are based on a systematic use of phase plane techniques. The remaining properties of the large solution $L_\l$ will be established in Sections \ref{sec:4}--\ref{sec:6}. Finally, in Section \ref{sec:7}, we prove Theorem \ref{th:1.2} based on Theorem \ref{th:1.1} and some additional comparison principles.

\section{Preliminaries}
\label{sec:2}
In this section we present some preliminary results, relying on comparison principles, that will be later applied throughout this paper. In this section, we are treating the general case in which $a\in\mc{C}([0,R];\R)$ satisfies $a(x)>0$ for all $x\in[0,R]$.
\par
Subsequently, for every $c, d\in\R$, with $c<d$, and any given
$$
   \mathcal{B}_{c}, \mathcal{B}_d \in \{ \mathcal{D}, \mathcal{N}, \mathcal{R}_{\beta} : \, \beta \in \mathbb{R} \},
$$
we are denoting by $\mathscr{B}$ the boundary operator defined in $\{c,d\}=\p[c,d]$ as
\begin{equation}
	\mathscr{B}u(x) = \begin{cases}
		\mathcal{B}_cu(x)& \hbox{if} \;\; x = c,\\
		\mathcal{B}_du(x)& \hbox{if} \;\; x = d.
			\end{cases}
\end{equation}
In agreement with the notations introduced in Section 1, as soon as $\mc{B}_d=\mc{R}_\b$, we are taking
$$
  \mc{B}_d u(d)= \mc{R}_\b u(d) = u'(d)+\b u(d).
$$
Moreover, for  every  $q \in \mathcal{C}([c,d];\mathbb{R})$, we will denote by $\sigma_1[-D^2+q, \mathscr{B}, (c,d)]$ the principal eigenvalue of the problem
\begin{equation}
\label{2.1}
	\left \{ \begin{array}{ll}
		(-D^2+q) \varphi=\sigma \varphi & \hbox{in} \;\; [c,d],\\
		\mathscr{B} \varphi=0\qquad   &\hbox{on} \;\; \{c,d \},  \end{array} \right.
\end{equation}
i.e. the lowest one. According to Section~7.4 and Proposition~8.3 in \cite{LG13}, one has that:
\begin{itemize}
	\item $\sigma_1[-D^2+q, \mathscr{B}, (c,d)]$ is the only eigenvalue of \eqref{2.1} that admits a positive eigenfunction;
	\item any other eigenvalue of \eqref{2.1} is larger than $\sigma_1[-D^2+q, \mathscr{B}, (c,d)]$. Thus,
it is strictly dominant;
	\item if $q_1\gneq q_2$,  then
$$
    \sigma_1[-D^2+q_1, \mathscr{B}, (c,d)]>\sigma_1[-D^2+q_2, \mathscr{B}, (c,d)].
$$
In other words, it is monotone with respect to the potential $q$.
\end{itemize}
Some of these properties, for general mixed boundary conditions, go back to Cano-Casanova and L\'opez-G\'omez
\cite{CCLG}.
\par
Throughout this paper, a \emph{subsolution} of the problem
\begin{equation}
\label{2.2}
	\left \{ \begin{array}{ll}
		-u''=\lambda u-a(x) |u|^{p-1}u& \;\;  \hbox{in} \;\; [c,d],\\[1ex]
		\mathscr{B}u=0 & \;\;  \hbox{on}\;\;\{c,d\},  \end{array} \right.
\end{equation}
is a function $\underline{u}\in\mc{C}^2([c,d];\R)$ such that
$$
  \left\{ \begin{array}{ll} -\underline{u}''\leq \l\underline{u}-a(x)|\underline u|^{p-1}\underline u & \;\;  \hbox{in} \;\;[c,d],\\[1ex]
  \mathscr{B}\underline{u}\leq 0&\;\;  \hbox{on}\;\;\{c,d\}.\end{array}\right.
$$
In such a case, $\underline{u}$ is said to be a \emph{strict subsolution} if some of these inequalities is strict. Similarly,  a \emph{supersolution} of \eqref{2.2} is any function $\overline{u}\in\mc{C}^2([c,d];\R)$ such that
$$
  \left\{ \begin{array}{ll} -\overline{u}''\geq \l\overline{u}-a(x)|\overline{u}|^{p-1}\overline u & \;\;  \hbox{in} \;\;[c,d],\\[1ex]
  \mathscr{B}\overline{u}\geq 0&\;\;  \hbox{on}\;\;\{c,d\},\end{array}\right.
$$
and $\overline{u}$ is said to be a \emph{strict supersolution} if some of these inequalities is strict.
\par
The next uniqueness result will be extremely useful later.

\begin{proposition}
\label{pr:2.1}
Let $u$, $v$ be two positive solutions of the equation
\begin{equation}
\label{2.3}
-u''=\lambda u-a(x) u^p
\end{equation}
in $[c,d]$ such that
$$
   \mathscr{B}u = \mathscr{B}v\geq 0\quad \hbox{on}\;\; \{c,d\}.
$$
Then, $u=v$ in $[c,d]$.
\end{proposition}
\begin{proof}
Note that, necessarily,
$$
  u(x)>0 \;\; \hbox{and}\;\; v(x)>0 \;\; \hbox{for all}\;\; x\in (c,d).
$$
Since $u$ and $v$ solve \eqref{2.3},
\begin{equation}
\label{2.4}
-(u-v)'' = \lambda\left(u-v\right) -a(x)\left( u^p - v^p \right).
\end{equation}
Thus, introducing the auxiliary function
$$
    \varphi(t) := \left( t u+(1-t)v \right)^p,\qquad t \in [0,1],
$$
we have that
\begin{equation}
u^p - v^p = \varphi(1)- \varphi(0) = \int_0^1 \varphi'(t) \, dt = (u-v) \int_0^1 p \left( tu+(1-t)v\right)^{p-1} \, dt.
\end{equation}
Hence, setting
\begin{equation}
H(x) := p\int_0^1 \left( tu(x)+(1-t)v(x) \right)^{p-1} \, dt,\qquad x\in [c,d],
\end{equation}
the identity \eqref{2.4} can be equivalently rewritten as
\begin{equation}
-(u-v)'' = \lambda\left(u-v\right) -a(x)\, H(x) \left(u-v\right).
\end{equation}
Suppose $u\neq v$. Then, $w:=u-v\neq 0$ is an eigenfunction of the following linear spectral problem
\begin{equation*}
	\left \{ \begin{array}{ll}
		(- D^2 + a(x)H(x))w = \lambda w& \hbox{in} \;\; [c,d],\\
		\mathscr{B}w = 0& \hbox{on} \;\; \{c,d \}, \end{array} \right.
\end{equation*}
associated with the eigenvalue $\l$. Then, since the principal eigenvalue is dominant,
we find that
$$
  \lambda \geq \sigma_1[-D^2+a(x)H(x), \mathscr{B}, (c,d)].
$$
On the other hand, since $v(x)>0$ for all $x\in (c,d)$, we obtain that
\begin{equation}
	\label{2.5}
H(x) > p \; u^{p-1}(x) \int_0^1 t^{p-1}\, dt = u^{p-1}(x).
\end{equation}
Therefore, by the monotonicity of the principal eigenvalue with respect to the potential, it follows from
\cite[Th. 7.10]{LG13} (going back to L\'opez-G\'omez and Molina-Meyer  \cite{LGMM} and Amann and L\'opez-G\'omez \cite{ALG})
that
$$
\lambda \geq \sigma_1[-D^2+a(x)H(x), \mathscr{B}, (c,d)] > \sigma_1[-D^2+a(x) u^{p-1}(x), \mathscr{B}, (c,d)] \geq \lambda.
$$
Indeed, if $\mathscr{B}u=0$ on $\{c,d\}$, then, by the uniqueness of the principal eigenvalue,
$$
  \sigma_1[-D^2+a(x) u^{p-1}(x), \mathscr{B}, (c,d)]=\l,
$$
while in the case when $\mathscr{B}u\gneq 0$ on $\{c,d\}$, $u$ provides us with a positive strict supersolution of
\eqref{2.2} and, hence, thanks to \cite[Th. 7.10]{LG13}, it becomes apparent that
$$
  \sigma_1[-D^2+a(x) u^{p-1}(x), \mathscr{B}, (c,d)]>\l.
$$
This contradiction shows that $u=v$ and ends the proof.
\end{proof}

The following result shows that any pair of positive sub- and supersolution of \eqref{2.2} must be ordered.

\begin{proposition}
\label{pr:2.2}
Let $u, v \in \mc{C}^2([c,d];\R)$ be two positive functions such that $u$ is a supersolution of \eqref{2.2} and $v$ satisfies
\begin{equation}
\label{2.6}
  \left \{ \begin{array}{ll}
		-v''\leq \lambda v-a(x) |v|^{p-1}v& \;\;  \hbox{in} \;\; [c,d],\\[1ex]
		\mathscr{B}v\leq \mathscr{B}u  & \;\;  \hbox{on}\;\;\{c,d\}.  \end{array} \right.
\end{equation}
Suppose, in addition, that either $-u''\gneq \l u- au^p$ in $[c,d]$, or some of the inequalities in \eqref{2.6}
is strict. Then,  $u(x) > v(x)$ for all $x\in(c,d)$. Moreover,
\begin{equation}
\label{2.7}
   \left\{ \begin{array}{ll} u(c)>v(c) & \;\; \hbox{if}\;\; \mc{B}_c\in \{\mc{N},\mc{R}_\b:\b\in\R\},\\[1ex]
     u'(c) > v'(c)& \;\; \hbox{if}\;\; \mc{B}_c = \mathcal{D}\;\; \hbox{and}\;\; u(c)=v(c),\end{array}\right.\end{equation}
and
\begin{equation}
\label{2.8}
   \left\{ \begin{array}{ll}
   u(d)>v(d)& \;\; \hbox{if}\;\; \mc{B}_d\in \{\mc{N},\mc{R}_\b:\b\in\R\},\\[1ex]
   u'(d) < v'(d) & \;\; \hbox{if}\;\;\mc{B}_d = \mathcal{D}\;\; \hbox{and}\;\;   u(d)=v(d).
   \end{array}\right.
\end{equation}
\end{proposition}

Here, we are assuming that
$$
  u(x)>0\;\;\hbox{and}\;\; v(x)>0 \;\; \hbox{for all}\;\; x\in (c,d),
$$
though this result holds true even when $v$ is non-negative.

\begin{proof}
Since $u$ satisfies
\begin{equation}
	\left \{ \begin{array}{ll}
		(-D^2+a(x)u^{p-1}(x)-\l)u \geq 0& \hbox{in} \;\; [c,d],\\[1ex]
		\mathscr{B}u \geq 0 & \hbox{on} \;\; \{c,d \},  \end{array} \right.
\end{equation}
if $u$ is a strict supersolution of \eqref{2.2}, it follows from \cite[Th. 7.10]{LG13} that
$$
	\sigma_1[-D^2+a(x) u^{p-1}(x)-\l, \mathscr{B}, (c,d)] > 0.
$$
Should $u$ be a positive solution of \eqref{2.2}, then, by the uniqueness of the principal eigenvalue,
$$
	\sigma_1[-D^2+a(x) u^{p-1}(x)-\l, \mathscr{B}, (c,d)] = 0.
$$
Thus, in any circumstances,
$$
	\sigma_1[-D^2+a(x) u^{p-1}(x)-\l, \mathscr{B}, (c,d)] \geq  0.
$$
On the other hand, by the assumptions, we have that
$$
-(u-v)'' \geq  \lambda\left(u-v\right) -a(x)\left( u^p - v^p \right)
$$
and, hence, arguing as in the proof of Proposition \ref{pr:2.1}, it is apparent that $w:=u-v$ satisfies
\begin{equation}
\label{2.9}
	\left \{ \begin{array}{ll}
		(- D^2 + a(x)H(x)-\l)w \geq 0& \hbox{in} \;\; [c,d],\\[1ex]
		\mathscr{B}w \geq 0  & \hbox{on} \;\; \{c,d \}, \end{array} \right.
\end{equation}
with, at least, one of these inequalities strict. Consequently, $w$ is a strict supersolution of $- D^2 + a(x)H(x)-\l$ in $(c,d)$ for the
boundary operator $\mathscr{B}$. On the other hand, combining \eqref{2.5} with the monotonicity of the principal eigenvalue with respect to the potential, yields to
$$
\sigma_1[-D^2+a(x)H(x)-\l , \mathscr{B}, (c,d)] >  \sigma_1[-D^2+a(x) u^{p-1}(x)-\l, \mathscr{B}, (c,d)] \geq 0.
$$
Therefore, once again by \cite[Th. 7.10]{LG13}, it follows from \eqref{2.9}, that $w$ is strongly positive in
$[c,d]$, i.e. $w\gg 0$, in the sense that $w(x)>0$ for all $x\in (c,d)$, and \eqref{2.7} and \eqref{2.8} are satisfied
(see the proof of \cite[Th. 2.1]{LG13} if necessary).
\end{proof}

The following results use the comparison principle established by Proposition \ref{pr:2.2} to infer a comparison result related to the positive large solutions of equation \eqref{2.3}.

\begin{corollary}
\label{co:2.3}
Let $u, v\in\mc{C}^2([c,d);\R)$ be two positive functions such that
\begin{equation}
\label{2.10}
 \left\{ \begin{array}{l} -u''\geq \l u -a(x)u^p \\[1ex] -v''\leq \l v -a(x)v^p \end{array} \;\; \hbox{in}\;\; [c,d),\right.
\end{equation}
\begin{equation}
\label{2.11}
\mathcal{B}_cu(c) = 0 = \mathcal{B}_cv(c),
\end{equation}
\begin{equation}
\label{2.12}
u(d) = +\infty \;\; \hbox{and} \;\; v(d) \in(0,+\infty).
\end{equation}
Then, $u(x) > v(x)$ for all $x \in (c,d)$, and
\begin{equation}
\label{2.13}
   \left\{ \begin{array}{ll} u(c)>v(c) & \;\; \hbox{if}\;\; \mc{B}_c\in \{\mc{N},\mc{R}_\b:\b\in\R\},\\[1ex]
     u'(c) > v'(c)& \;\; \hbox{if}\;\; \mc{B}_c = \mathcal{D}.\end{array}\right.
\end{equation}
\end{corollary}
\begin{proof}
By the continuity of $v$, it follows from \eqref{2.12} that there exists $\delta_0 > 0$ such that
$$
   u(x) > v(x)\;\; \hbox{for all}\;\; x\in[d-\delta_0,d).
$$
Pick any $\delta \in (0,\delta_0)$ and consider $\mathcal{B}_{d-\delta}:=\mathcal{D}$. Then, $u$ is a supersolution of \eqref{2.2} in $[c,d-\d]$, $v$ is a subsolution of \eqref{2.2} in $[c,d-\d]$, and
$$
  \mc{B}_{d-\d}u(d-\d)=u(d-\d)> v(d-\d)=\mc{B}_{d-\d}v(d-\d).
$$
Thus, $\mathscr{B}u\gneq \mathscr{B}v$ on $\{c,d-\d\}$. Consequently, by \eqref{2.10} and
\eqref{2.11}, the  conclusion follows applying Proposition \ref{pr:2.2} in $[c,d-\d]$. Note that \eqref{2.13} is a direct consequence of  \eqref{2.7}.
\end{proof}

Similarly, the next result holds.

\begin{corollary}
\label{co:2.4}
Let $u\in \mc{C}^2((c,d);\R)$ and $v\in\mc{C}^2([c,d];\R)$ be two positive functions satisfying
\begin{equation*}
 \left\{ \begin{array}{l} -u''\geq \l u -a(x)u^p \\[1ex] -v''\leq \l v -a(x)v^p \end{array} \;\; \hbox{in}\;\; (c,d),\right.
\end{equation*}
and such that
\begin{equation}
	\label{2.14}
u(c) = u(d) = +\infty\;\; \hbox{and}\;\; v(c), v(d) \in (0,+\infty).
\end{equation}
Then, $u(x) > v(x)$ for all $x \in (c,d)$.
\end{corollary}
\begin{proof}
By the continuity of $v$, it follows from \eqref{2.14} that there exists $\delta_0 > 0$ such that
$$
   u(x)>v(x)\;\; \hbox{for all}\;\; x\in(c,c+\d_0]\cup[d-\d_0,d).
$$
Applying Proposition \ref{pr:2.2} in $[c+\d,d-\d]$ with $\delta \in (0,\delta_0)$ and
$\mathcal{B}_{c+\delta}=\mathcal{B}_{d- \delta}:=\mathcal{D}$, the proof is complete.
\end{proof}

\section{Proof of Part (a) of Theorem \ref{th:1.1}}
\label{sec:3}
Throughout the proof of Theorem \ref{th:1.1}, we are assuming that $a(x)$ is a positive constant $a$. Thus, phase portrait techniques can be applied to study the equation
\begin{equation}
\label{eq:3.1}
-u''=\lambda u-a |u|^{p-1}u
\end{equation}
to establish the existence and uniqueness of the positive solution of \eqref{eq:1.1}, $L_{\l}$, for
$$
   \mathcal{B} \in \{ \mathcal{D}, \mathcal{N}, \mathcal{R}_{\beta}\;:\; \beta \in \mathbb{R}\}.
$$
As the constructions in the phase plane vary with the boundary conditions, we will focus attention on each
of the boundary conditions separately.

\subsection{The Dirichlet case $\mathcal{B} = \mathcal{D}$}
\label{sec:3.1}

In this case, the problem \eqref{eq:1.1} becomes
\begin{equation}
	\label{eq:3.2}
	\left \{ \begin{array}{l}
		-u''=\lambda u-a |u|^{p-1}u\qquad \hbox{in} \;\; [0,R),\\[1ex]
		u(0)=0,\quad            u(R) = +\infty.          \end{array} \right.
\end{equation}
Setting  $v = u'$, we have that
$$
   v' = a |u|^{p-1}u- \l u.
$$
Thus,
\begin{equation}
vv' + \l u u' - a |u|^{p-1}uu' = 0,
\end{equation}
and, integrating the previous expression, we obtain that all the positive solutions of  \eqref{eq:3.1} lie on the level lines of the energy function defined by
\begin{equation}
	\label{eq:3.3}
E(u,v) :=\dfrac{1}{2}v^2+ \dfrac{\l}{2}u^2-  \dfrac{a}{p+1}|u|^{p+1}.
\end{equation}
In order to ascertain the phase diagram of \eqref{eq:3.1}, we will study the \emph{potential energy}
\begin{equation}
\varphi (u) :=  \dfrac{\l}{2}u^2-  \dfrac{a}{p+1}|u|^{p+1}.
\end{equation}
Its critical points are given by the zeroes of
\begin{equation}
\varphi'(u)= \lambda u-a |u|^{p-1}u,
\end{equation}
whose nature changes according to the sign of $\lambda$.

\begin{figure}[h!]
	\centering
		\begin{overpic}[scale=0.75]{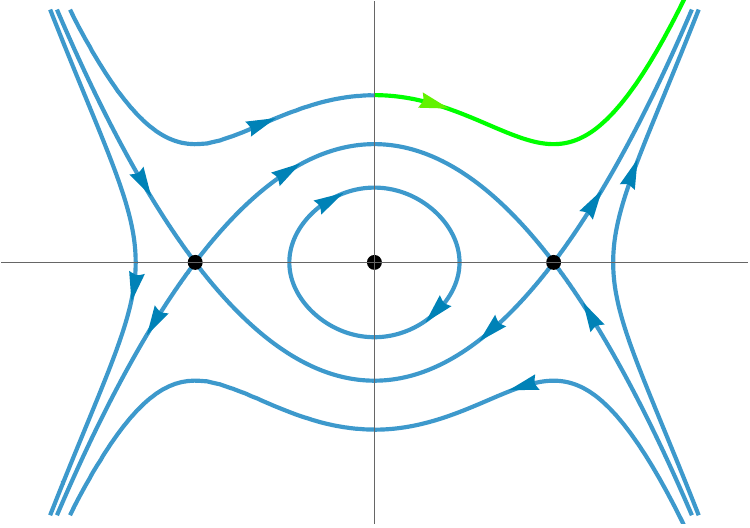}
		\put (50.5,59) {\tiny$v_0$}
		\put (50.5,52) {\tiny$v_0^*$}
		\put (72,30) {\tiny$u_0^*$}
		\put (23.5,30) {\tiny$- u_0^*$}
		\put (102,35) {\scriptsize$u$}
		\put (49,72) {\scriptsize$v$}	
	\end{overpic}
	\caption{Phase diagram of \eqref{eq:3.1} for $\l>0$. The green trajectory corresponds to a positive solution of \eqref{eq:3.1} satisfying $u(0)=0$.}
	\label{fig1}
\end{figure}
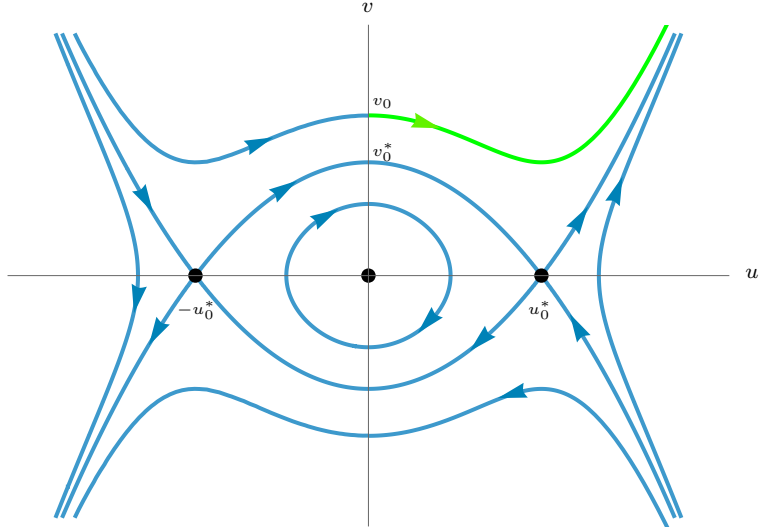

\medskip
\noindent \emph{Case $\l > 0$.} Then, $\v$ has three critical points. Namely, $u=0$,
which is a local minimum and gives rise to a center in the phase plane, and $u=\pm u_0^*$, where
\begin{equation}
u_0^* := \left(\frac{\l}{a}\right)^{\frac{1}{p-1}},
\end{equation}
that are local maxima and correspond to saddle points in the phase portrait of \eqref{eq:3.1}, which has been
sketched in Figure \ref{fig1}. To obtain positive solutions of \eqref{eq:3.2}, among all the trajectories in
Figure \ref{fig1}, we will consider those entirely lying in $\{u>0\}$ that, starting at a point $(0, v_0)$, $v_0 > 0$, are unbounded and are traveled exactly in a time equal to $R$. Necessarily, $v_0 > v_0^*$, where $v_0^*>0$ is the \lq\lq shooting critical speed\rq\rq\, of the heteroclinic connection between the two equilibria $\pm u_0^*$, which can be calculated from $E(0, v_0^*) = E(u_0^*,0)$ and  is given by
\begin{equation}
v_0^* =\sqrt{\l\frac{p-1}{p+1}}\,u_0^*.
\end{equation}

The next step is to compute the blow up time, $\mathcal{T}_{\mathcal{D}} \equiv \mc{T}_\mc{D}(v_0)$, of the solution of
\begin{equation}
	\label{eq:3.4}
	\left \{ \begin{array}{l}
		-u''=\lambda u-a |u|^{p-1}u, \\[1ex]
		u(0)=0,\quad           u'(0) = v_0,          \end{array} \right.
\end{equation}
with $v_0>v_0^*$. According to the Cauchy--Lipschitz theory (see Chapter 6 of \cite{LGT}, if necessary), the Cauchy problem \eqref{eq:3.4} has a unique (maximal) solution $u$, defined in some (maximal) interval
$I=[0,T_\mathrm{max})$ for some $T_\mathrm{max}\leq +\infty$. Moreover, if $T_\mathrm{max}<+\infty$, then
$$
  \lim_{t\ua T_\mathrm{max}}u(t)=+\infty.
$$
Since the energy defined in \eqref{eq:3.3} is conserved along trajectories, such a solution satisfies
\begin{equation}
\dfrac{1}{2}v^2(t)+ \dfrac{\l}{2}u^2(t)-  \dfrac{a}{p+1}u^{p+1}(t)= \dfrac{v_0^2}{2}
\end{equation}
for all $t\in I=[0,T_\mathrm{max})$. From this identity, it is easily seen that
$v(t)$ also blows up at $T_\mathrm{max}$ if $T_\mathrm{max}<+\infty$. Moreover, by simply having a glance at
Figure \ref{fig1}, it is easily realized that $u'(\tau)=v(\tau)>0$ for all $\tau\in I=[0,T_\mathrm{max})$. Thus, for every
$t\in (0,T_\mathrm{max})$, we find that
\begin{equation}
t = \int_0^t  d\tau  = \int_0^{t} \dfrac{u'(\tau)}{v(\tau)}\,d\tau = \int_0^{t} \dfrac{u'(\tau)}{\sqrt{v_0^2- \l u^2(\tau)+ \dfrac{2a}{p+1}u^{p+1}(\tau)}}\,d\tau
\end{equation}
and, as $u(\tau)$ is increasing, the change of variable $u(\tau) = s$ gives
\begin{equation}
\label{eq:3.5}
t =\int_0^{u(t)} \dfrac{ds}{\sqrt{v_0^2- \l s^2+ \dfrac{2a}{p+1}s^{p+1}}}.
\end{equation}
Consequently, for every $t\in (0,T_\mathrm{max})$,
$$
  t < \int_0^\infty \dfrac{ds}{\sqrt{v_0^2- \l s^2+ \dfrac{2a}{p+1}s^{p+1}}}= \imath <+\infty,
$$
because
$$
  1- \frac{p+1}{2}=-\frac{p-1}{2}<0.
$$
Note that, by the construction itself, the denominator of the previous integrand reaches its minimum at $s = u_0^*$, which is positive because we are taking $v_0 > v_0^*$. Hence, the integral on the right hand side of \eqref{eq:3.5} is convergent, which entails $T_\mathrm{max}\leq \imath$. Therefore, $u(t)$  blows up at $T_\mathrm{max}$ and, actually, letting
$t\ua T_\mathrm{max}$ in \eqref{eq:3.5}, shows that
\begin{equation}
\label{3.6}
 \mc{T}_\mc{D}(v_0) \equiv T_\mathrm{max}(v_0) =\int_0^{+\infty} \dfrac{ds}{\sqrt{v_0^2- \l s^2+ \dfrac{2a}{p+1}s^{p+1}}}= \imath.
\end{equation}
From this identity, it is easily seen that $\mathcal{T}_{\mathcal{D}}(v_0)$ is continuous and decreasing with respect to $v_0>v_0^*$. Moreover,
\begin{equation}
	\label{3.7}
	\lim_{v_0 \uparrow +\infty} \mathcal{T}_{\mathcal{D}}(v_0) = 0.
\end{equation}
Furthermore,  since the piece of heteroclinic trajectory starting from $(0,v_0^*)$ is traveled in infinite time (see Figure \ref{fig1}), the continuous dependence of the solution of \eqref{eq:3.4} on $v_0$, established by
\cite[Th. 5.26]{LGT} shows that
\begin{equation}
\lim_{v_0 \downarrow v_0^*} \mathcal{T}_{\mathcal{D}}(v_0) = + \infty.
\end{equation}

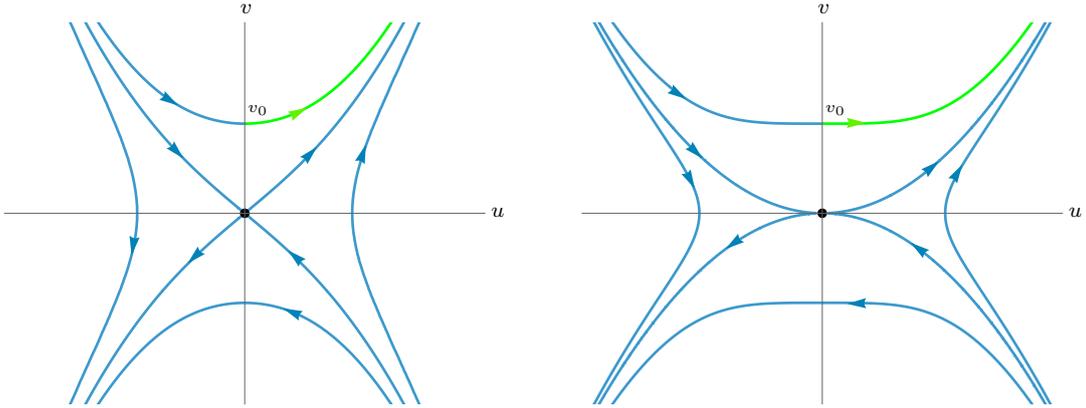
\begin{figure}[h!]
\centering
	\begin{tabular}{cc}
	\begin{overpic}[scale=0.5]{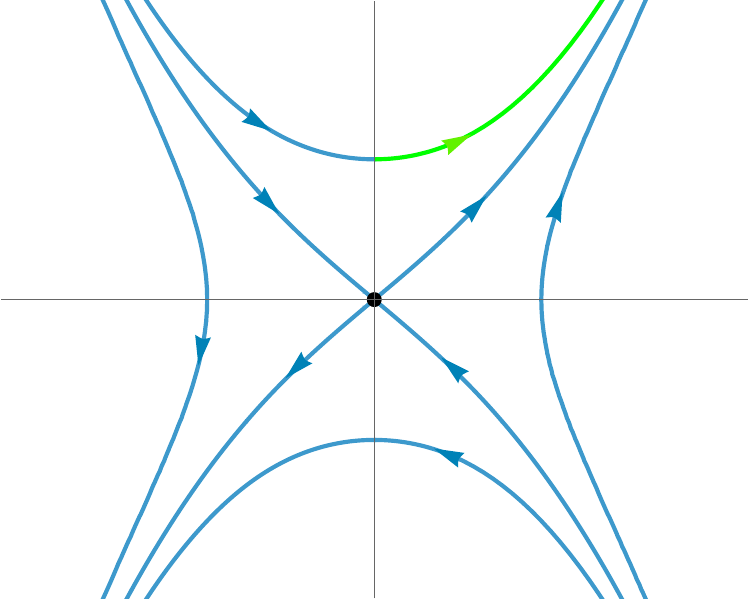}
		\put (50.5,61.5) {\tiny$v_0$}
		\put (101,39.5) {\scriptsize$u$}
		\put (49,82) {\scriptsize$v$}
	\end{overpic} &
	\hspace{7mm}
	\begin{overpic}[scale=0.5]{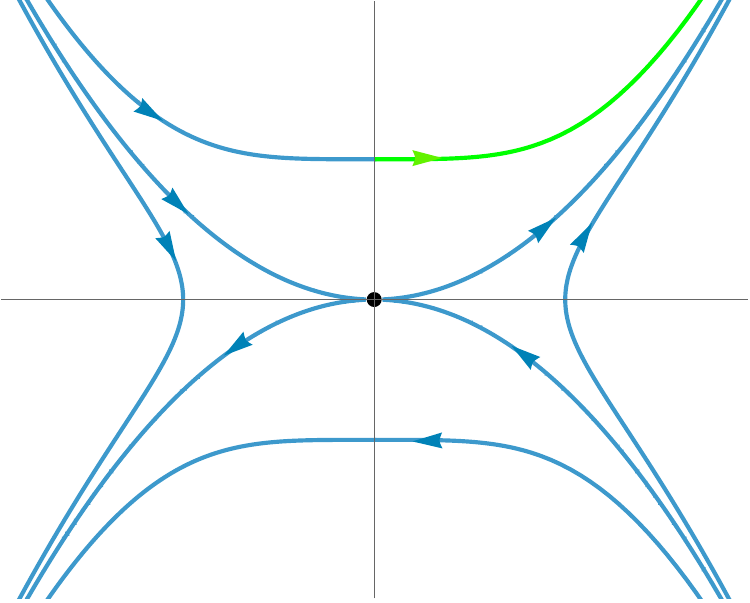}
		\put (50.5,61.5) {\tiny$v_0$}
		\put (101,39.5) {\scriptsize$u$}
		\put (49,82) {\scriptsize$v$}
	\end{overpic}
\end{tabular}
\caption{Phase diagrams of \eqref{eq:3.1} for $\l <0$ (left) and $\l=0$ (right). The green trajectory corresponds to a positive solution of \eqref{eq:3.1} satisfying $u(0)=0$.}
\label{fig2}
\end{figure}

\noindent \emph{Case $\l \leq 0$}. In this situation, $u=0$ is the unique critical point of $\varphi$, and
it is a maximum. Thus, the corresponding phase diagram looks like shown by Figure \ref{fig2},
for $\l<0$ (left) and $\l=0$ (right).

\begin{figure}[h!]
	\centering
		\begin{overpic}[scale=0.53]{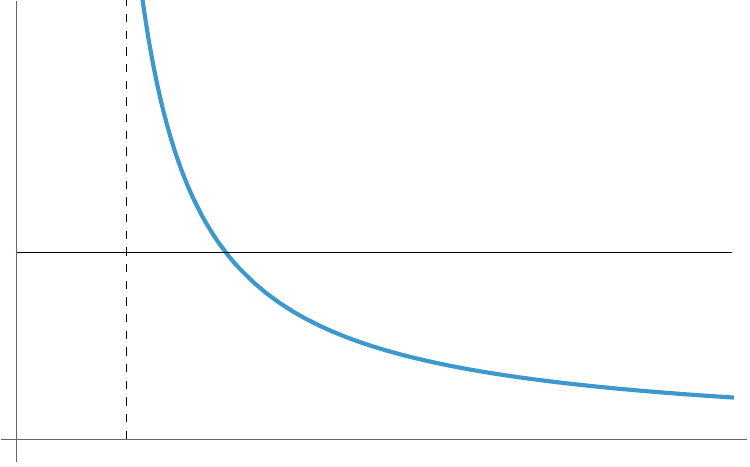}
		\put(-1.2,-0.7) {\small$0$}
		\put(102,3) {\small$v_0$}
		\put(-15.2,60) {\small$\mathcal{T}_{\mathcal{D}}(v_0)$}
		\put(100,28) {\small$R$}
		\put(15.5,-1) {\small$v_0^*$}
	\end{overpic}\hspace{1cm}
		\begin{overpic}[scale=0.53]{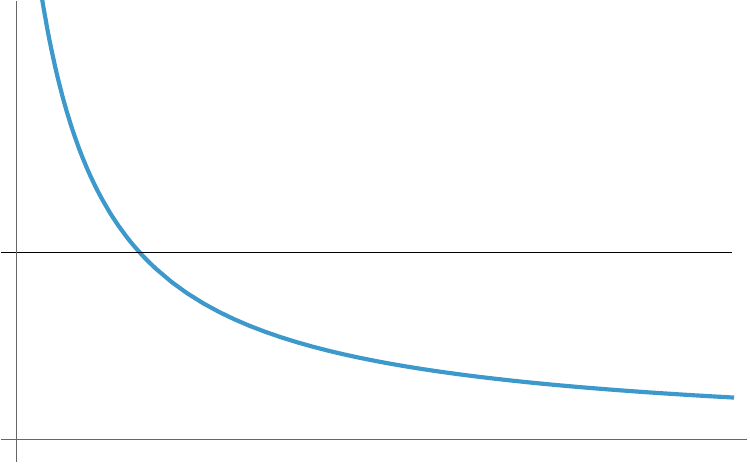}
			\put(-1.2,-0.7) {\small$0$}
			\put(102,3) {\small$v_0$}
			\put(-15.2,60) {\small$\mathcal{T}_{\mathcal{D}}(v_0)$}
			\put(100,28) {\small$R$}
\end{overpic}
	\caption{Graph of the blow-up time $\mathcal{T}_{\mathcal{D}}(v_0)$ defined in \eqref{3.6} for $\l > 0$ (left) and $\l\leq 0$ (right).}
	\label{fig3}
\end{figure}
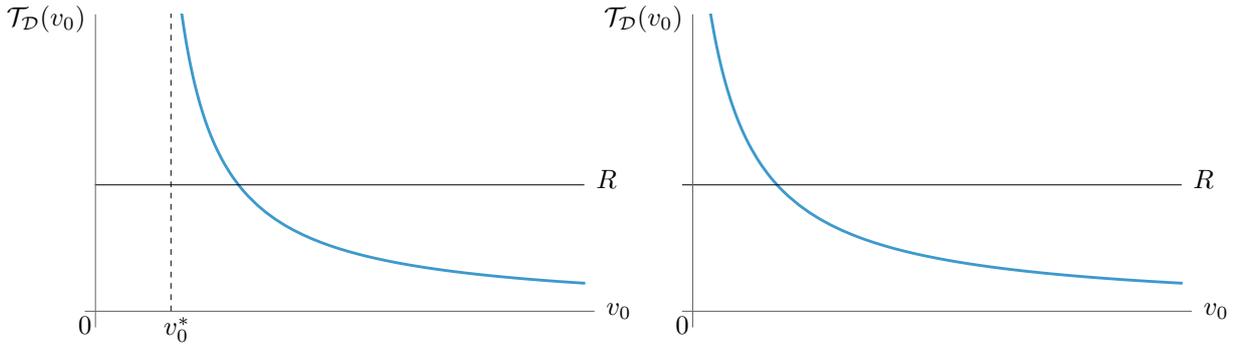

As in the previous case, we are interested in the trajectories lying in $\{u>0\}$ such that, starting at a point $(0, v_0)$, with $v_0 > 0$, are fully traveled in time $R$. By adapting the argument given in the case $\l>0$, it is easily seen that, for every $v_0>0$, the solution of the Cauchy problem \eqref{eq:3.4} also blows up after a finite time
given by \eqref{3.6}. As in the previous case, $\mc{T}_\mc{D}(v_0)$ is continuous and decreasing with respect to $v_0>0$, and \eqref{3.7} holds. Moreover, also by continuous dependence,
\begin{equation}
\lim_{v_0 \downarrow 0} \mathcal{T}_{\mathcal{D}}(v_0) = + \infty.
\end{equation}
Figure \ref{fig3} shows the plots of the functions $\mathcal{T}_{\mathcal{D}}(v_0)$ in both cases
$\l>0$ (left) and $\l \leq 0$ (right).
Thanks to these properties, it is easily realized that, for every $R > 0$, there exists a unique initial value $v_0 =v_0(R) > 0$ such that $\mathcal{T}_{\mathcal{D}}(v_0) = R$. Moreover, $v_0(R)>v_0^*$ if $\l>0$. Therefore, for every $R > 0$, the problem \eqref{eq:3.2} has a unique positive solution. This concludes the proof of Part (a) for Dirichlet boundary conditions.

\subsection{The Neumann case $\mathcal{B} = \mathcal{N}$}
\label{sec:3.2}

In this section, we consider problem \eqref{eq:1.1} with $\mathcal{B} = \mathcal{N}$, which reads
\begin{equation}
	\label{3.8}
	\left \{ \begin{array}{l}
		-u''=\lambda u-a |u|^{p-1}u\qquad \hbox{in} \;\; [0,R),\\[1ex]
		u'(0)=0,\quad    u(R) = +\infty.          \end{array} \right.
\end{equation}
As in Section \ref{sec:3.1}, we are differentiating the cases $\l >0$ and $\l\leq 0$.

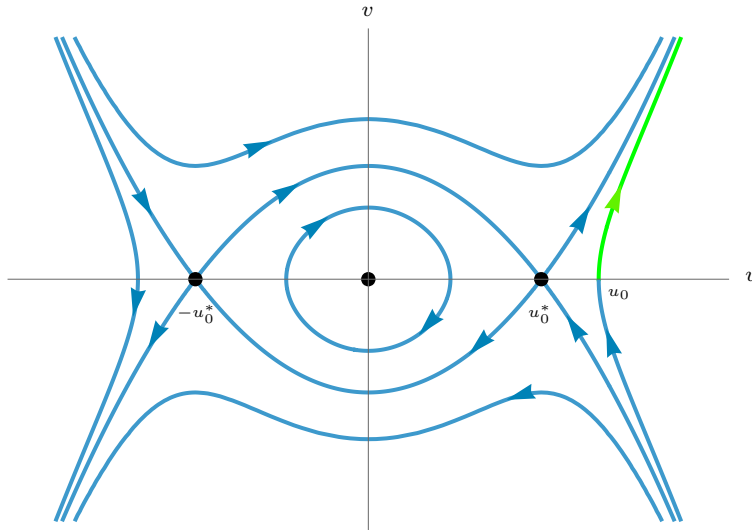
\begin{figure}[h!]
	\centering
		\begin{overpic}[scale=0.75]{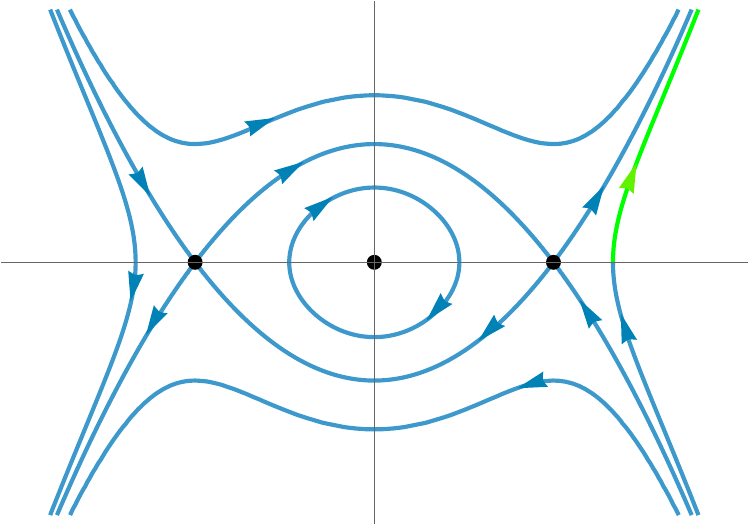}
		\put (83,33) {\tiny$u_0$}
		\put (72,30) {\tiny$u_0^*$}
		\put (23.5,30) {\tiny$- u_0^*$}
		\put (102,35) {\scriptsize$u$}
		\put (49,72) {\scriptsize$v$}
	\end{overpic}
	\caption{Phase diagram of the equation \eqref{eq:3.1} for $\l>0$. The green trajectory corresponds to a positive solution satisfying $u'(0)=v(0)=0$.}
	\label{fig4}
\end{figure}

\medskip
\noindent \emph{Case $\l >0$.} Although the phase portrait of the differential equation has been already sketched in
Figure \ref{fig1}, in Figure \ref{fig4} we are highlighting the positive trajectories  lying in $\{u>0\}$ such that, starting at a point $(u_0,0)$ with $u_0>0$, are unbounded and completely traveled in time $R >0$. A necessary condition for this is $u_0>u_0^*$.

According to the Cauchy--Lipschitz theory, for every $u_0>u_0^*$, the underlying Cauchy problem
\begin{equation}
\label{3.9}
	\left \{ \begin{array}{l}
		-u''=\lambda u-a |u|^{p-1}u,\\[1ex]
		u(0)=u_0,\quad        u'(0) = 0,          \end{array} \right.
\end{equation}
has a unique (maximal) solution $u$, defined in some (maximal) interval
$I=[0,T_\mathrm{max})$ for some $T_\mathrm{max}\leq +\infty$. Moreover, if $T_\mathrm{max}<+\infty$, then
$$
  \lim_{t\ua T_\mathrm{max}}u(t)=+\infty.
$$
Since the energy defined in \eqref{eq:3.3} is conserved along trajectories, such a solution satisfies
\begin{equation}
\frac{1}{2}v^2(t)+ \frac{\l}{2}u^2(t)-  \frac{a}{p+1}u^{p+1}(t)= \dfrac{\l}{2}u_0^2-  \dfrac{a}{p+1}u_0^{p+1}
\end{equation}
for all $t \in I=[0,T_\mathrm{max})$. Thus, for every $t\in (0,T_\mathrm{max})$, since $v(\tau)=u'(\tau)>0$ for all
$\tau \in (0,T_\mathrm{max})$,
\begin{equation}
t = \int_0^t d\tau  = \int_0^{t} \dfrac{u'(\tau)}{v(\tau)}\,d\tau = \int_0^t
 \dfrac{u'(\tau)}{\sqrt{\dfrac{2a}{p+1}\left( u^{p+1}(\tau)- u_0^{p+1} \right)- \l \left( u^2(\tau)-u_0^2 \right)} }\,d\tau.
\end{equation}
Thus, doing the change of variable $s=u(\tau)$, we find that
\begin{equation}
t= \int_{u_0}^{u(t)}  \dfrac{ds}{\sqrt{\dfrac{2a}{p+1}\left( s^{p+1}-u_0^{p+1}\right) - \l \left(s^2-u_0^2\right)}}.
\end{equation}
Hence, setting $s:=u_0 \t$, we are lead to
\begin{equation}
\label{3.10}
t = \int_1^{u(t)/u_0}  \dfrac{d \theta}{\sqrt{ \dfrac{2a}{p+1}u_0^{p-1}\left( \theta^{p+1}-1\right)-\l
(\t^2-1)}}\;\; \hbox{for all}\;\; t \in I.
\end{equation}
Consequently, since
$$
    \jmath := \int_1^{+\infty}  \dfrac{d \theta}{\sqrt{ \dfrac{2a}{p+1}u_0^{p-1}\left( \theta^{p+1}-1\right)-\l
(\t^2-1)}}<+\infty,
$$
because $u_0>u_0^*$, it becomes apparent that $T_\mathrm{max}\leq \jmath$ and, actually, letting $t\ua T_\mathrm{max}$ in \eqref{3.10},
\begin{equation}
\label{3.11}
  \mc{T}_\mc{N}(u_0)\equiv T_\mathrm{max}(u_0)=\int_1^{+\infty}  \dfrac{d \theta}{\sqrt{ \dfrac{2a}{p+1}u_0^{p-1}\left( \theta^{p+1}-1\right)-\l
(\t^2-1)}}=\jmath <+\infty.
\end{equation}
Therefore, $\mathcal{T}_{\mathcal{N}}(u_0)$ is continuous and decreasing with respect to $u_0>u_0^*$. Moreover, reasoning as in the Dirichlet case treated in Section \ref{sec:3.1}, it readily follows that
\begin{equation}
\label{3.12}
\lim_{u_0 \downarrow u_0^*} \mathcal{T}_{\mathcal{N}}(u_0) = + \infty, \qquad \lim_{u_0 \uparrow +\infty} \mathcal{T}_{\mathcal{N}}(u_0) = 0.
\end{equation}
\par

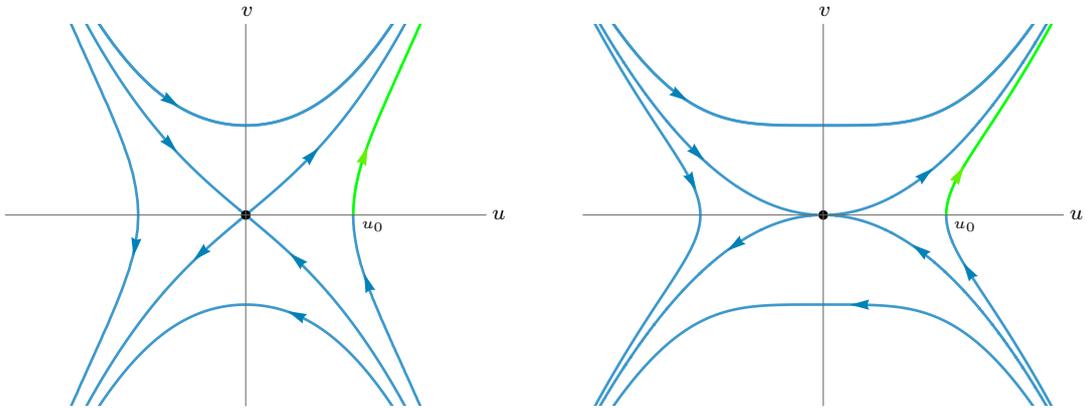
\begin{figure}[h!]
\centering
\begin{tabular}{cc}
	\begin{overpic}[scale=0.5]{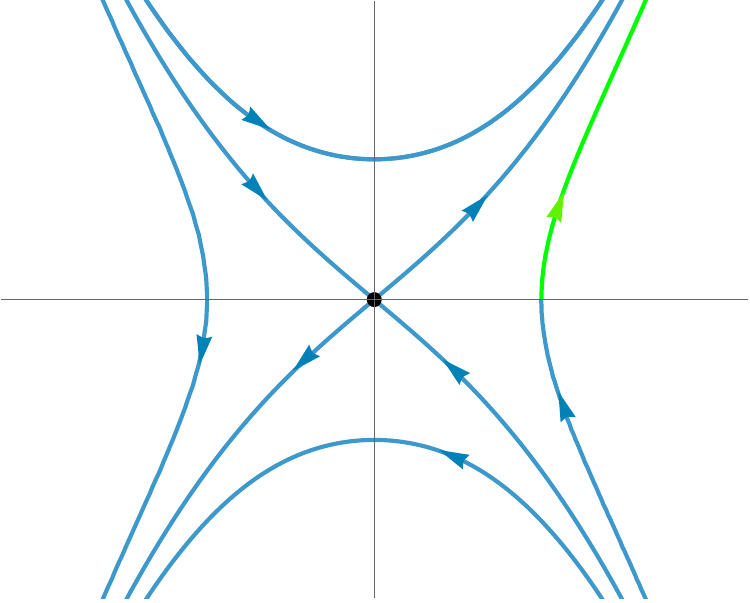}
		\put (74,37.5) {\tiny$u_0$}
		\put (101,39.5) {\scriptsize$u$}
		\put (49,82) {\scriptsize$v$}
	\end{overpic} &
	\hspace{7mm}
	\begin{overpic}[scale=0.5]{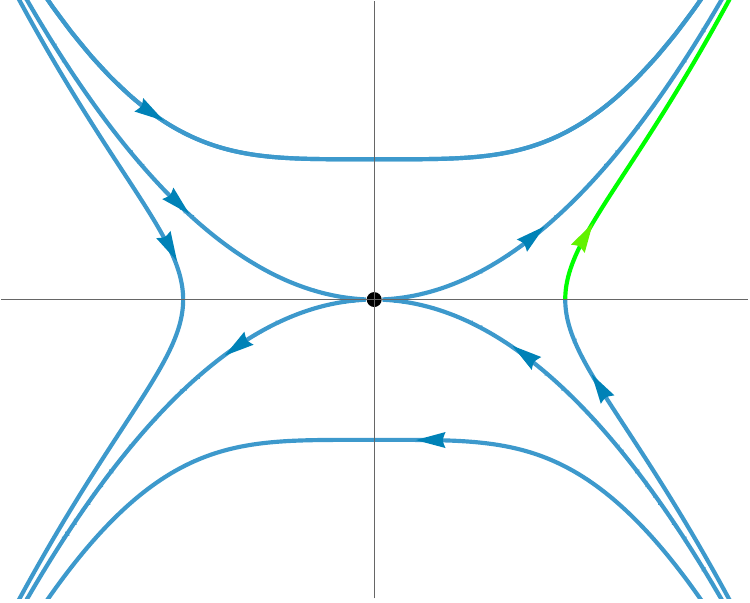}
		\put (77,37.5) {\tiny$u_0$}
		\put (101,39.5) {\scriptsize$u$}
		\put (49,82) {\scriptsize$v$}
	\end{overpic}\\
\end{tabular}
\caption{Phase diagram of the equation \eqref{eq:3.1} for $\l <0$ (left) and $\l=0$ (right). The green trajectory corresponds to a positive solution satisfying $u'(0)=0$.}
\label{fig5}
\end{figure}

\begin{figure}[h!]
	\centering
	\begin{overpic}[scale=0.53]{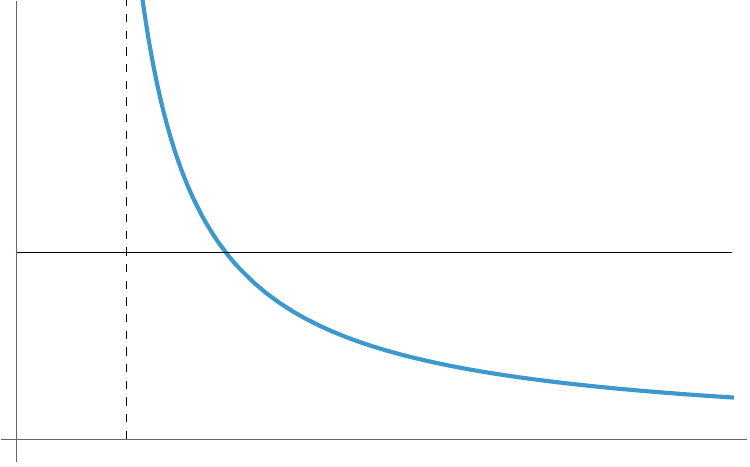}
		\put(-1.2,-0.7) {\small$0$}
	\put(102,3) {\small$u_0$}
	\put(-15.2,60) {\small$\mathcal{T}_{\mathcal{N}}(u_0)$}
	\put(100,28) {\small$R$}
	\put(15.5,-1) {\small$u_0^*$}
\end{overpic} \hspace{1cm}
\begin{overpic}[scale=0.53]{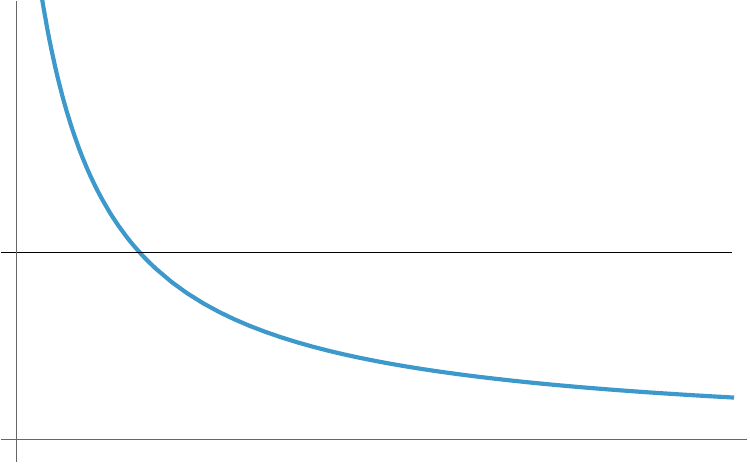}
	\put(-1.2,-0.7) {\small$0$}
	\put(102,3) {\small$u_0$}
	\put(-15.2,60) {\small$\mathcal{T}_{\mathcal{N}}(u_0)$}
	\put(100,28) {\small$R$}
\end{overpic}
	\caption{Graph of the blow-up time $\mathcal{T}_{\mathcal{N}}(u_0)$ defined in \eqref{3.11}, for  $\l > 0$ (left) and $\l\leq 0$ (right).}
	\label{fig6}
\end{figure}

\noindent \emph{Case $\l \leq 0$.} Figure \ref{fig5} represents the phase diagram of equation \eqref{eq:3.1} in this case. As in the previous situation, we are interested in the trajectories that, starting at a point $(u_0,0)$, with $u_0>0$, are completely traveled in a finite time $R >0$.
Reasoning as in the case $\l>0$, it turns out that, for every $u_0>0$, the solution of \eqref{3.9} in the case $\l\leq 0$ also blows up at the value of $\mc{T}_\mc{N}$ given by \eqref{3.11}. Thus, $\mathcal{T}_{\mathcal{N}}(u_0)$ is continuous, decreasing with respect to $u_0>0$, and it satisfies
\begin{equation}
	\label{3.13}
\lim_{u_0 \downarrow 0} \mathcal{T}_{\mathcal{N}}(u_0) = + \infty, \qquad
\lim_{u_0 \uparrow +\infty} \mathcal{T}_{\mathcal{N}}(u_0) = 0.
\end{equation}
Figure \ref{fig6} shows the graphs of $\mc{T}_\mc{N}(u_0)$  for each of these cases. Note that they satisfy \eqref{3.12} and \eqref{3.13}. According to them, for every $R > 0$,  there exists a unique $u_0 =u_0(R)> 0$ such that $\mathcal{T}_{\mathcal{N}}(u_0) = R$.
Moreover, $u_0>u_0^*$ if $\l>0$. Therefore, for every $\l\in\R$, the problem \eqref{3.8} has a unique positive solution. This ends the proof of Part (a) for Neumann boundary conditions.

\subsection{The Robin case $\mathcal{B} = \mathcal{R}_{\beta}$ with $\beta > 0$}
\label{sec:3.3}
In this section, we will focus the attention on problem \eqref{eq:1.1} choosing $\mathcal{B} = \mathcal{R}_{\beta}$ for some $\beta >0$. In such case \eqref{eq:1.1} becomes
\begin{equation}
	\label{3.14}
	\left \{ \begin{array}{ll}
		-u''=\lambda u-a |u|^{p-1}u\qquad \hbox{in} \;\; [0,R),\\[1ex]
		u'(0)=\beta u(0),\quad  u(R) = +\infty.          \end{array} \right.
\end{equation}
As in the previous sections, we will differentiate the cases $\l >0$ and $\l \leq 0$.
\par
\medskip

\begin{figure}[h!]
	\centering
		\begin{overpic}[scale=0.75]{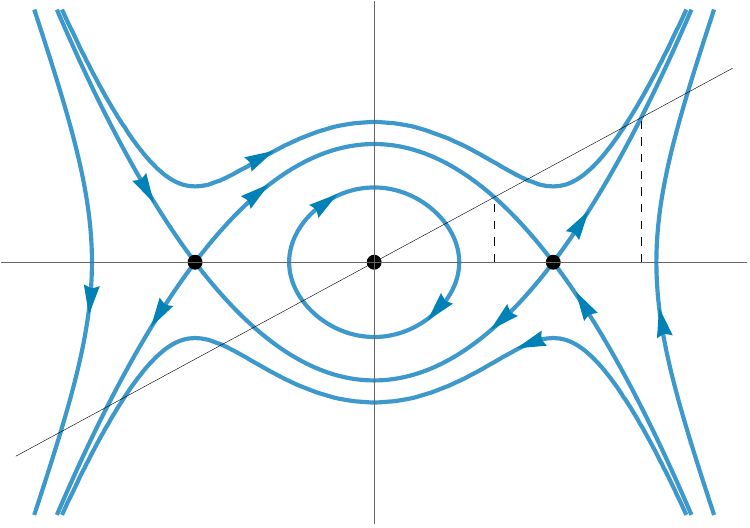}
			\put (100,61) {\scriptsize $v = \beta u$}
			\put (72,30) {\tiny$u_0^*$}
			\put (65,33) {\tiny$u_{-}$}
			\put (83.7,33) {\tiny$u_{+}$}
			\put (23.5,30) {\tiny$- u_0^*$}
			\put (102,35) {\scriptsize$u$}
			\put (49,72) {\scriptsize$v$}
	\end{overpic}
	\caption{Phase diagram of the equation \eqref{eq:3.1} for $\l>0$ taking $\mathcal{B} = \mathcal{R}_{\beta}$ with $\beta > 0$.}
	\label{fig7}
\end{figure}
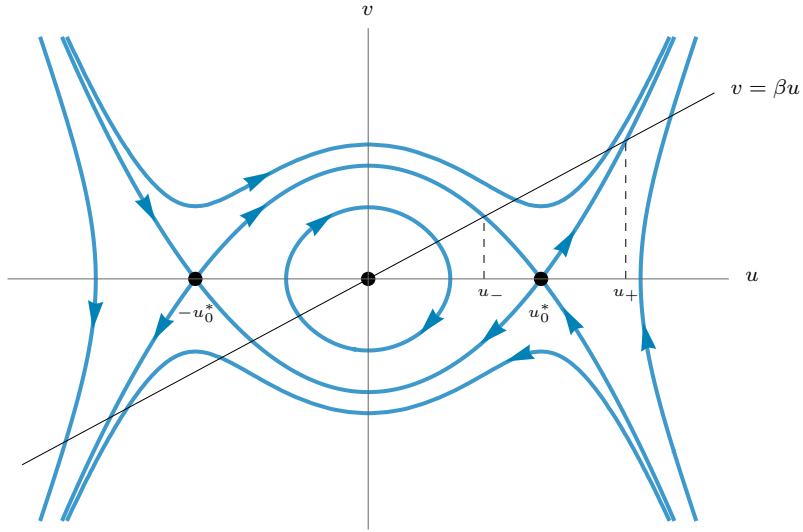

\noindent \emph{Case $\l >0$.} In Figure \ref{fig7} we are superimposing the phase diagram of \eqref{eq:3.1} in this case with the straight line $v=\b u$. Now, the line $v=\b u$ crosses, in the region $\{u>0\}$,
 the integral curve containing the heteroclinic connections  linking $(-u_0^*,0)$ to $(u_0^*,0)$ at exactly two points with abscisas $u_-$ and $u_+$ such that
$u_-<u_0^*<u_+$.

The existence and uniqueness of $u_-$ is a direct consequence of the fact that the
heteroclinic connection looks like shows the figure and $v=\b u$ is unbounded. Thus, $v=\b u$ must leave the region enclosed by the two heteroclinics. To show the existence and uniqueness of $u_+$ we should prove that there is a unique
value of $u>u_0^*$, $u_+$, such that
$$
   \frac{\beta^2+\l}{2}u^2-\frac{a}{p+1}|u|^{p+1}= \frac{\l}{2}(u_0^*)^2-\frac{a}{p+1}(u_0^*)^{p+1}\equiv
   \v_*.
$$
Indeed, setting
$$
  f(u):=  \frac{\beta^2+\l}{2}u^2-\frac{a}{p+1}|u|^{p+1}-\v_*,
$$
we have that $f(u_-)=0$, $f(u_0^*)=\frac{\b^2}{2}(u_0^*)^2>0$, and, for every $u> u_-$,
$$
  f'(u) = (\b^2+\l)u-au^p =u(\b^2+\l-au^{p-1}).
$$
Hence, $f'(u)>0$ for all $u\in (u_-,u_c)$, where we are denoting
$$
   u_c:= \left(\tfrac{\b^2+\l}{a}\right)^\frac{1}{p-1},
$$
while $f'(u)<0$ if $u> u_c$. Note that
$$
  u_0^* = \left( \tfrac{\l}{a}\right)^\frac{1}{p-1}<  \left(\tfrac{\b^2+\l}{a}\right)^\frac{1}{p-1}=u_c.
$$
In particular, $f(u)>0$ for all $u\in (u_-,u_c]$. Moreover,
$$
  \lim_{u\ua+\infty}f(u)=-\infty.
$$
Therefore, there is a unique $u_+ >u_c$ such that $f(u)>0$ for all $u\in (u_-,u_+)$, $f(u_+)=0$, and
$f(u)<0$ for all $u>u_+$, as already illustrated in Figure \ref{fig7}.
\par
\medskip

\begin{figure}[h!]
\centering
\begin{tabular}{cc}
	\begin{overpic}[scale=0.59]{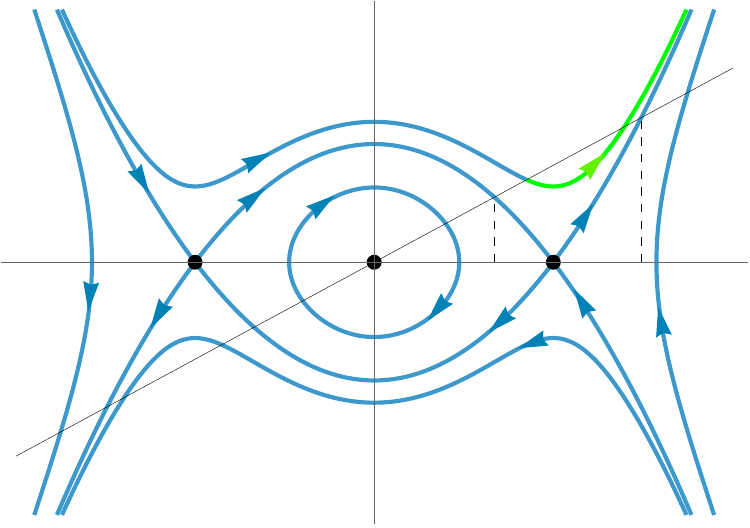}
			\put (100,61) {\scriptsize $v = \beta u$}
			\put (72,29) {\tiny$u_0^*$}
			\put (64.5,33) {\tiny$u_{-}$}
			\put (83,33) {\tiny$u_{+}$}
			\put (22.8,29) {\tiny$- u_0^*$}
			\put (102,35) {\scriptsize$u$}
			\put (49,72) {\scriptsize$v$}
	\end{overpic} &
	\hspace{7mm}
	\begin{overpic}[scale=0.59]{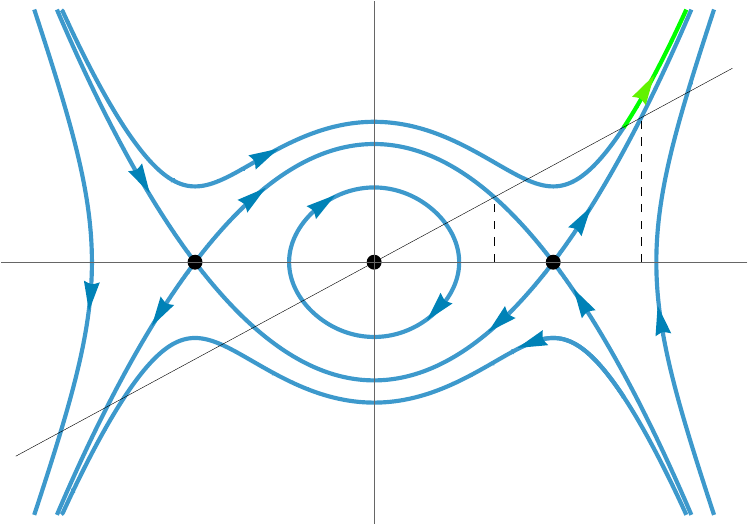}
			\put (100,61) {\scriptsize $v = \beta u$}
			\put (72,29) {\tiny$u_0^*$}
			\put (64.5,33) {\tiny$u_{-}$}
			\put (83,33) {\tiny$u_{+}$}
			\put (22.8,29) {\tiny$- u_0^*$}
			\put (102,35) {\scriptsize$u$}
			\put (49,72) {\scriptsize$v$}	
	\end{overpic} \\
	\vspace{10mm}
	(a) Case $u_0 \in (u_-,u_0^*)$. & (b) Case $u_0 \in (u_0^*,u_+)$. \\
	\begin{overpic}[scale=0.59]{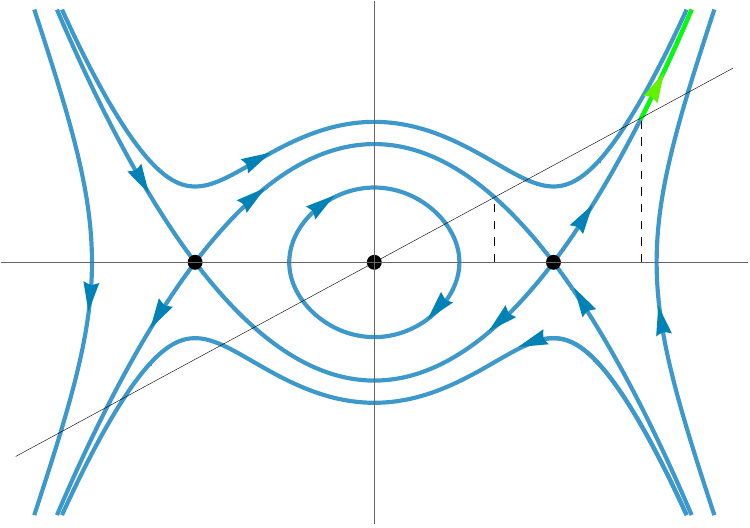}
			\put (100,61) {\scriptsize $v = \beta u$}
			\put (72,29) {\tiny$u_0^*$}
			\put (64.5,33) {\tiny$u_{-}$}
			\put (83,33) {\tiny$u_{+}$}
			\put (22.8,29) {\tiny$- u_0^*$}
			\put (102,35) {\scriptsize$u$}
			\put (49,72) {\scriptsize$v$}
	\end{overpic} &
	\hspace{4mm}
	\begin{overpic}[scale=0.59]{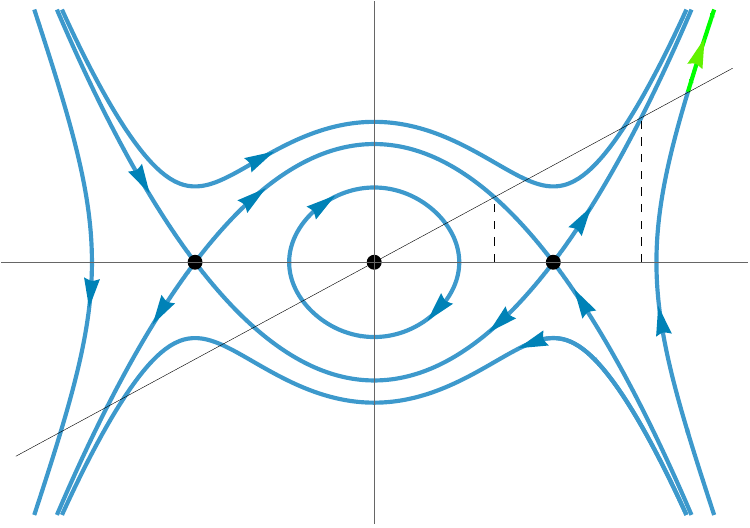}
			\put (100,61) {\scriptsize $v = \beta u$}
			\put (72,29) {\tiny$u_0^*$}
			\put (64.5,33) {\tiny$u_{-}$}
			\put (83,33) {\tiny$u_{+}$}
			\put (23,29) {\tiny$- u_0^*$}
			\put (102,35) {\scriptsize$u$}
			\put (49,72) {\scriptsize$v$}
	\end{overpic} \\
	(c) Case $u_0 = u_+$. & (d) Case $u_0 >u_+$. \\
\end{tabular}
\caption{The green trajectories on each of the phase portraits correspond to the positive solutions of \eqref{3.14} provided they are run in a time $R$.}
\label{fig8}
\end{figure}
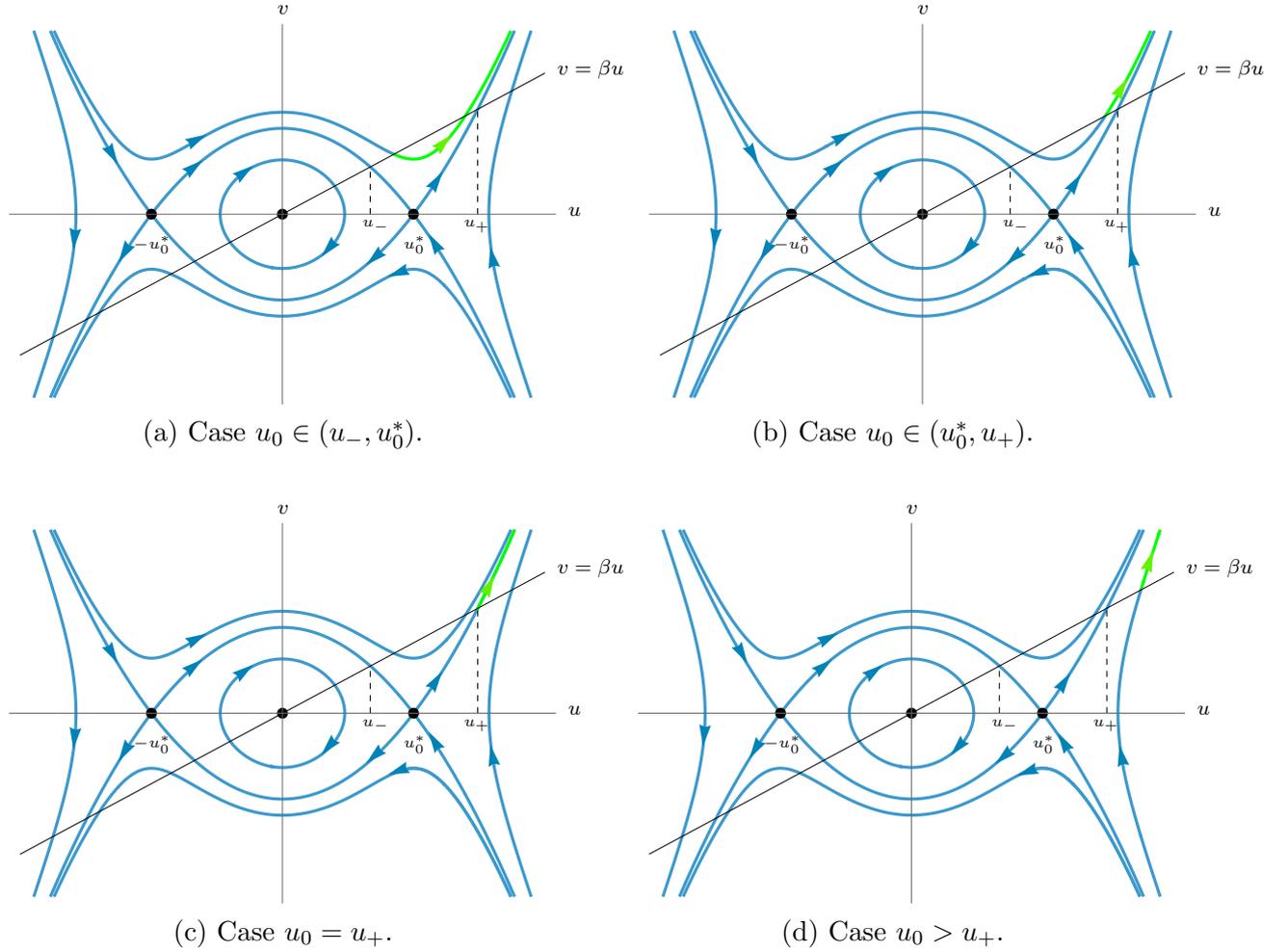

As we aim to solve problem \eqref{3.14}, among the trajectories collected in Figure \ref{fig7}, we are only interested in those lying in $\{u>0\}$, starting at a point $(u_0,v_0)$ with $v_0=\b u_0>0$, which are unbounded and traveled in a time $R$. They have been plotted in green in Figure \ref{fig8}. Note that, in any circumstances,
$u_0>u_-$.

As in the previous cases, we should find out the time, $\mathcal{T}_{\mathcal{R}_{\b}}\equiv  \mathcal{T}_{\mathcal{R}_{\b}}(u_0)$, it takes for the solution of
\begin{equation}
	\label{3.15}
	\left \{ \begin{array}{ll}
		-u''=\lambda u-a |u|^{p-1}u,\\[1ex]
		u(0)=u_0,\quad u'(0)=v_0=\beta u_0,       \end{array} \right.
\end{equation}
to blow up for every $u_0>u_-$. According to the Cauchy--Lipschitz theory, the  problem \eqref{3.15} has a unique (maximal) solution $u$ defined in some (maximal) interval $I=[0,T_\mathrm{max})$ for some $T_\mathrm{max}\leq +\infty$. Moreover, since the system is conservative, if $T_\mathrm{max}<+\infty$, then
$$
  \lim_{t\ua T_\mathrm{max}}u(t)=+\infty.
$$
As the energy defined in \eqref{eq:3.3} is preserved along trajectories, such a solution satisfies
\begin{equation}
\frac{1}{2}v^2(t)+ \frac{\l}{2}u^2(t)-  \frac{a}{p+1}u^{p+1}(t)=\left(\frac{\b^2+\l}{2}\right)u_0^2-\frac{a}{p+1}u_0^{p+1}
\end{equation}
for all $t \in I = [0,T_\mathrm{max})$. Reasoning as in the Neumann case presented in Section \ref{sec:3.2}, it is
easily seen that $T_\mathrm{max}<+\infty$ and that, actually,
\begin{equation}
\label{3.16}
\mathcal{T}_{\mathcal{R}_{\b}}(u_0) \equiv T_\mathrm{max}(u_0) =  \int_1^{+\infty}  \dfrac{d \theta}{\sqrt{
\dfrac{2a}{p+1}u_0^{p-1}\left(\theta^{p+1}-1\right) - \l \left(\theta^2-1\right)  +\beta^2}} < + \infty
\end{equation}
for all $u_0>u_-$. Thus, $\mathcal{T}_{\mathcal{R}_{\b}}(u_0)$ is continuous and decreasing with respect to $u_0>u_{-}$. Moreover, adapting the argument given in the Neumann case yields
\begin{equation}
\label{3.17}
\lim_{u_0 \downarrow u_{-}} \mathcal{T}_{\mathcal{R}_{\b}}(u_0) = + \infty, \qquad
\lim_{u_0 \uparrow +\infty} \mathcal{T}_{\mathcal{R}_{\b}}(u_0) = 0.
\end{equation}
\par
\begin{figure}[h!]
\centering
	\begin{tabular}{cc}
	\begin{overpic}[scale=0.59]{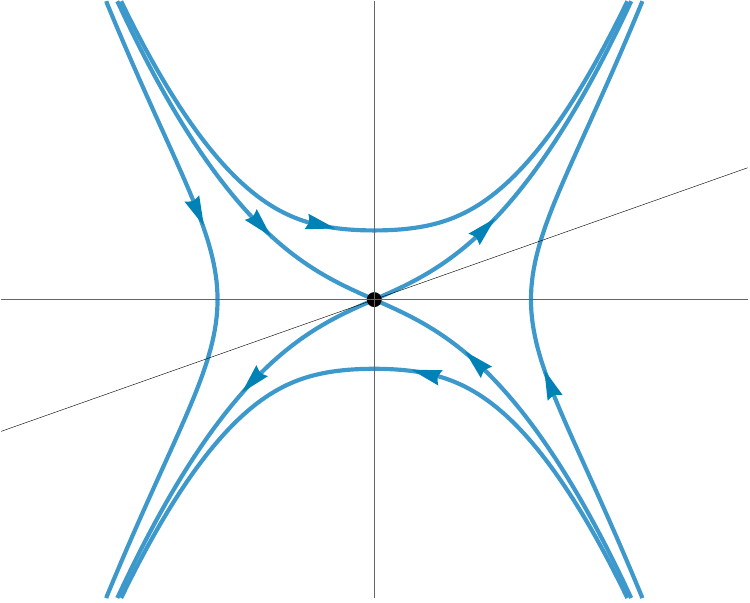}
		\put (101,39.5) {\scriptsize$u$}
		\put (49,82) {\scriptsize$v$}
		\put (102,58) {\scriptsize $v = \beta u$}
	\end{overpic} &
	\hspace{7mm}
	\vspace{3mm}
	\begin{overpic}[scale=0.59]{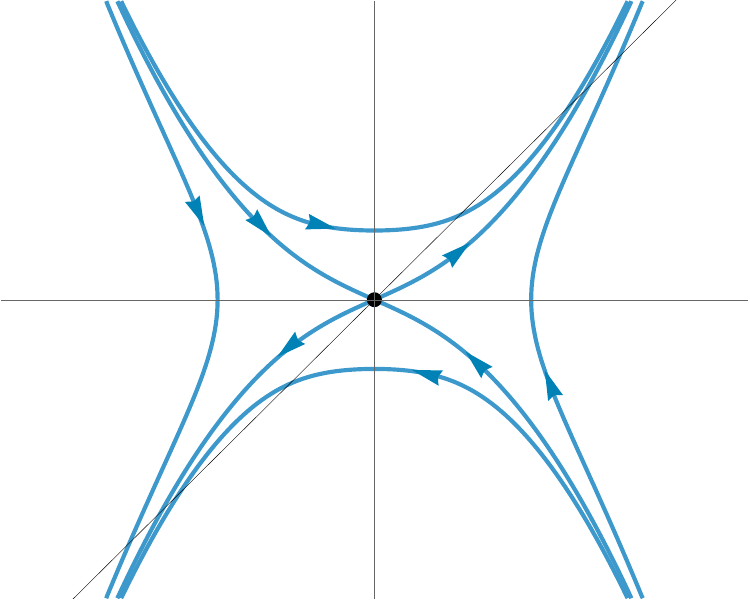}
		\put (101,39.5) {\scriptsize$u$}
		\put (49,82) {\scriptsize$v$}
		\put (93,80) {\scriptsize $v = \beta u$}
	\end{overpic}\\
	(a) Case $\b^2 \in (0,-\l]$. & (b) Case $\b^2 >- \l$. \\
\end{tabular}
\caption{Phase diagram of the equation \eqref{eq:3.1} for $\l<0$ taking $\mathcal{B} = \mathcal{R}_{\beta}$ with $\beta > 0$.}
\label{fig9}
\end{figure}
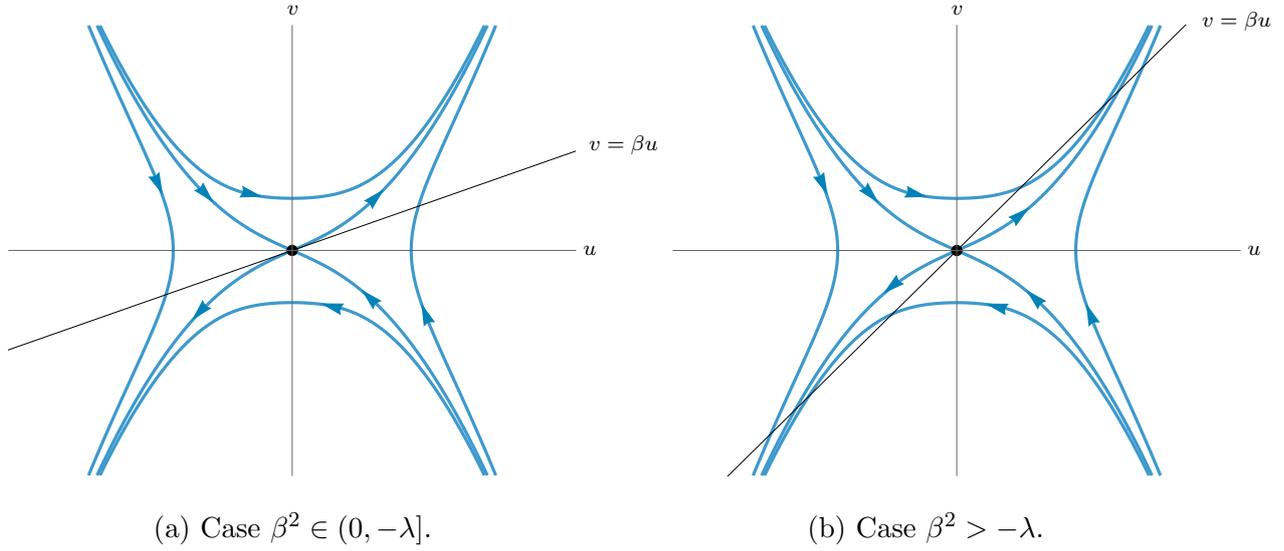

\medskip
\noindent \emph{Case $\l \leq 0$.} Figure \ref{fig9} sketches the phase diagram of equation \eqref{eq:3.1} in this case.
Now, it consists of a unique equilibrium, $(0,0)$, which  is a \lq global\rq \, saddle point,
whose branch of unstable manifold within the first quadrant is given by
$$
  v=\sqrt{\frac{2a}{p+1}u^{p+1}-\l u^2},\qquad u>0.
$$

Thus, it has a unique crossing point  with the straight line $v=\b u$ if and only if $\b^2>-\l$, which is the case represented in Figure \ref{fig9}(b). Figure \ref{fig9}(a) corresponds with a situation when $\b^2\in (0,-\l]$, where the branch of the unstable manifold in the first quadrant does not meet $v=\b u$, except at the origin.
\par
\begin{figure}[h!]
	\centering
		\begin{overpic}[scale=0.65]{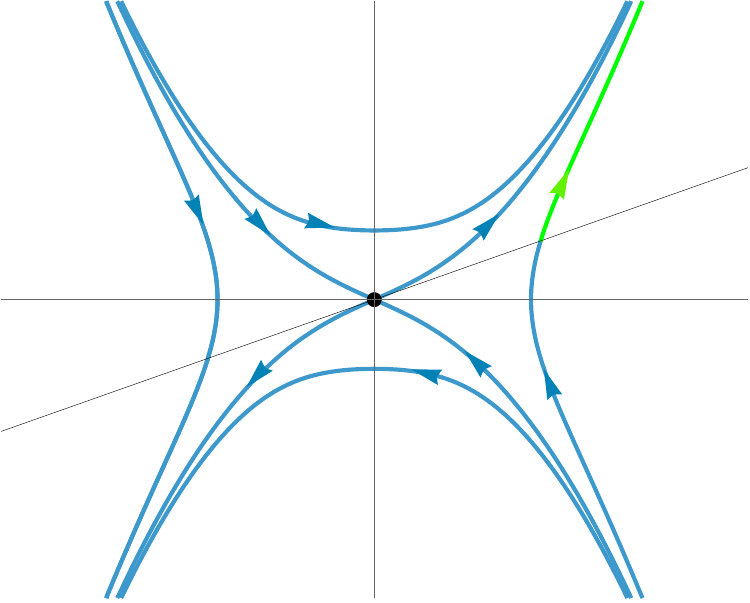}
		\put (101,39.5) {\scriptsize$u$}
		\put (49,82) {\scriptsize$v$}
		\put (102,58) {\scriptsize $v = \beta u$}
	\end{overpic}
	\caption{The green trajectory corresponds to a positive solution of \eqref{3.14} in the case $\b^2 \in (0,- \l]$.}
	\label{fig10}
\end{figure}
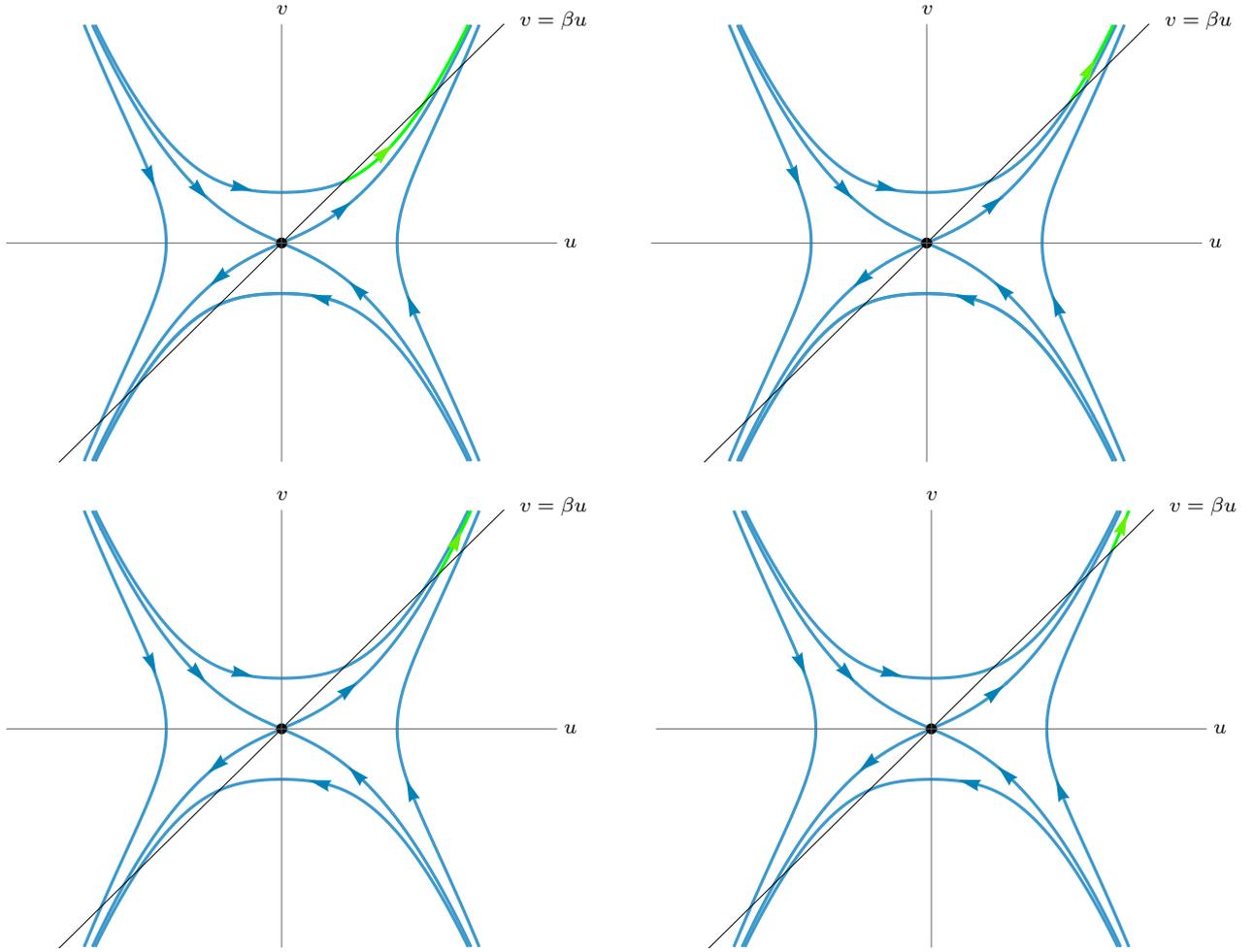
\begin{figure}[h!]
\centering
\begin{tabular}{cc}
	\begin{overpic}[scale=0.59]{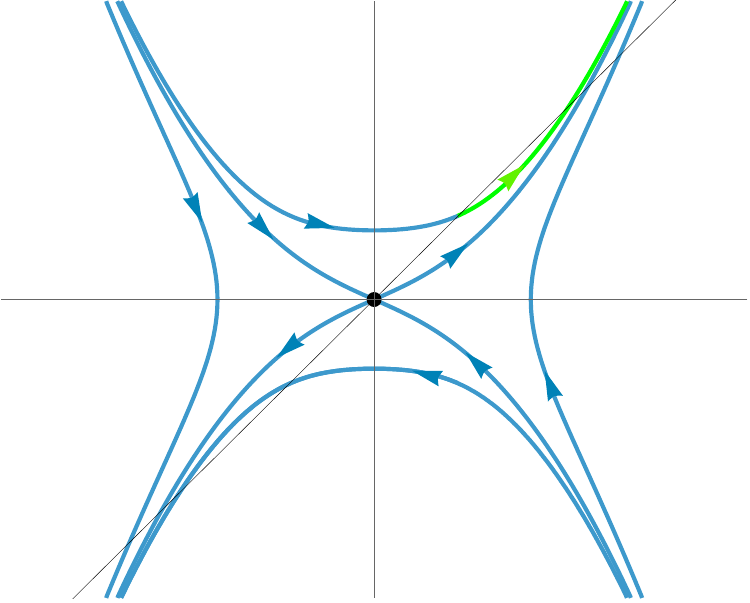}
		\put (101,39.5) {\scriptsize$u$}
		\put (49,82) {\scriptsize$v$}
		\put (93,80) {\scriptsize $v = \beta u$}
	\end{overpic} &
	\hspace{7mm}
	\begin{overpic}[scale=0.59]{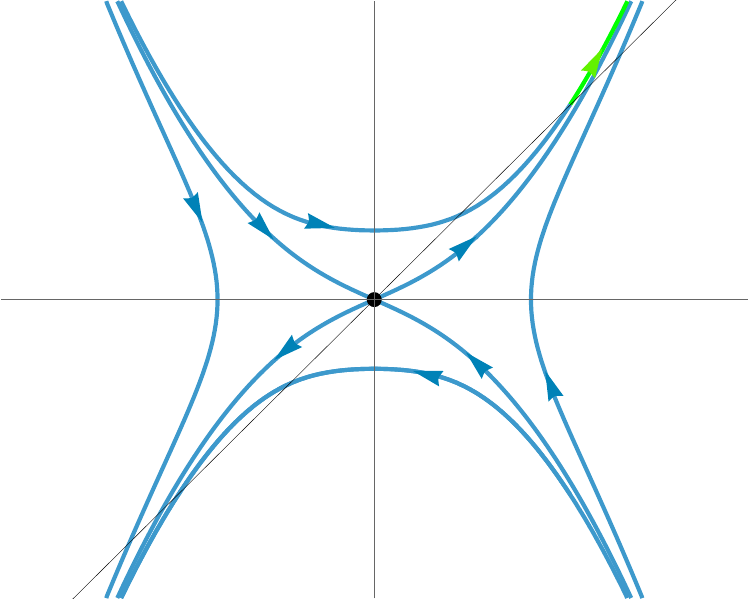}
		\put (101,39.5) {\scriptsize$u$}
		\put (49,82) {\scriptsize$v$}
		\put (93,80) {\scriptsize $v = \beta u$}
	\end{overpic} 	\vspace{5mm} \\
	\begin{overpic}[scale=0.59]{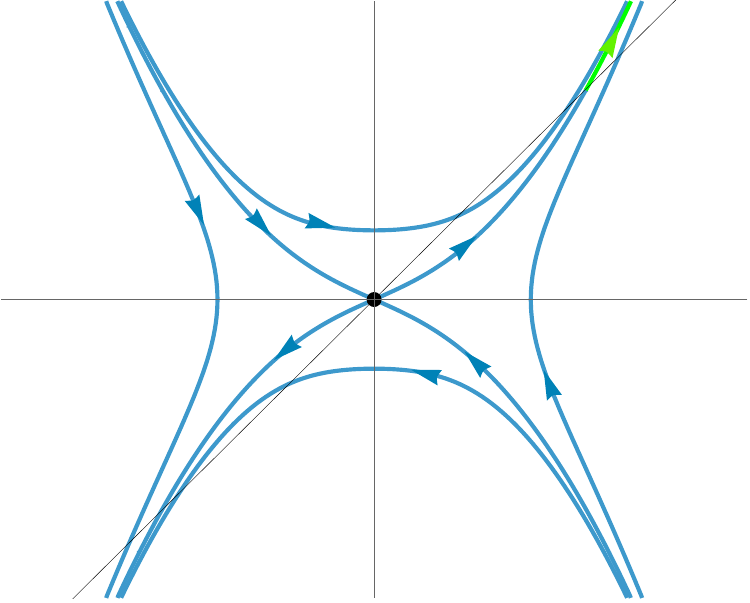}
		\put (101,39.5) {\scriptsize$u$}
		\put (49,82) {\scriptsize$v$}
		\put (93,80) {\scriptsize $v = \beta u$}
	\end{overpic} &
	\hspace{7mm}
	\begin{overpic}[scale=0.59]{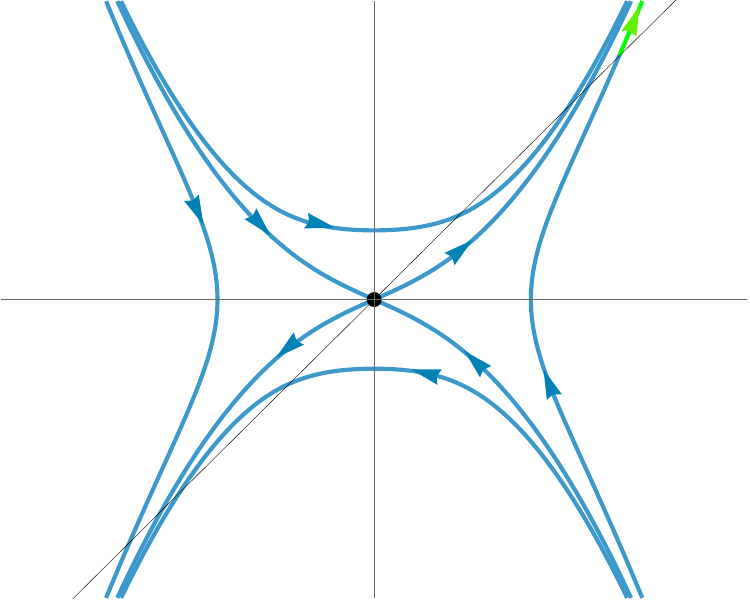}
		\put (101,39.5) {\scriptsize$u$}
		\put (49,82) {\scriptsize$v$}
		\put (93,80) {\scriptsize $v = \beta u$}
	\end{overpic} \\
\end{tabular}
\caption{The green trajectories on each of the phase portraits correspond to the positive solutions of \eqref{3.14} in the case $\b^2 > -\l$.}
\label{fig11}
\end{figure}

As above, we are interested in the trajectories that, starting at a point $(u_0,v_0) = (u_0, \beta u_0)$, with $u_0 >0$, blow up in a finite time $R >0$. In the case $\b^2 \in (0, -\l]$, as illustrated by Figure \ref{fig10}, the unique admissible trajectory providing us with this type of solutions is the exterior one, which has been plotted in green. Similarly, having a look at Figure \ref{fig11}, as soon as $\b^2 > -\l$ there are four admissible trajectories that have been also represented in green.

\par
By adapting the argument already given in the case $\l >0$, it is easily seen that the blow-up time of the solution of \eqref{3.15} is also given by \eqref{3.16} with $u_0>0$. Thus, $\mathcal{T}_{\mathcal{R}_{\b}}(u_0)$ is continuous and decreasing with respect to $u_0>0$ and it satisfies
\begin{equation}
\label{3.18}
\lim_{u_0 \downarrow 0} \mathcal{T}_{\mathcal{R}_{\b}}(u_0) = + \infty, \qquad \lim_{u_0 \uparrow +\infty} \mathcal{T}_{\mathcal{R}_{\b}}(u_0) =0.
\end{equation}
Consequently,  in both cases, the graph of $\mathcal{T}_{\mathcal{R}_{\b}}(u_0)$ looks like the one shown in Figure \ref{fig6}, replacing $u_0^*$ by $u_{-}$ in the left figure. Note that $\mathcal{T}_{\mathcal{R}_{\b}}(u_0)$ satisfies \eqref{3.17} if $\l>0$ and \eqref{3.18} if $\l\leq 0$. Therefore, for every $R >0$, there exists a unique $u_0 = u_0(R) >0$ such that
$$
\mathcal{T}_{\mathcal{R}_{\b}}(u_0) = R.
$$
Moreover, $u_0(R) > u_{-}$ if $\l >0$. As a direct consequence, for every $\l \in \mathbb{R}$, the
problem \eqref{3.14} has a unique positive solution. This ends the proof of Part (a) for Robin boundary condition with $\b >0$.
 \par
The solution of \eqref{3.14} can have its trajectory on any of the green integral curves plotted in Figure  \ref{fig8}. In particular, its trajectory might be a piece of the half branch of the unstable manifold of $u_0^*$ within the first quadrant. Indeed, this can be reached by imposing
$$
\mathcal{T}_{\mathcal{R}_{\b}}(u_+) \equiv  \int_1^{+\infty}  \dfrac{d \theta}{\sqrt{
\dfrac{2a}{p+1}u_+^{p-1}\left(\theta^{p+1}-1\right) - \l \left(\theta^2-1\right)  +\beta^2}} =R,
$$
where it should not be forgotten that $u_+$ is, actually, a function of $\b$, $u_+=u_+(\b)$, which is
increasing with respect to $\b$. We claim that there is a unique value of $\b>0$ satisfying this identity. Indeed, by continuous dependence,
$$
  \lim_{\b\da 0} \mathcal{T}_{\mathcal{R}_{\b}}(u_+(\b)) = \mathcal{T}_{\mathcal{N}}(u_0^*)=+\infty.
$$
Moreover,
$$
  \lim_{\b\ua +\infty} \mathcal{T}_{\mathcal{R}_{\b}}(u_+(\b)) = 0.
$$
Therefore, since $\mathcal{T}_{\mathcal{R}_{\b}}(u_+(\b))$ is decreasing with respect to $\b$, there is a unique
value of $\b>0$ for which the previous claim holds. This argument can be also adapted to cover the case  $\l\leq 0$ and $\b^2 >-\l$.

\subsection{The Robin case $\mathcal{B} = \mathcal{R}_{\beta}$ with $\beta < 0$}
\label{sec:3.4}
In this section, we will focus the attention on the problem \eqref{eq:1.1} for the choice $\mathcal{B} = \mathcal{R}_{\beta}$ with $\beta <0$,
\begin{equation}
	\label{3.19}
	\left \{ \begin{array}{ll}
		-u''=\lambda u-a |u|^{p-1}u\quad \hbox{in} \;\; [0,R),\\[1ex]
		u'(0)=\beta u(0),\quad        u(R) = +\infty.          \end{array} \right.
\end{equation}
In this case, the line $v=\b u$ has a negative slope. This will entail some substantial differences with
respect to the case $\b>0$. Also in the case  $\b<0$, we will distinguish two cases according to the sign of $\l\in\R$.
\par
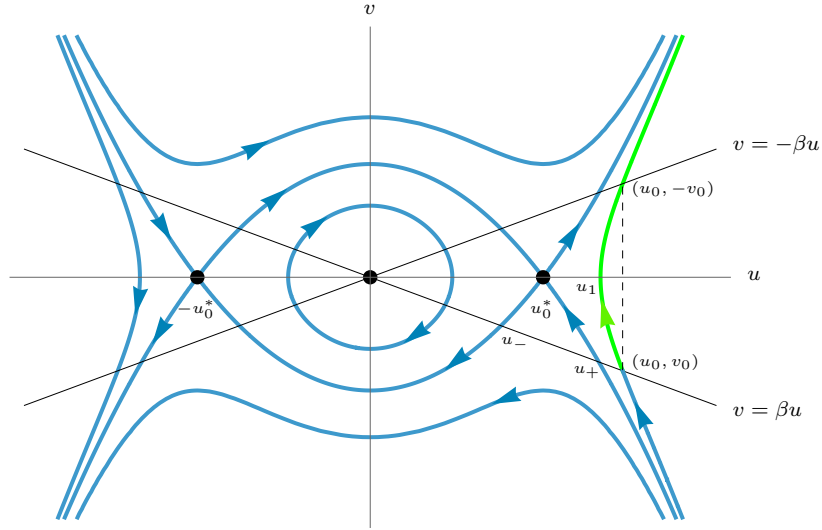
\begin{figure}[h!]
	\centering
	\begin{overpic}[scale=0.75]{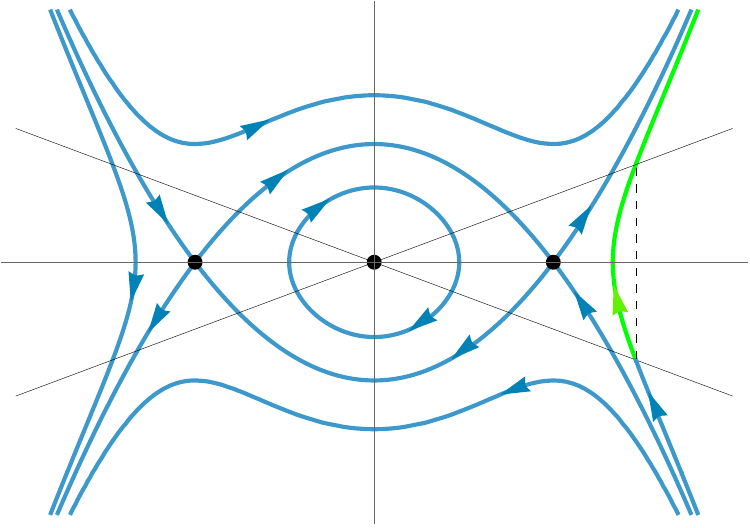}
		\put (86,47) {\tiny$(u_0,-v_0)$}
		\put (86,23) {\tiny$(u_0,v_0)$}
		\put (100,16) {\scriptsize $v = \beta u$}
		\put (100,53) {\scriptsize $v = -\beta u$}
		\put (72,30.4) {\tiny$u_0^*$}
		\put (68,26) {\tiny$u_{-}$}
		\put (78,22) {\tiny$u_{+}$}
		\put (78.3,33.4) {\tiny$u_{1}$}
		\put (23.2,30.4) {\tiny$- u_0^*$}
		\put (102,35) {\scriptsize$u$}
		\put (49,72) {\scriptsize$v$}
	\end{overpic}
	\caption{Phase diagram of the equation \eqref{eq:3.1} for $\l>0$ taking $\mathcal{B} = \mathcal{R}_{\beta}$ with $\beta < 0$. The green trajectory corresponds to a positive solution satisfying the condition $u'(0)=\b u(0)$ with $\b<0$.}
	\label{fig12}
\end{figure}

\medskip
\noindent \emph{Case $\l >0$.}  Figure \ref{fig12} shows the phase diagram of \eqref{eq:3.1} superimposed
with the straight lines $v=\pm \b u$, as well as the type of trajectories that might provide us with a
positive solution of \eqref{3.19}, i.e. the one completely lying in $\{u>0\}$, starting at a point $(u_0,v_0)$ with $v_0=\b u_0<0$, which is unbounded and completely traveled in time $R$, whose trajectory has been plotted in green in Figure \ref{fig12}.  As in the previous case, the line $v=\b u$ crosses the energy level
of the heteroclinic connections,
$$
    E(u,v)=E(\pm u_0^*,0),
$$
at exactly two different points, $u_-<u_+$, with $u>0$. Note that this energy level consists of the
equilibria $(\pm u_0^*,0)$ together with their stable and unstable manifolds, and that, necessarily, $u_0>u_+$.
\par
According to the Cauchy--Lipschitz theory, for every $u_0>u_+$, the Cauchy problem
\begin{equation}
	\label{3.20}
	\left \{ \begin{array}{ll}
		-u''=\lambda u-a |u|^{p-1}u,\\[1ex]
		u(0)=u_0,\quad 	u'(0)=v_0=\beta u_0, \end{array} \right.
\end{equation}
has a unique (maximal) solution $u$, defined in some (maximal) interval $I=[0,T_\mathrm{max})$ for some $T_\mathrm{max}\leq +\infty$, and
$$
  \lim_{t\ua T_\mathrm{max}}u(t)=+\infty\;\; \hbox{if}\;\; T_\mathrm{max}(u_0)<+\infty.
$$
Arguing as in the previous cases, it readily follows that
$$
   \mathcal{T}_{\mathcal{R}_{\b}}(u_0)\equiv T_\mathrm{max}(u_0) <+\infty.
$$

Indeed, as it is apparent from Figure \ref{fig12}, taking into account the symmetries of the phase portrait, for every $u_0>u_+$, such a time it is given by
\begin{equation}
	\label{3.21}
	\mathcal{T}_{\mathcal{R}_{\b}}(u_0)= 2\mathcal{T}(u_0) + \mathcal{T}_{\mathcal{R}_{-\b}}(u_0)
\end{equation}
where $\mathcal{T}\equiv \mc{T}(u_0)$ is the time that the solution of \eqref{3.20} takes to reach the unique point
$$
   (u_1,0)=(u_1(u_0),0), \quad u_1>0,
$$
such that
$$
   E(u_0,v_0)=E(u_1,0),
$$
and $\mathcal{T}_{\mathcal{R}_{-\b}}(u_0)$ is the blow-up time
of the Robin case for $-\b >0$, starting at $(u_0,-v_0)$, which is given by \eqref{3.16}  (see Figure \ref{fig12}).
Since $\mathcal{T}_{\mathcal{R}_{-\b}}(u_0)<+\infty$, it follows from \eqref{3.21} that $\mathcal{T}_{\mathcal{R}_{\b}}(u_0)<+\infty$, because any compact arc of integral curve without equilibria must
be run in a finite time and, hence, $\mc{T}<+\infty$. Actually,
\begin{equation}
\label{3.22}
\mathcal{T}(u_0) <  \mathcal{T}_{\mathcal{N}}(u_1(u_0)),
\end{equation}
where $\mathcal{T}_\mathcal{N}$ stands for the blow-up time of the Neumann problem already given by \eqref{3.11}.
\par
Arguing as in Section \ref{sec:3.3}, we obtain that
\begin{align}
	\mathcal{T}(u_0) &= \int_{u_1}^{u_0}  \dfrac{d\xi}{\sqrt{\l \left(u_0^2 - \xi^2\right) -  \dfrac{2a}{p+1}\left( u_0^{p+1} - \xi^{p+1}\right)+ \beta^2 u_0^2}}\\
	&=\int_{\frac{u_1}{u_0}}^{1}  \dfrac{d \theta}{\sqrt{\l \left(1- \theta^2\right) -  \dfrac{2a}{p+1}u_0^{p-1}\left( 1- \theta^{p+1}\right)+ \beta^2}}.
\label{3.23}
\end{align}
Thus, $\mathcal{T}(u_0)$ is continuous with respect to $u_0>u_+$ and, by continuous dependence with respect to initial data,
\begin{equation}
\lim_{u_0 \downarrow u_{+}} \mathcal{T}(u_0) = + \infty.
\end{equation}
Moreover, by the right limit of \eqref{3.12}, it follows from \eqref{3.22} that
\begin{equation}
\label{3.24}
\lim_{u_0 \uparrow +\infty} \mathcal{T}(u_0)= 0,
\end{equation}
because
\begin{equation}
\label{3.25}
  \lim_{u_0\ua +\infty} u_1(u_0) = +\infty.
\end{equation}
Indeed, the energy conservation gives
$$
\frac{\b^2+\l}{2}u_0^2-\frac{a}{p+1}u_0^{p+1}=\frac{\l}{2}u_1^2-\frac{a}{p+1}u_1^{p+1}
$$
and, since $p>1$,  the left-hand side of this relation converges to $-\infty$ as $u_0\ua+\infty$. By examining the right-hand side, for the equality to hold, necessarily \eqref{3.25} is satisfied.
\par
\begin{figure}[h!]
	\centering
	\begin{overpic}[scale=0.60]{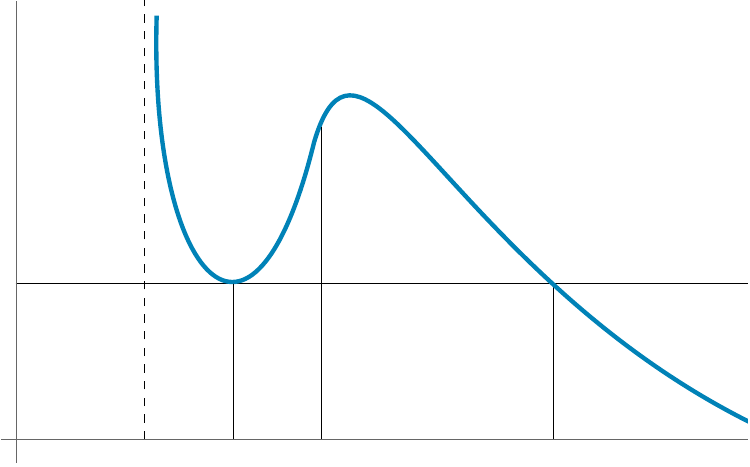}
		\put(-1,-0.5) {\small$0$}
		\put(103,3) {\small$u_0$}
		\put(-13,60) {\small$\mathcal{T}(u_0)$}
		\put(103,23) {\small$S$}
		\put(15.5,-1) {\small$u_{+}$}
		\put(29,-1) {\small$u_{0,1}$}
		\put(40.5,-1) {\small$u_{0,2}$}
		\put(71.5,-1) {\small$u_{0,3}$}
	\end{overpic}
	\caption{Scheme of the proof of the monotonicity of $\mathcal{T}(u_0)$.}
	\label{fig13}
\end{figure}
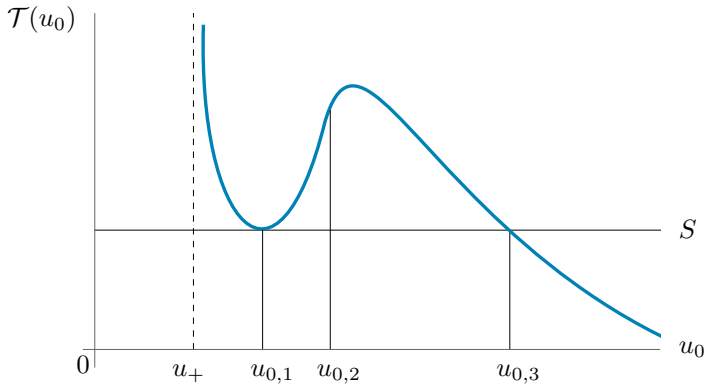

\medskip
The most delicate part of the proof of Part (a) in the case we are dealing with is to show that $\mathcal{T}(u_0)$ is non-increasing with respect to $u_0>u_+$. The proof of this feature relies  on the comparison principles established in Section \ref{sec:2}. Assume, by contradiction, that
$$
  \hbox{there exist}\;\; u_{0,1}<u_{0,2}\;\; \hbox{such that}\;\; \mathcal{T}(u_{0,1})<\mathcal{T}(u_{0,2}).
$$
Then, thanks to \eqref{3.24}, it follows from the continuity of $\mc{T}$, that there is a further $u_{0,3}>u_{0,2}$ such that
\begin{equation}
 S:=\mathcal{T}(u_{0,1})=\mathcal{T}(u_{0,3})<\mc{T}(u_{0,2}),
\end{equation}
as illustrated by Figure \ref{fig13}.

Then, by definition of $\mathcal{T}$, it becomes apparent that the boundary value problem
\begin{linenomath}
	\begin{equation}
		\left \{ \begin{array}{l}
			-u''=\lambda u-a |u|^{p-1}u\quad   \hbox{in} \;\; [0,S],\\[1ex]
			u'(0)=\beta u(0),\quad 			u'(S) =0,          \end{array} \right.
	\end{equation}
\end{linenomath}
has two positive solutions: one with $u(0)=u_{0,1}$, and another one with $u(0)=u_{0,3}$. This contradicts Proposition \ref{pr:2.1}, and concludes the proof that $\mathcal{T}(u_{0})$ is non-increasing as a function of $u_0>u_+$.
Consequently, since $\mc{T}_{\mc{R}_{-\b}}(u_0)$ is continuous and decreasing with respect to $u_0$, we find from \eqref{3.21} that  also $\mc{T}_{\mc{R}_{\b}}(u_0)$ is continuous and decreasing with respect to $u_0>u_+$. Moreover,
thanks to \eqref{3.17} and \eqref{3.24}, we can infer from \eqref{3.21} that
\begin{equation}
\label{3.26}
	\lim_{u_0 \downarrow u_+} \mathcal{T}_{\mathcal{R}_{\b}}(u_0) = + \infty, \qquad 		\lim_{u_0 \uparrow +\infty} \mathcal{T}_{\mathcal{R}_{\b}}(u_0) = 0.
\end{equation}
Therefore, by \eqref{3.26}, for every $R>0$, there exists a unique $u_0>u_+$ such that $\mathcal{T}_{\mathcal{R}_{\b}}(u_0) =R$. So,
\eqref{3.14} also has a unique positive solution for each $\l>0$ and $\b<0$.
Note that the graph of $\mathcal{T}_{\mathcal{R}_{\b}}(u_0)$ looks like the one shown in the left plot of Figure \ref{fig6}, replacing $u_0^*$ by $u_{+}$.
\par
\begin{figure}[h!]
\centering
\begin{tabular}{ccc}
	\begin{overpic}[scale=0.39]{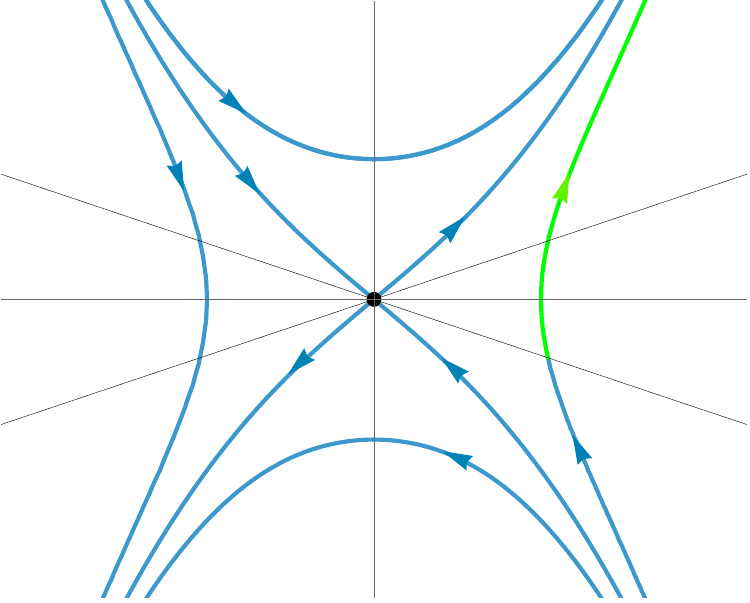}
		\put (73.5,34) {\tiny$(u_0,v_0)$}
		\put (66.5,42) {\tiny$u_{1}$}
		\put (95,19) {\scriptsize $v = \beta u$}
		\put (95,59) {\scriptsize $v = -\beta u$}
		\put (101,39.5) {\scriptsize$u$}
		\put (49,82) {\scriptsize$v$}
	\end{overpic} & 	
\hspace{3mm}
	\begin{overpic}[scale=0.39]{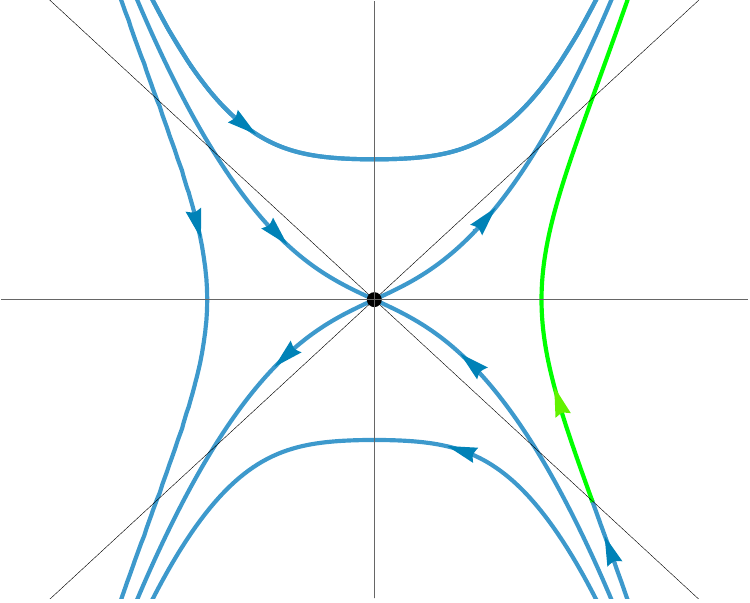}
		\put (80,14) {\tiny$(u_0,v_0)$}
		\put (66.5,37) {\tiny$u_{1}$}
		\put (67.5,18) {\tiny$\tilde{u}$}
		\put (95,0) {\scriptsize $v = \beta u$}
		\put (95,80) {\scriptsize $v =- \beta u$}
		\put (101,39.5) {\scriptsize$u$}
		\put (49,82) {\scriptsize$v$}
	\end{overpic} & 	
\hspace{3mm}
	\begin{overpic}[scale=0.39]{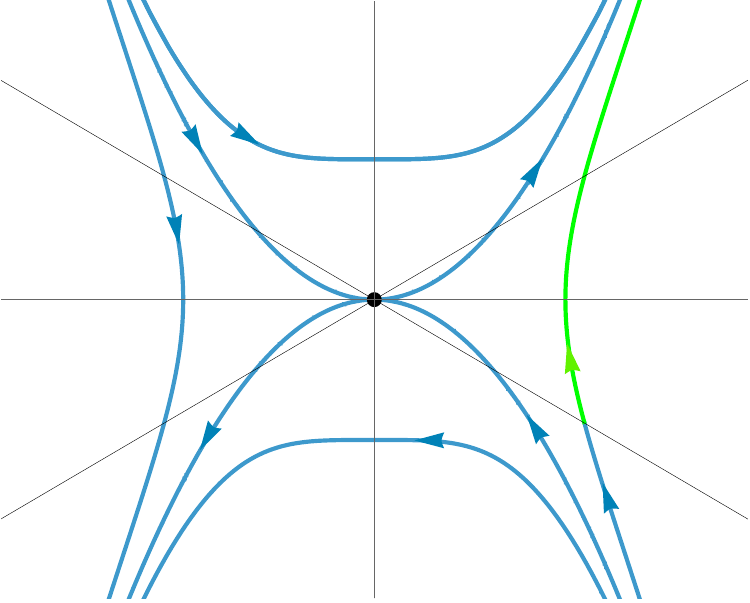}
		\put (79,24) {\tiny$(u_0,v_0)$}
		\put (69.5,37) {\tiny$u_{1}$}
		\put (62,28) {\tiny$\tilde{u}$}
		\put (95,6) {\scriptsize $v = \beta u$}
		\put (95,71) {\scriptsize $v = -\beta u$}
		\put (101,39.5) {\scriptsize$u$}
		\put (49,82) {\scriptsize$v$}
	\end{overpic}
\end{tabular}
\caption{Phase diagrams of \eqref{eq:3.1}, for $\l \leq 0$. The left plot shows the case in which $\l < 0$
and condition \eqref{3.27} does not hold. The central plot shows an admissible case with $\l < 0$ satisfying
\eqref{3.27}. The right plot shows the case $\l=0$. The green trajectories correspond to a positive solution satisfying the condition $u'(0)=\b u(0)$ with $\b<0$.}
\label{fig14}
\end{figure}
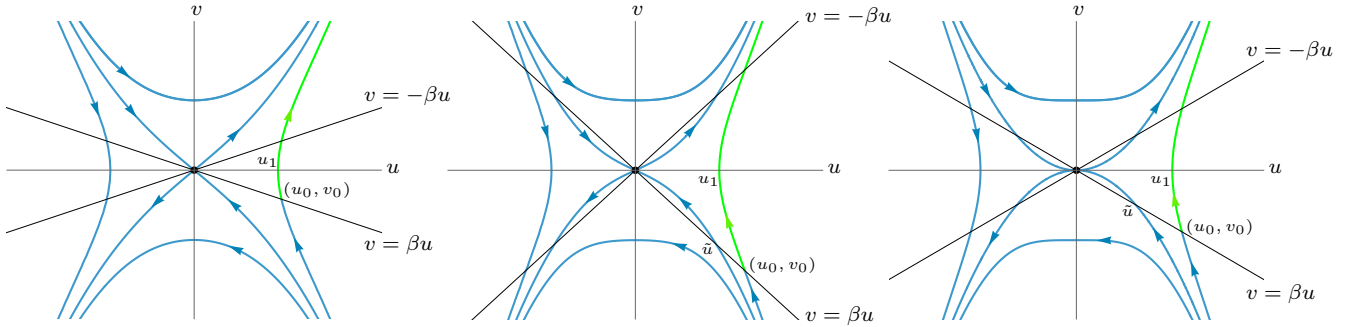

\medskip
\noindent \emph{Case $\l \leq 0$.} Figure \ref{fig14} shows the phase diagram of \eqref{eq:3.1} in this case. Among all the trajectories represented in Figure \ref{fig14} we are interested in those completely lying in $\{u>0\}$, starting at a point $(u_0,v_0)$ with $v_0=\b u_0<0$, and that are unbounded and completely traveled in time $R$. In this case, we need to establish whether, or not, the straight line $v=\b u$ intersects the integral curve $E(u,v)=E(0,0)$ in other points apart from $(0,0)$. From \eqref{eq:3.3}, the condition $E(u,\b u)=E(0,0)$ can be expressed, equivalently, as
\begin{equation}
	u^2\left(\frac{\b^2+\l}{2}-\frac{a}{p+1}u^{p-1}\right)=0.
	\end{equation}
Hence, if
\begin{equation}
	\label{3.27}
	\beta^2 + \l >0
\end{equation}
the line $v=\b u$ intersects to that integral curve at $(\tilde u,\b \tilde u)$, where
$$
 \tilde u= \left( \tfrac{p+1}{2a}(\b^2+\l)\right)^\frac{1}{p-1}>0,
$$
while  then the only intersection point between these two curves is $(0,0)$ if $\b^2+\l\leq 0$. In such case,
for notational convenience, we will set $\tilde u:=0$.  Note that \eqref{3.27} holds for all $\b<0$ if $\l = 0$.

In any circumstances, $u_0 > \tilde u$ is a necessary condition for the existence of a positive solution for the singular problem \eqref{3.19}. By adapting the argument already given in case $\l >0$, it becomes apparent that, for every $u_0>\tilde u$, the unique solution of the Cauchy problem \eqref{3.20}  blows up after a finite time given by \eqref{3.21} with $\mc{T}(u_0)$ given though \eqref{3.23}. Thus, also in this case, $\mathcal{T}_{\mathcal{R}_{\b}}(u_0)$ is continuous and decreasing with respect to $u_0> \tilde u$ and satisfies
\begin{equation}
\label{3.28}
\lim_{u_0 \downarrow \tilde u} \mathcal{T}_{\mathcal{R}_{\b}}(u_0) = + \infty, \qquad \lim_{u_0 \uparrow +\infty} \mathcal{T}_{\mathcal{R}_{\b}}(u_0) =0.
\end{equation}
Consequently, the graph of $\mathcal{T}_{\mathcal{R}_{\b}}(u_0)$ in the case $\l\leq 0$ looks like the one shown in the left picture of Figure \ref{fig6}, replacing $u_0^*$ by $\tilde u$. Therefore, for every $R >0$, there exists a unique
$u_0 = u_0(R) >\tilde u$ such that
$$
\mathcal{T}_{\mathcal{R}_{\b}}(u_0) = R.
$$
Summarizing, for every $\l \in \mathbb{R}$, the problem \eqref{3.19} admits a unique positive solution. This ends the proof of Part (a) for Robin boundary condition with $\b <0$.

\section{Proof of Part (b) of Theorem \ref{th:1.1}}
\label{sec:4}
The proof is based on the construction of the unique solution, $L_\l$,  of \eqref{eq:1.1}, for $a(x)\equiv a$ constant, already carried out in Section \ref{sec:3} for all $\l\in\R$, on the comparison results of Section \ref{sec:2}, as well as on the next pivotal result.

\begin{proposition}
\label{pr:4.1}
For every $M >0$, the boundary value problem
\begin{equation}
\label{4.1}
		\left \{ \begin{array}{l}
			- \theta''=\lambda \theta -a |\theta|^{p-1} \theta\quad \hbox{in} \;\; [0,R],\\[1ex]
			\mc{B}\theta(0)=0,\quad			\theta(R) = M >0,          \end{array} \right.
\end{equation}
has a unique positive solution, $\theta_M$. Moreover, the map $M \mapsto \theta_M$ is increasing, in the sense that
\begin{equation}
\label{4.2}
  \theta_{M_1}(x) < \theta_{M_2}(x)\;\;\hbox{for all}\;\;x \in (0,R] \;\; \hbox{if}\;\; 0<M_1 < M_2.
\end{equation}
Furthermore, besides \eqref{4.2}, the following relations are satisfied as soon as $0<M_1 < M_2$:
\begin{equation}
\label{4.3}
   \left\{ \begin{array}{ll} \theta_{M_1}(0) < \theta_{M_2}(0) & \;\; \hbox{if}\;\; \mc{B}\in \{\mc{N},\mc{R}_\b:\b\in\R\},\\[1ex]
     \theta_{M_1}'(0) < \theta_{M_2}'(0) & \;\; \hbox{if}\;\; \mc{B}= \mathcal{D}. \end{array}\right.
\end{equation}
Finally, if we make the dependence of $\t_M$ on $\l$ explicit by setting $\t_M=\t_{M,\l}$, then the map $\l\mapsto \t_{M,\l}$ is increasing, in the sense that
\begin{equation}
\label{4.4}
  \theta_{M,\l_1}(x) < \theta_{M,\l_2}(x)\;\; \hbox{for all}\;\; x \in (0,R)\;\;\hbox{if}\;\; \l_1 < \l_2.
\end{equation}
In addition, besides \eqref{4.4}, as soon as $\l_1<\l_2$, one has that
\begin{equation}
\label{4.5}
   \left\{ \begin{array}{ll} \theta_{M,\l_1}(0) < \theta_{M,\l_2}(0) & \;\; \hbox{if}\;\; \mc{B}\in \{\mc{N},\mc{R}_\b:\b\in\R\},\\[1ex]
     \theta_{M,\l_1}'(0) < \theta_{M,\l_2}'(0) & \;\; \hbox{if}\;\; \mc{B}= \mathcal{D}. \end{array}\right.
\end{equation}
\end{proposition}
\begin{proof}
The proof of the existence is based on the method of sub- and supersolution.  By a subsolution of \eqref{4.1} we mean a $\mc{C}^2$ function $\un\t$ in $[0,R]$ such that
$$
	\left \{ \begin{array}{l}
			- \un\theta'' \leq \lambda \un\theta -a |\un\theta|^{p-1} \un\theta\quad \hbox{in} \;\; [0,R],\\[1ex] \mc{B}\un\theta(0)\leq 0,\quad \un\theta(R) \leq M.         \end{array} \right.
$$
Similarly, any function $\bar\t$ of class $\mc{C}^2$ such that
$$
	\left \{ \begin{array}{l}
			- \bar\theta'' \geq \lambda \bar\theta -a |\bar\theta|^{p-1} \bar\theta\quad \hbox{in} \;\; [0,R],\\[1ex] \mc{B}\bar\theta(0)\geq 0,\quad \bar\theta(R) \geq M,          \end{array} \right.
$$
is said to be a supersolution of \eqref{4.1}. Moreover, it is said that $\un\t$ is a strict subsolution (resp. $\bar\t$ is a strict supersolution) of \eqref{4.1} if some of the previous inequalities $\leq$ (resp. $\geq$) is satisfied in the strict sense, i.e. replacing $\leq$ by $\lneq$ (resp. replacing $\geq$ by $\gneq$).
\par
Obviously, $\un\t:=0$ is a strict subsolution of \eqref{4.1} for all $M>0$, regardless the nature of the boundary conditions. Moreover, for a sufficiently large positive constant $C$,
$$
\bar \t:= C \;\; \hbox{whenever}\;\; \mc{B}\in\left\{\mc{D},\mc{N},\mc{R}_\b: \b >0\right\}
$$
provides us with  a positive strict supersolution of \eqref{4.1} such that $\bar \t \geq \underline \t$. In the case $\mc{B}=\mc{R}_\b$ with $\b <0$, the construction of the supersolution is far more delicate. Suppose we are in this case, choose $\e \in (0,R)$, and consider a function
$z \in \mc{C}^2([0,\varepsilon];\R)$ such that $z'(0)<0$, and
$$
\varphi(x) := e^{\kappa z(x)} \; \;\hbox{for all}\;\; x\in [0,\e], \;\;\hbox{where}\;\; \kappa > \frac{\b}{z'(0)} >0.
$$
Now, pick any function $\psi \in \mc{C}^2([\varepsilon,R];\R)$ such  that $\psi(x) >0$ for all $x \in [\e,R]$ and
$$
\varphi(\varepsilon) = \psi(\varepsilon), \qquad \varphi'(\varepsilon) = \psi'(\varepsilon), \qquad   \varphi''(\varepsilon) = \psi''(\varepsilon),
$$
so that the function
\begin{equation}
  \Phi(x) := \left\{ \begin{array}{ll}  \varphi(x)  & \;\; \hbox{if}\;\; x \in [0, \varepsilon],\\[1ex]
     \psi(x)  & \;\; \hbox{if}\;\; x \in (\varepsilon, R], \end{array}\right.
\end{equation}
be of class $\mc{C}^2$ in $[0,R]$. We claim that $\bar \t:=C\Phi$ is a positive supersolution of
\eqref{4.1} for sufficiently large $C>0$. In particular, it satisfies $\un\t \leq \bar\t$ in $[0,R]$.
To prove it, at $x=0$ we should have that
$$
  -\bar\t'(0)+\b \bar \t(0)\geq 0,
$$
or, equivalently,
$$
  -C\kappa z'(0)e^{\kappa z(0)}+\b C e^{\kappa z(0)}\geq 0,
$$
which can be attaint by choosing any function $z(x)$ with $z'(0)<0$ and
\begin{equation}
   \kappa \geq \frac{\b}{z'(0)}>0.
\end{equation}
In the interval $[0,R]$ we should have
$$
  -\bar \t''\geq \l \bar \t-a \bar \t^p.
$$
Thus, in $[0,\e]$, we should impose
\begin{equation}
\label{4.6}
  a C^{p-1} e^{(p-1)\kappa z(x)}\geq \l +\kappa z''(x)+\kappa^2 (z'(x))^2 \quad
  \hbox{for all}\;\; x\in [0,\e],
\end{equation}
while in $[\e,R]$, the following inequality is required
\begin{equation}
\label{4.7}
  a C^{p-1} \psi^{p-1}\geq \psi'' +\l \psi.
\end{equation}
Since $z \in \mc{C}^2([0,\e];\R)$, it is easily
realized that the inequality \eqref{4.6} holds for sufficiently large $C>0$. Similarly, since $\psi\in \mc{C}^2([\e,R];\R)$ and it is bounded away from zero, \eqref{4.7} is as well satisfied for sufficiently large $C>0$.
\par
Finally, at $x=R$ we should impose that $\bar \t(R)\geq M$, i.e. $C \psi(R)\geq M$, which can be attaint by
enlarging $C$, if necessary. Therefore, according to a classical result of Amann \cite{Am},  the problem
\eqref{4.1} admits, at least, one positive solution, $\t$, such that $\t\leq \bar\t$.
The uniqueness of the positive solution of \eqref{4.1} is a direct consequence of Proposition \ref{pr:2.1}, and  \eqref{4.2} and \eqref{4.3} follow from Proposition \ref{pr:2.2}.
\par
For the monotonicity in $\l$, observe that, if $\l_1<\l_2$, then $\t_{M,\l_2}$ is a positive strict supersolution of \eqref{4.1} with $\l = \l_1$. Therefore, applying Proposition \ref{pr:2.2} to the equation \eqref{2.3} for $\l = \l_1$, with $v = \t_{M,\l_1}$ and $u= \t_{M,\l_2}$, \eqref{4.4} and \eqref{4.5} hold. The proof is complete.
\end{proof}
\par
\medskip
As a next step, we prove that the solution $\t_{M}$ constructed in the previous lemma converges to the unique large solution $L_\l$ of problem \eqref{eq:1.1} as $M\ua +\infty$. The main ingredients of the proof
are the monotonicity in $M$ and the bounds obtained in the proof of the following result.

\begin{proposition}
\label{pr:4.2}
Let $\t_M$ be the unique positive solution of \eqref{4.1}. Then,
\begin{equation}
\label{4.8}
	\lim_{M\ua+\infty}\t_M=L_\l \quad \hbox{uniformly in compact subsets of} \;\; [0,R),
\end{equation}
where $L_\l$ is the unique positive solution of \eqref{eq:1.1}.
\end{proposition}
\begin{proof}
The existence of the point-wise limit in \eqref{4.8} follows from the monotonicity with respect to $M$ established by Proposition \ref{pr:4.1}. Now, we will show that the limit is finite in $[0,R)$. Since
$$
   L_\l(R)\equiv \lim_{x\ua R} L_{\l}(x)=+\infty,
$$
we have that, for every $M>0$, there exists $\e = \e(M)>0$ such that
\begin{equation}
		\theta_M \leq \max_{x\in[0,R]}{\theta_M}(x) < L_{\l} \qquad \text{in $[R- \varepsilon, R)$}.
\end{equation}
On the other hand, since $L_\l$ is a strict supersolution of
\begin{equation}
		\left \{ \begin{array}{l}
			- \theta''=\lambda \theta -a |\theta|^{p-1} \theta\quad \hbox{in} \;\; [0,R-\e],\\[1ex]
			\mc{B}\theta(0)=0,\quad \theta(R-\e) = \t_{M}(R-\e),          \end{array} \right.
\end{equation}
and $\t_M$ is a positive solution such that, by construction,
$$
  \t_M(R-\e)<L_\l(R-\e),
$$
it follows from Proposition \ref{pr:2.2}  that, in case  $\mc{B}\neq \mc{D}$,
\begin{equation}
		\theta_M(x) < L_{\l}(x) \quad \hbox{for all} \;\; x\in[0, R- \varepsilon],
\end{equation}
whereas, in case $\mc{B}=\mc{D}$,
$$
  \theta_M(x) < L_{\l}(x) \quad \hbox{for all} \;\; x\in(0, R- \varepsilon],\;\;\hbox{and}
  \;\; \theta_M'(0) < L_{\l}'(0).
$$
Thus,  we have proved that, for every $M>0$,
\begin{equation}
		\theta_M(x) \leq L_{\l}(x) \quad \hbox{for all}\;\; x\in[0, R),
\end{equation}
	and, as a consequence,
	\begin{equation}
		\lim_{M \uparrow +\infty} \theta_M(x) \leq L_{\l}(x) \quad \hbox{for all} \;\; x \in [0,R).
	\end{equation}	
	To show that the limit in \eqref{4.8} coincides with $L_\l$, fix a $\delta\in (0,R)$. As the solutions $\theta_M$ are positive and uniformly bounded above by $L_{\l}$ in $[0, R - \delta]$ for  all $M>0$, we obtain from \eqref{4.1} that $\left|\theta_M''\right|$ is uniformly bounded in $[0, R- \delta]$ for all $M>0$. Moreover, for every $x \in [0,R- \delta]$,
	\begin{equation}
		\left|\theta_M'(x)\right| = \left|\int_0^x \theta_M''(s) \, ds + \theta_M'(0)\right| \leq  \int_0^x \left|\theta_M''(s)\right| \, ds + \left|\t_M'(0)\right| \leq C_1\left(R-\d\right) + C_\l=:C_2,
	\end{equation}
where $C_1$ is a bound of $\left|\theta_M''\right|$ in $[0,R- \delta]$, which does not depend on $M$,  and
$$
	C_\l=\left\{ \begin{array}{ll} L_\l'(0) & \quad \hbox{if}\;\; \mc{B}=\mc{D}, \\[1ex]
    \left|\b\right| L_\l(0) & \quad \hbox{if}\;\; \mc{B}=\mc{R}_\b\;\;\hbox{for some}\;\; \b\in\R.\end{array} \right.
$$
This shows that $\left|\theta_M'\right|$ are bounded in $[0,R- \delta]$ uniformly for $M>0$. Thus, by the Mean Value Theorem, for every  $x,y \in [0, R- \delta]$ with $x >y$, there exists $\xi \in (y,x)$ such that
	\begin{equation}
		\left|\theta_M(x) - \theta_M(y)\right| \leq \left|\theta_M'(\xi)\right|\left(x-y\right)\leq C_2\left(x-y\right).
	\end{equation}
Consequently, the family  $ \left\{ \theta_M \right\}_{M >0}$ is uniformly bounded and equicontinuos in $[0, R- \delta]$ and, therefore, by the Ascoli--Arzel\`a Theorem, there exist $\t_\infty\in\mc C([0,R-\d];\R)$ and a subsequence $ \left\{ \theta_{M_n} \right\}_{n \geq 1}$ such that
	\begin{equation}
		\label{4.9}
		\lim_{n\ua +\infty}M_n=+\infty\;\;\hbox{and}\;\;
\lim_{n \uparrow +\infty} \theta_{M_n} = \theta_{\infty} \quad \text{uniformly in } [0, R- \delta].
	\end{equation}
Note that, thanks to the monotonicity of $M\mapsto\theta_M(x)$ obtained in Proposition \ref{pr:4.1}, \eqref{4.9}
actually entails  that
\begin{equation}
		\lim_{M \uparrow+ \infty} \theta_{M} = \theta_{\infty} \quad \text{  uniformly  in }  [0, R- \delta].
\end{equation}
Once again by the Ascoli-Arzel\`a Theorem, there exist $v_\infty\in\mc C([0,R-\d])$ and a subsequence, relabeled again by $M_n$, $\left\{ \theta_{M_n}' \right\}_{n \geq 1}$ such that
\begin{equation}
\label{4.10}
	\lim_{n \uparrow +\infty} \theta_{M_n}' = v_{\infty} \quad \text{uniformly in } [0, R- \delta].
\end{equation}
We claim that $\t_\infty$ is actually differentiable, and that
\begin{equation}
\label{4.11}
	\theta_{\infty}'(x) = v_{\infty}(x) \quad \text{for all } x \in [0, R - \delta].
\end{equation}
Indeed, since
\begin{equation}
	\theta_{M_n}(x) = \theta_{M_n}(0) + \int_0^x \theta_{M_n}'(s) \, ds \quad \text{for every } x \in [0, R - \delta],
\end{equation}
letting $n \uparrow +\infty$ in this relation and using \eqref{4.9} and \eqref{4.10} leads to
\begin{equation}
	\theta_{\infty}(x) = \theta_{\infty}(0) + \int_0^x v_{\infty}(s) \, ds.
\end{equation}
Therefore, by the Fundamental Theorem of Calculus, we can conclude that $\t_\infty'$ exists in $[0,R-\d]$ and
satisfies \eqref{4.11}. This shows the previous claim.
\par
Next, we will show that $\theta_{\infty}$ solves the differential equation of \eqref{eq:1.1} in $[0,R- \delta]$. Indeed, since for every $x \in [0, R - \delta]$
	\begin{equation}
		\theta_{M_n}'(x) - \theta_{M_n}'(0) = \int_0^x \theta_{M_n}''(s) \hspace{0.3mm} ds = \int_0^x  \left(a\theta_{M_n}^{p}(s)- \lambda \theta_{M_n}(s)\right)  \, ds,
	\end{equation}
	by letting $n \uparrow +\infty$ yields to
	\begin{equation}
		\theta_{\infty}'(x) - \theta_{\infty}'(0) = \int_0^x  \left(a\theta_{\infty}^{p}(s)- \lambda \theta_{\infty}(s)\right)  \, ds.
	\end{equation}
	Hence, $\theta_{\infty}'$ is differentiable in $[0,R-\d]$, and
	\begin{equation}
		\theta_{\infty}''(x) =a\theta_{\infty}^{p}(x)- \lambda \theta_{\infty}(x) \;\; \hbox{for all} \;\;  x \in [0, R - \delta].
\end{equation}
This shows that $\theta_{\infty}$ solves the differential equation of \eqref{eq:1.1} in $[0,R- \delta]$ for all $\d\in (0,R)$.
\par
On the other hand, since $\mc{B}\theta_{M_n}(0) = 0$ for all $n$, letting $n \uparrow +\infty$, together with the pointwise convergence of $\t_{M_n}(0)$ to $\t_{\infty}(0)$ and of $\t_{M_n}'(0)$ to $\t_{\infty}'(0)$,  gives $\mc{B}\theta_{\infty}(0) = 0$. Thus, $\theta_{\infty}$ satisfies the boundary condition of \eqref{eq:1.1} at $x = 0$.
\par
Therefore, since $\d$ was arbitrary, the function $\theta_{\infty}$ satisfies the differential equation of \eqref{eq:1.1} in $[0,R)$ and the boundary condition of \eqref{eq:1.1} for $x = 0$. It only remains to show that $\theta_{\infty}(R) = +\infty$, in the sense that
\begin{equation}
\label{4.12}
	\lim_{x \uparrow R} \theta_{\infty}(x) = +\infty.
\end{equation}
Indeed, fix $M>0$. Since $\t_M(R)=M$ and $\t_M$ is continuous, there exists $\delta = \delta(M) >0$ such that
\begin{equation}
	\theta_M(x) > \frac{M}{2} \quad \text{for all } x \in [R- \delta, R].
\end{equation}
Moreover, the monotonicity of $M\mapsto \t_M$ established in Proposition \ref{pr:4.1} guarantees that, for every $N > M$,
\begin{equation}
	\theta_N(x) > \theta_M(x) > \frac{M}{2} \quad \text{for all } x \in [R- \delta, R].
\end{equation}
Letting  $N\ua +\infty$ in this relation, leads to
\begin{equation}
	\theta_\infty(x) > \frac{M}{2} \quad \text{for all } x \in [R- \delta, R),
\end{equation}
which concludes the proof of \eqref{4.12}. Summarizing, the function $\theta_{\infty}$ solves \eqref{eq:1.1}. By the uniqueness obtained in Theorem \ref{th:1.1} (a), necessarily $\t_\infty = L_{\l}$ in $[0,R)$. This ends the proof.
\end{proof}

We are finally ready to prove Theorem \ref{th:1.1}~(b). In the case of Dirichlet boundary conditions, the precise statement of Theorem \ref{th:1.1}~(b) reads as follows.

\begin{proposition}
\label{pr:4.3}
For every $x\in [0,R)$, the mapping $\l \mapsto L_{\l}(x)$, where $L_\l$ is the unique positive solution of  \eqref{eq:3.2}, is differentiable and increasing with respect to $\l$, in the sense that
\begin{equation}
\label{4.13}
    L_{\l_1}(x) < L_{\l_2}(x) \;\;\hbox{for all}\;\; x \in (0,R)\;\;\hbox{if}\;\; \l_1 < \l_2.
\end{equation}
Moreover, $L_{\l_1}'(0) < L_{\l_2}'(0)$.
\end{proposition}
\begin{proof}
Subsequently we denote by $\mc{T}_\mc{D}(\l,v_0)\equiv \mc{T}_\mc{D}(v_0)$ the blow-up time
of the solution of \eqref{eq:3.4} given in \eqref{3.6}. So, in this proof we are making explicit, in addition,  the dependence of $\mc{T}_\mc{D}$ on $\l$. By the properties analyzed in Section \ref{sec:3}, we already know that, for every
$\l\in\R$, $\mc{T}_\mc{D}(\l,v_0)$ is decreasing with respect to $v_0$ and that
$$
  \lim_{v_0\da v_0^*}\mc{T}_\mc{D}(\l,v_0)=+\infty \; \; \hbox{if} \; \; \l >0,\quad
  \lim_{v_0\da 0}\mc{T}_\mc{D}(\l,v_0)=+\infty \; \; \hbox{if} \; \; \l \leq 0, \quad  \lim_{v_0\ua +\infty}\mc{T}_\mc{D}(\l,v_0)=0.
$$
Thus, for every $\l\in\R$, there exists
a unique $v_{0}(\l)\in (v_0^*,+\infty)$ if $\l >0$, and $v_0(\l) \in (0, +\infty)$ if $\l \leq 0$, such that
 \begin{equation}
 	\label{4.14}
\mathcal{T}_{\mathcal{D}}(\l, v_0(\l)) = R.
\end{equation}
On the other hand, by the  Dominated Convergence Theorem, it is easily seen that
$\mc{T}_\mc{D}(\l,v_0)$ is of class $\mc{C}^1$-regularity with respect to $\l$ and $v_0$. Since $\mathcal{T}_{\mathcal{D}}(\l,v_0)$ is differentiable with respect to $\l$ and $v_0$, and, for every $\l\in\R$ and $v_0\in (v_0^*,+\infty)$ if $\l >0$ and $v_0 \in (0, +\infty)$ if $\l \leq 0$,
$$
  \frac{\p \mc{T}_\mc{D}}{\p v_0}(\l,v_0)=-\int_0^{+\infty} v_0 \left( v_0^2-\l s^2+\frac{2a}{p+1}s^{p+1}\right)^{-\frac{3}{2}}\,ds <0,
$$
the Implicit Function Theorem entails that $v_0(\l)$ is of class $\mc{C}^1$ with respect to $\l$. Thus,
differentiating \eqref{4.14} with respect to $\l$, we find that
\begin{equation}
\label{4.15}
v_0'(\l)=-\frac{\partial \mathcal{T}_{\mathcal{D}}}{\partial \l}(\l,v_0(\l))\left(\frac{\partial \mathcal{T}_{\mathcal{D}}}{\partial v_0}(\l, v_0(\l))\right)^{-1} >0,
\end{equation}
because, according to \eqref{3.6}, we have that
$$
   \frac{\p \mc{T}_\mc{D}}{\p \l}(\l,v_0)= \int_0^{+\infty} \frac{s^2}{2} \left( v_0^2-\l s^2+\frac{2a}{p+1}s^{p+1}\right)^{-\frac{3}{2}}\,ds>0.
$$
By definition of $v_0(\l)$, \eqref{4.15} implies that
$$
   L_{\l_1}'(0)<L_{\l_2}'(0)\;\;\hbox{if}\;\; \l_1<\l_2.
$$
As $L_{\l}$ solves the initial value problem \eqref{eq:3.4} in $[0,R)$ with $v_0=v_0(\l)$, the differentiability of $\l\mapsto L_{\l}(x)$, for each $x\in[0,R)$, follows easily from Peano's Differentiation Theorem (see, e.g.,  Hartman \cite{Hart02}).
\par
Finally, to show \eqref{4.13}, we recall that, due to Proposition \ref{pr:4.1}, we already know that
\begin{equation}
	\label{4.16}
	\t_{M,\l_1}(x)<\t_{M,\l_2}(x) \quad \text{for all } x\in(0,R)
\end{equation}
provided $\l_1<\l_2$, where we are denoting by $\t_{M,\l}$ the unique positive solution of \eqref{4.1}.
Thus, letting  $M\ua+\infty$ in \eqref{4.16} it follows from Proposition \ref{pr:4.2} that
\begin{equation}
 L_{\l_1}(x)\leq L_{\l_2}(x) \qquad \text{for all } x\in[0,R).
\end{equation}
To prove \eqref{4.13}, assume, by contradiction, that there exists a point $\tilde x\in(0,R)$ such that
$$
   L_{\l_1}(\tilde x)= L_{\l_2}(\tilde x)\equiv m>0.
$$
Then, $L_{\l_2}-L_{\l_1}$ has an interior minimum, $0$, at $\tilde x$.
But, it follows from the differential equations that
\begin{align*}
   (L_{\l_2}-L_{\l_1})''(\tilde x) & =-\l_2 L_{\l_2}(\tilde x)  + a L_{\l_2}^{p+1} (\tilde x) +\l_1
   L_{\l_1}(\tilde x) -a L_{\l_1}^{p+1}(\tilde x)\\ & =(\l_1-\l_2)m +a\left[
   m^{p+1}-m^{p+1} \right] = (\l_1-\l_2)m  <0,
\end{align*}
which is impossible. This contradiction completes the proof.
\end{proof}

Instead, in the case of Neumann or Robin boundary condition (observe that $\mc{N}=\mc{R}_0$), Theorem \ref{th:1.1}~(b) reads as follows.

\begin{proposition}
\label{pr:4.4}
For every $x\in [0,R)$, the mapping $\l \mapsto L_{\l}(x)$, where $L_\l$ is the unique positive solution of  \eqref{3.14} with $\b\in\R$, is differentiable and increasing with respect to $\l$, in the sense that
\begin{equation}
\label{4.17}
    L_{\l_1}(x) < L_{\l_2}(x) \;\;\hbox{for all}\;\; x \in [0,R)\;\;\hbox{if}\;\; \l_1 < \l_2.
\end{equation}
\end{proposition}
\begin{proof}
Suppose that $\b\geq 0$, and denote by
$$
   \mc{T}_{\mc{R}_{\b}}(\l,u_0)\equiv \mc{T}_{\mc{R}_{\b}}(u_0)
$$
the blow-up time of the solution of \eqref{3.14}, which was given in \eqref{3.16}. Hence, in this context we are making explicit the dependence of $\mc{T}_{\mc{R}_{\b}}$ on $\l$. By the properties already analyzed in Section \ref{sec:3}, for every
$\l\in\R$, $\mc{T}_{\mc{R}_{\b}}(\l,u_0)$ is decreasing with respect to $u_0$ and
\begin{equation}
   \left\{ \begin{array}{ll}
\lim_{u_0\da u_0^*}\mc{T}_{\mc{R}_{\b}}(\l,u_0)=+\infty & \;\; \hbox{if}\;\; \mc{B}= \mc{R}_0=\mc{N}  \; \; \hbox{and} \; \; \l >0,\\[1ex]
\lim_{u_0\da u_{-}}\mc{T}_{\mc{R}_{\b}}(\l,u_0)=+\infty & \;\; \hbox{if}\;\; \mc{B}=\mc{R}_{\b} \; \; \hbox{with} \; \; \b > 0 \; \; \hbox{and} \; \; \l >0,\\[1ex]
\lim_{u_0\da 0}\mc{T}_{\mc{R}_{\b}}(\l,u_0)=+\infty & \;\; \hbox{if}\;\; \mc{B} =\mc{R}_{\b} \; \; \hbox{with} \; \; \b \geq 0 \; \; \hbox{and} \; \; \l \leq 0,\\[1ex]
\lim_{u_0\ua +\infty}\mc{T}_{\mc{R}_{\b}}(\l,u_0)=0 & \;\; \hbox{if}\;\; \mc{B} =\mc{R}_{\b} \; \; \hbox{with} \; \; \b \geq 0 \; \; \hbox{and} \; \; \l \in\R.
\end{array}\right.
\end{equation}
Thus, for every $\l\in\R$, there exists a unique $u_{0}(\l)$ such that
 \begin{equation}
 	\label{4.18}
\mc{T}_{\mc{R}_{\b}}(\l, u_0(\l)) = R.
\end{equation}
On the other hand, as a  direct consequence of the Dominated Convergence Theorem, it becomes apparent that
$\mc{T}_{\mc{R}_{\b}}(\l,u_0)$ is of class $\mc{C}^1$-regularity with respect to $\l$ and $u_0$. Since $\mc{T}_{\mc{R}_{\b}}(\l,u_0)$ is differentiable with respect to $\l$ and $u_0$, and, for every $\l\in\R$ and appropriate ranges of $u_0$'s, we have that
$$
  \frac{\p \mc{T}_{\mc{R}_{\b}}}{\p u_0}(\l,u_0)=-a\frac{p-1}{p+1}u_0^{p-2}\int_1^{+\infty} (\t^{p+1}-1) \left( \frac{2a}{p+1}u_0^{p-1}(\t^{p+1}-1)-\l(\t^2-1)+ \b^2 \right)^{-\frac{3}{2}}\,d\t <0,
$$
by the Implicit Function Theorem, the function $u_0(\l)$ is of class $\mc{C}^1$ with respect to $\l$. Thus,
differentiating \eqref{4.18} with respect to $\l$, we find that
\begin{equation}
\label{4.19}
u_0'(\l)=-\frac{\partial \mc{T}_{\mc{R}_{\b}}}{\partial \l}(\l,u_0(\l))\left(\frac{\partial \mc{T}_{\mc{R}_{\b}}}{\partial u_0}(\l, u_0(\l))\right)^{-1} >0,
\end{equation}
because, according to \eqref{3.16}, we have that
$$
   \frac{\p \mc{T}_{\mc{R}_{\b}}}{\p \l}(\l,u_0)= \int_1^{+\infty} \frac{\t^2 -1}{2} \left( \frac{2a}{p+1}u_0^{p-1}(\t^{p+1}-1)-\l(\t^2-1)+ \b^2 \right)^{-\frac{3}{2}}\,d\t >0.
$$
Since $u_0(\l):=L_\l(0)$, \eqref{4.19} implies that
$$
   L_{\l_1}(0)<L_{\l_2}(0)\;\;\hbox{if}\;\; \l_1<\l_2.
$$
As $L_{\l}$ solves the initial value problem \eqref{3.15} in $[0,R)$ with
$$
   u_0=u_0(\l)=L_\l(0),\quad v_0=v_0(\l)=L_\l'(0),
$$
the differentiability of the map $\l\mapsto L_{\l}(x)$, for each $x\in[0,R)$, follows directly from Peano's Differentiation Theorem,  because $v_0(\l) = \b u_0(\l)$ is also of class $\mc{C}^1$ with respect to $\l$. Finally, by adapting the last part of the proof of Proposition \ref{pr:4.3}, \eqref{4.17} holds, which
concludes the proof of Proposition \ref{pr:4.4} for $\b \geq 0$.
\par
Subsequently, we suppose that $\b<0$. In this case, by \eqref{3.21}, we already know that
\begin{equation}
\label{4.20}
	\mathcal{T}_{\mathcal{R}_{\b}}(u_0)= 2\mathcal{T}(u_0) + \mathcal{T}_{\mathcal{R}_{-\b}}(u_0),
\end{equation}
where $\mathcal{T}\equiv \mc{T}(u_0)$ is the time that the solution of \eqref{3.20} takes to reach the unique point
$(u_1,0)=(u_1(u_0),0)$, $u_1>0$, such that $E(u_0,v_0)=E(u_1,0)$ (see Figures \ref{fig12} and \ref{fig14}),
and $\mathcal{T}_{\mathcal{R}_{-\b}}(u_0)$ is the blow-up time
of the Robin case for $-\b >0$ starting at $(u_0,-v_0)$. Thus, differentiating with respect to $u_0$ in \eqref{4.20}
 we find that
$$
  \frac{\p \mathcal{T}_{\mathcal{R}_{\b}}}{\p u_0}(u_0)= 2
  \frac{\p \mathcal{T}}{\p u_0}(u_0) + \frac{\p \mathcal{T}_{\mathcal{R}_{-\b}}}{\p u_0}(u_0)<0
$$
because, according to the analysis of Section \ref{sec:3}, $\frac{\p \mathcal{T}}{\p u_0}(u_0)\leq 0$ and, due to the previous
analysis, we have that $\frac{\p \mathcal{T}_{\mathcal{R}_{-\b}}}{\p u_0}(u_0)<0$ because $-\b>0$.
\par
As the monotonicity  of $\mc{T}(u_0)=\mc{T}(\l,u_0)$ with respect to $\l$ is unclear, to show that
$L_{\l_1}(0)<L_{\l_2}(0)$ if $\l_1<\l_2$ we are going to use a different approach. According to
Propositions \ref{pr:4.1} and \ref{pr:4.2}, we already know that
\begin{equation}
\label{4.21}
  L_{\l_1}\leq L_{\l_2} \;\ \hbox{in}\;\; [0,R) \;\; \hbox{if}\;\; \l_1<\l_2.
\end{equation}
So, arguing by contradiction, suppose that $L_{\l_1}(0)=L_{\l_2}(0)=:C>0$ for some $\l_1<\l_2$.
Then, by \eqref{4.21}, $L_\l(0)=C$ for all $\l \in [\l_1,\l_2]$. Hence, by the nature of the boundary conditions,
$$
	L_{\l}(0)= C>0 \;\; \text{and} \;\; L_{\l}'(0)=\b C<0 \;\; \text{for all} \;\; \l\in[\l_1,\l_2].
$$
Thus, by the differential equation in \eqref{eq:1.1}, we obtain that, for every $\l \in [\l_1,\l_2]$,
\begin{equation}
	-L_{\l}''(0) = \l L_{\l}(0) - a L_{\l}^p(0) = \l C - a C^p.
\end{equation}
Hence, $\l\mapsto L_{\l}''(0)$ is decreasing for every $\l \in [\l_1,\l_2]$. In particular,
\begin{equation}
\label{4.22}
	L_{\l_2}''(0) < L_{\l_1}''(0).
\end{equation}
On the other hand, by Taylor's Theorem, for every $\l\in[\l_1,\l_2]$ and $x$ in a right neighborhood of 0, there exists $ \xi_{\l}=\xi_\l(x) \in (0,x)$ such that
\begin{equation}
	L_{\l}(x) = L_{\l}(0)+ L_{\l}'(0)x + \dfrac{L_{\l}''(\xi_{\l})}{2}x^2=C\left(1+\b x\right)+\dfrac{L_{\l}''(\xi_{\l})}{2}x^2.
\end{equation}
Thanks to \eqref{4.21}, this relation gives
\begin{equation}
	L_{\l_1}''(\xi_{\l_1}(x)) \leq L_{\l_2}''(\xi_{\l_2}(x)).
\end{equation}
Therefore, letting $x \da 0$ in this inequality, we conclude that
\begin{equation}
	L_{\l_1}''(0) \leq L_{\l_2}''(0),
\end{equation}
which contradicts \eqref{4.22} and completes the proof of the proposition also for $\b<0$.
\end{proof}

\section{Proof of Part (c) of Theorem \ref{th:1.1}}
\label{sec:5}
In this section we will study the limiting behavior of $L_{\l}$ as  $\l \uparrow + \infty$. Thus, we assume, without loss of generality, that $\l$ is positive and large enough. Doing the change of variable
\begin{equation}
	u = \lambda^{\frac{1}{p-1}}\, w,
\end{equation}
we have that $u$ is a positive large solution for \eqref{eq:1.1} if an only if $w$ solves
\begin{equation}
		\label{eq:5.1}
		\left \{ \begin{array}{ll}
			-\dfrac{1}{\l}w''= w-a w^p\quad \hbox{in} \;\; [0,R),\\[1ex]
			\mathcal{B}w(0)=0,\quad 	w(R) = +\infty.          \end{array} \right.
\end{equation}
According to the results of Section \ref{sec:3}, \eqref{eq:5.1} has a unique positive solution, $w_\l$.
Subsequently, we denote by $q_\l$ the unique positive solution of the auxiliary problem
\begin{equation}
		\left \{ \begin{array}{l}
			-\dfrac{1}{\l}q''= q-a q^p \quad \hbox{in} \;\; [0,R],\\[1ex]
			\mathcal{B}q(0)=0,\quad 		q(R) =0.          \end{array} \right.
\end{equation}
The solution $q_\l$ exists as long as $\l > \sigma_1[-D^2; \mathscr{B};(0,R)]$ (see, e.g.,
Aleja et al. \cite{AAL}), where we are setting
\begin{equation}
		\mathscr{B}q(x) := \begin{cases}
			\mathcal{B}q(x) & \text{if } x = 0,\\
			q(x) & \text{if } x = R. \end{cases}
\end{equation}
According to \cite[Th. 1.2]{FRLG19}, we have that
\begin{equation}
	\lim_{\l \uparrow +\infty} q_{\l} = \left(\dfrac{1}{a}\right)^{\frac{1}{p-1}}
\end{equation}
uniformly on compact subsets of $(0,R)$ if $\mathcal{B}=\mathcal{D}$, and uniformly on compact subsets of $[0,R)$
if $\mathcal{B} \in \{ \mathcal{N}, \mathcal{R}_{\beta} \}$ with $\beta \in \mathbb{R}$.  By adapting the proof of Corollary \ref{co:2.3} to \eqref{eq:5.1} it is easily seen that
\begin{equation}
	w_{\l}(x) \geq q_{\l}(x) \qquad  \text{for all } x \in [0,R).
\end{equation}
Therefore, recalling that we are denoting by $L_{\l}$ the unique positive solution of \eqref{eq:1.1},
\begin{equation}
	\lim_{\l \uparrow +\infty} L_{\l} = \lim_{\l \uparrow +\infty} \left( \lambda^{\frac{1}{p-1}}\,w_{\l}\right) \geq \lim_{\l \uparrow +\infty} \left( \lambda^{\frac{1}{p-1}}\,q_{\l}\right) =+\infty,
\end{equation}
uniformly on compact subsets of $(0,R)$ if $\mathcal{B}=\mathcal{D}$, and
uniformly on compact subsets of the interval $[0,R)$ if $\mathcal{B} \in \{ \mathcal{N}, \mathcal{R}_{\beta} \}$ with $\beta \in \mathbb{R}$. This concludes the proof of Theorem \ref{th:1.1}~(c).

\section{Proof of Part (d) of Theorem \ref{th:1.1}}
\label{sec:6}
In this section we will study the limiting behavior of $L_{\l}$ as $\l \downarrow - \infty$. The following result is pivotal for the proof.

\begin{figure}[h!]
\centering
\begin{overpic}[scale=0.8]{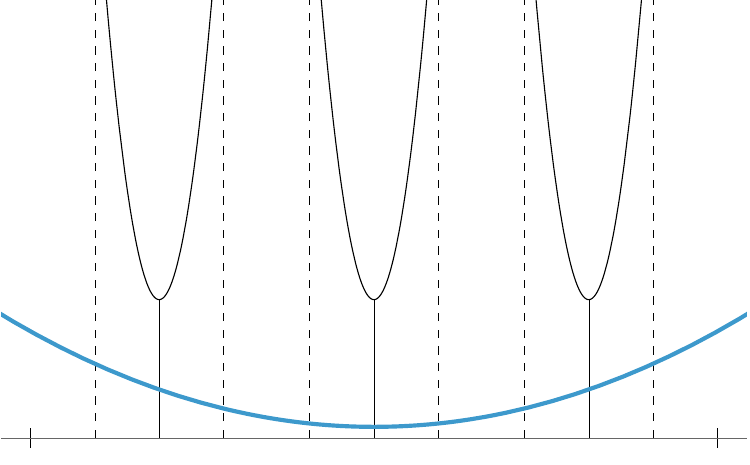}
\put(103,3) {\small$x$}
\put(102,21) {\small$u_\l$}
\put(48.5,-1) {\scriptsize$x_0$}
\put(36.5,-1) {\scriptsize$x_0 - \delta $}
\put(55,-1) {\scriptsize$x_0 + \delta $}
\put(46.5,65) {\small$v_{\l,x_{0}}$}
\put(1,-1) {\scriptsize$c + \varepsilon $}
\put(92.5,-1) {\scriptsize$d - \varepsilon $}
\end{overpic}
\caption{Scheme of the proof of Lemma \ref{le:6.1}, with the large solutions $v_{\l, x_0}$, where $x_0 \in [c+ \varepsilon, d- \varepsilon]$ is arbitrary.}
\label{fig15}
\end{figure}
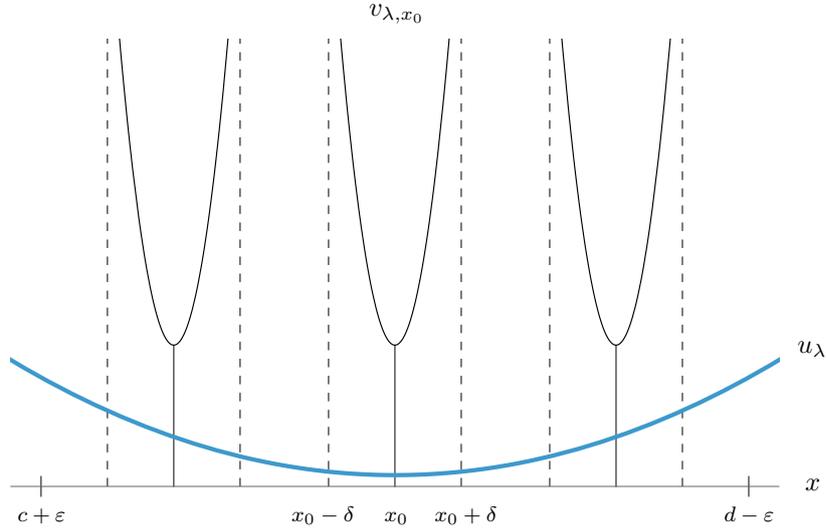

\begin{lemma}
\label{le:6.1}
For every $\l\in\R$, the singular boundary value problem
\begin{equation}
	\label{6.1}
	\left \{ \begin{array}{l}
		-u''=\lambda u-a u^{p}\quad \hbox{in} \;\; (c,d),\\
		u(c)= u(d)= +\infty,    \end{array} \right.
\end{equation}
has a unique positive solution, $u_{\l}$. Moreover, for every $\e>0$,
\begin{equation}
\lim_{\l \downarrow - \infty} u_{\l}= 0 \quad \hbox{uniformly in}\;\; [c+\e,d-\e].
\end{equation}

\end{lemma}
\begin{proof}
The existence and uniqueness of $u_{\l}$ are due to the existence and uniqueness of positive solution for the problem \eqref{3.8}, because $u_{\l}$ is the reflection  of the solution of \eqref{3.8}, considered in $[\frac{c+d}{2},d)$ instead of $[0,R)$, about the midpoint $\frac{c +d}{2}$ of $(c,d)$. Moreover, from the expression of $\mathcal{T}_{\mathcal{N}}(\l, u_0)$ given in \eqref{3.11}, we find that $\mathcal{T}_{\mathcal{N}}(\l,u_0)$ is increasing with respect to $\l$ and, for every $u_0>0$,
\begin{equation}
\label{6.2}
	\lim_{\l \downarrow -\infty} \mathcal{T}_{\mathcal{N}}(\l,u_0) =0.
\end{equation}
With the notation of Section \ref{sec:3}, this implies that
\begin{equation}
\label{6.3}
\lim_{\l \downarrow -\infty} u_0(\l)= \lim_{\l \downarrow -\infty} L_\l(0) =0
\end{equation}
where $L_\l$ is the unique solution of \eqref{3.8}. Indeed,
arguing by contradiction, assume that there exist $\o>0$ and a sequence $\{\l_n\}_{n\geq 1}$ such that
$$
  \lim_{n\ua +\infty}\l_n=-\infty\;\;\hbox{and}\;\; u_0(\l_n)\geq \o\;\;\hbox{for all}\;\; n\geq 1.
$$
Then, by construction, it follows that, for every $n\geq 1$,
$$
  R=\mc{T}_\mc{N}(\l_n,u_0(\l_n))\leq \mc{T}_\mc{N}(\l_n,\o).
$$
Thus, letting $n\ua +\infty$ in this inequality, \eqref{6.2} implies that $R\leq 0$, which is a
contradiction. Hence, \eqref{6.3} holds. Consequently,
\begin{equation}
\label{6.4}
\lim_{\l \downarrow - \infty}u_{\l} \left(\tfrac{c + d}{2}\right) = 0.
\end{equation}
Now, fix $\varepsilon >0$, take a $\delta\in(0,\varepsilon)$, and, for every $x_0 \in [c+ \varepsilon, d - \varepsilon]$, let us denote by  $v_{\l, x_0}$ the unique positive solution of the problem
\begin{equation}
	\left \{ \begin{array}{l}
		-v''=\lambda v-a v^{p} \qquad \hbox{in} \;\;  (x_0 - \delta,x_0 +\delta),\\[1ex]
		v(x_0 - \delta)= v(x_0 + \delta)= +\infty,    \end{array} \right.
\end{equation}
which has been plotted in Figure \ref{fig15}. According to \eqref{6.4},
\begin{equation}
	\label{6.5}
\lim_{\l \downarrow - \infty}v_{\l, x_0}(x_0)= 0.
\end{equation}
Moreover,  Corollary \ref{co:2.4} entails that
\begin{equation}
\label{6.6}
u_{\l}(x) < v_{\l,x_0}(x)\;\;\hbox{for all}\;\; x \in (x_0- \d, x_0 + \d),
\end{equation}
as illustrated by Figure \ref{fig15}.
Therefore, by \eqref{6.5} and \eqref{6.6}, we can infer that
\begin{equation}
\lim_{\l \downarrow - \infty}u_{\l}(x_0)= 0.
\end{equation}
Since $x_0\in[c+\e,d-\e]$ is arbitrary, this shows the point-wise convergence
of $u_\l$ to zero as $\l\da -\infty$ in $(c,d)$.
\par
To show the uniform convergence in compact subsets, fix $\e >0$. Since $u_\l$ is symmetric about $\frac{c+d}{2}$, and it is decreasing in the interval  $\left[c+\e, \frac{c+d}{2}\right]$ and increasing in $\left[\frac{c+d}{2}, d-\e \right]$, we have that
$$
\max_{[c+\e,d-\e]} u_\l = u_\l(c+\e) = u_\l(d-\e).
$$
Therefore, by the point-wise convergence at $c+\e$, we can conclude that
$$
\lim_{\l \da -\infty} \left(\max_{[c+\e,d-\e]} u_\l \right) = 0.
$$
This ends the proof.
\end{proof}
\par
\par
Now, we are ready to complete the proof of Theorem \ref{th:1.1}~(d). Let $u_\l$ be the unique positive solution given by Lemma \ref{le:6.1} with $(c,d)=(0,R)$. Owing to Corollary \ref{co:2.4}, the unique solution of \eqref{eq:1.1}, $L_\l$, satisfies
\begin{equation}
L_{\l}(x) < u_{\l}(x) \quad \hbox{for all} \; \;  x \in (0,R).
\end{equation}
From this relation and Lemma \ref{le:6.1}, we can conclude that
\begin{equation}
	\label{6.7}
\lim_{\l \downarrow - \infty} L_{\l}(x)= 0 \quad \hbox{for all} \; \; x \in (0,R).
\end{equation}
It remains to prove the uniform convergence to zero of $L_\l$ as $\l\da -\infty$ in compact subsets of $[0,R)$. In the case  $\mathcal{B} \in \{ \mc{D}, \mc{N}, \mc{R}_{\b} : \b >0 \}$, according to the construction of $L_\l$ carried out in Section \ref{sec:3}, it becomes apparent that $L_{\l}(x)$ is increasing for all $x \in [0,R)$. Thus, for every $\e\in (0,R)$,
$$
\max_{[0,R-\e]} L_\l = L_\l (R- \e)
$$
and, owing to \eqref{6.7}, we can infer that
$$
\lim_{\l \downarrow - \infty} \left( \max_{[0,R-\e]} L_\l \right)= 0,
$$
which proves Theorem \ref{th:1.1}~(d) for $\mathcal{B} \in \{ \mc{D}, \mc{N}, \mc{R}_{\b} : \b >0 \}$.
\par
\begin{figure}[h!]
	\centering
	\begin{overpic}[scale=0.75]{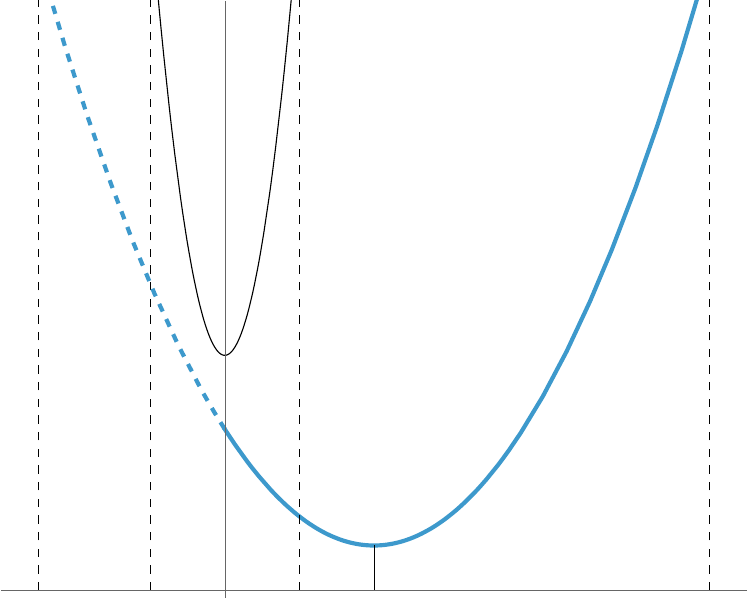}
		\put(103,1) {\small$x$}
		\put(48.5,-3) {\scriptsize$x_1$}
		\put(18,-3) {\scriptsize$-\delta$}
		\put(39,-3) {\scriptsize$\delta$}
		\put(29.2,-3) {\scriptsize $0$}
		\put(42,75) {\small$v_{\l,0}$}
		\put(80,70) {\small$L_{\l}$}
		\put(9,50) {\small$P_{\l}$}
		\put(-1.5,-3) {\scriptsize$2x_1-R$}
		\put(93.5,-3) {\scriptsize$R$}
	\end{overpic}
\vspace{3mm}
	\caption{Scheme of the extended large solution, $P_{\l}$, for $\mathcal{B} = \mathcal{R}_{\beta}$, $\beta < 0$.}
	\label{fig16}
\end{figure}
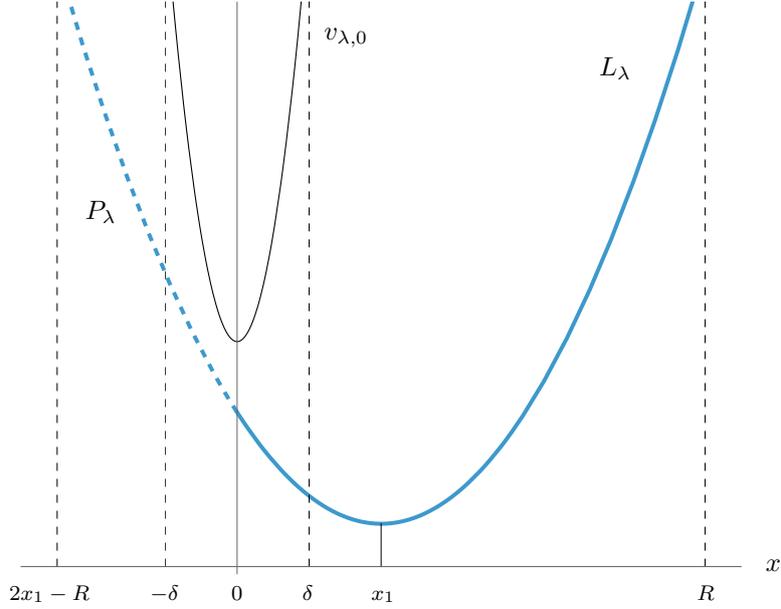

Suppose that $\mathcal{B} = \mathcal{R}_{\beta}$ for some $\beta < 0$. In such a case, according to the analysis
done in Section \ref{sec:3}, we already know the existence of $x_1\in (0,R)$ such that $L_\l$ is decreasing in $[0,x_1]$ and increasing in $[x_1,R)$. Naturally, $L_\l'(x_1)=0$. Next, we consider the following extension of $L_{\l}$ to a large solution in the reflected interval $(2x_1-R,R)$ about $x_1$:
\begin{equation}
	P_\l(x)=\left\{ \begin{array}{ll} L_\l(2x_1-x) & \;\; \hbox{if} \;\; x \in (2x_1-R,0),\\[1ex]
   L_\l(x) & \;\; \hbox{if} \; \;  x \in [0,R),\end{array}\right.
\end{equation}
which has been plotted in Figure \ref{fig16}. $P_\l$ is an even extension of $L_\l$ about $x_1$. Actually, by the phase plane symmetry about $v=0$, $P_\l$ is the unique solution of the Cauchy problem
\begin{equation*}
	\left \{ \begin{array}{l}
		-u''=\lambda u-a |u|^{p-1}u, \\[1ex]
		u(x_1)=L_\l(x_1), \;\; u'(x_1)=0.  \end{array} \right.
\end{equation*}
By uniqueness, $P_\l$ must be the unique positive solution of \eqref{6.1} with $c=2x_1-R$ and $d=R$.
\par
Pick a sufficiently small $\d>0$ such that $(-\d,\d)\subset (2x_1-R,x_1)$,  and let $v_{\l,0}$ be the unique positive solution of
\begin{equation*}
	\left \{ \begin{array}{l}
		-v''=\lambda v-a v^{p} \qquad \hbox{in} \;\;  (-\delta,\delta),\\[1ex]
		v(-\delta)= v(\delta)= +\infty,    \end{array} \right.
\end{equation*}
as illustrated in Figure \ref{fig16}. According to Corollary \ref{co:2.4}, we have that
$P_\l(x)< v_{\l,0}(x)$ for all $x\in (-\d,\d)$. Moreover, by Lemma \ref{le:6.1}, we already know that
$$
  \lim_{\l\da -\infty}v_{\l,0}(x)=0 \;\; \hbox{for all}\;\; x \in (-\d,\d).
$$
Therefore,
\begin{equation}
	\label{6.8}
\lim_{\l \downarrow - \infty} L_{\l}(x)= 0\;\; \hbox{for all}\;\; x\in [0,\d).
\end{equation}
\par
Lastly, for any given $\e >0$, since $L_\l$ is decreasing in the interval $[0, x_1]$ and increasing in $[x_1,R-\e]$, we have that
$$
    \max_{[0,x_1]} L_\l = L_\l (0), \qquad \max_{[x_1,R-\e]} L_\l = L_\l (R- \e).
$$
Thus,
$$
   \max_{[0,R-\e]} L_\l=\max\{L_\l(0),L_\l(R-\e)\}.
$$
Therefore, by \eqref{6.7} and \eqref{6.8}, we can conclude that
$$
\lim_{\l \downarrow - \infty} \left( \max_{[0,R-\e]} L_\l \right) = 0,
$$
which ends the proof of Theorem \ref{th:1.1}~(d).
\par
\medskip

\section{Proof of Theorem \ref{th:1.2}}
\label{sec:7}
Throughout this section, we will assume that $a(x)$ is a continuous function in $[0,R]$ such that $a(x) > 0$ for all $x \in [0,R]$. Thus,
\begin{equation}
	\label{vii.1}
0< a_{\ell} := \min_{[0,R]} a \leq a(x) \leq \max_{[0,R]} a =: a_m \qquad \text{for all } x \in [0,R].
\end{equation}
We will prove Theorem \ref{th:1.2} into several steps.
\par
\medskip
\noindent\emph{Step 1: Existence of positive solutions for problem \eqref{eq:1.1}.} The proof of this fact will rely on the study of the auxiliary boundary value problem
\begin{linenomath}
\begin{equation}
	\label{vii.2}
	\left \{ \begin{array}{ll}
		- \theta''=\lambda \theta -a(x) |\theta|^{p-1} \theta\quad \hbox{in} \;\; [0,R],\\[1ex]
		\mathcal{B}\theta(0)=0,\quad    \theta(R) = M >0.          \end{array} \right.
\end{equation}
\end{linenomath}
Precisely, we will show that, for every $M >0$, \eqref{vii.2} has a unique positive solution, $\theta_{M,a}$. Moreover, for every $M >0$, the map $M \mapsto \theta_{M,a}$ is increasing, in the sense that  $\theta_{M_1,a}(x) < \theta_{M_2,a}(x)$ for each $x \in (0,R]$ if $M_1 < M_2$. Note that, letting $x\da 0$, this also provides us with
$\theta_{M_1,a}(0) \leq \theta_{M_2,a}(0)$. Once proven that, we will show that
\begin{equation}
\label{vii.3}
L_{\l,a}(x):=\lim_{M \uparrow +\infty} \theta_{M,a}(x) \quad\hbox{for all}\;\;  x\in[0,R),
\end{equation}
provides us with a positive solution of \eqref{eq:1.1} such that
\begin{equation}
\label{vii.4}
L_{\l, a_m}(x) \leq L_{\l,a}(x) \leq L_{\l, a_{\ell}}(x) \qquad \text{for all } x \in [0,R),
\end{equation}
where $L_{\l, a_m}$ and $L_{\l, a_{\ell}}$ are, respectively, the unique positive solutions of \eqref{eq:1.1} for $a \equiv a_m$ and $a \equiv a_{\ell}$, whose existence and uniqueness are guaranteed by Theorem \ref{th:1.1}.
\par
To prove the existence of $\theta_{M,a}$, we consider $\theta_{M, a_{\ell}}$ and $\theta_{M, a_m}$, the unique positive solutions of \eqref{vii.2} with $a \equiv a_{\ell}$ and $a \equiv a_m$, respectively, whose existence and uniqueness has been proved in Proposition \ref{pr:4.1}. Thanks to \eqref{vii.1}, $\theta_{M, a_{\ell}}$ and $\theta_{M, a_m}$ are, respectively, a positive supersolution and a positive subsolution of \eqref{vii.2}. Moreover, thanks to Proposition
\ref{pr:2.2}, $\theta_{M, a_m}\leq \theta_{M, a_\ell}$. Thus, by the method of sub and supersolutions (see e.g. Amann
\cite{Am},\cite{AmR}) the problem \eqref{vii.2} has a positive solution, $\theta_{M,a}$, such that
\begin{equation}
\label{vii.5}
\theta_{M, a_m}(x) \leq \theta_{M, a}(x) \leq \theta_{M, a_{\ell}}(x) \quad \text{for all } x \in [0,R].
\end{equation}
The uniqueness of the solution $\t_{M,a}$ is a direct consequence of Proposition \ref{pr:2.1}.
\par
To show the monotonicity of the map $M \mapsto \theta_{M,a}$ with respect to $M$, let $M_1 < M_2$ and consider the corresponding solutions of \eqref{vii.2}, $\theta_{M_1,a}$ and $\theta_{M_2,a}$. Then, since
\begin{equation}
\mathcal{B}\theta_{M_2,a}(0) - \mathcal{B}\theta_{M_1,a}(0) =0, \quad \theta_{M_2,a}(R) - \theta_{M_1,a}(R) =M_2- M_1>0,
\end{equation}
$\t_{M_2,a}$ is a strict positive supersolution of
\begin{linenomath}
\begin{equation}
	\left \{ \begin{array}{ll}
		- \theta''=\lambda \theta -a(x) |\theta|^{p-1} \theta\quad \hbox{in} \;\; [0,R],\\[1ex]
		\mathcal{B}\theta(0)=0,\quad           \theta(R) = M_1 >0,          \end{array} \right.
\end{equation}
\end{linenomath}
and, hence, by Proposition \ref{pr:2.2},
\begin{equation}
\left(\theta_{M_2,a} - \theta_{M_1,a}\right)(x) > 0 \qquad \text{for all } x \in (0,R).
\end{equation}
As
$$
    \theta_{M_2,a}(R) = M_2 > M_1 = \theta_{M_1,a}(R),
$$
the proof that $M \mapsto \theta_{M,a}$ is increasing is complete.
\par
To conclude the proof of this step, we show that $L_{\l,a}$ defined in \eqref{vii.3} provides us with a positive solution of \eqref{eq:1.1}. First of all, we observe that the monotonicity of the map $M \mapsto \theta_{M,a}$ guarantees that the point-wise limit in \eqref{vii.3} is well defined in $[0,R)$. Then, letting  $M \uparrow +\infty$ in \eqref{vii.5}, Proposition \ref{pr:4.2} yields
\begin{equation*}
L_{\l, a_m}(x) \leq L_{\l,a}= \lim_{M \uparrow +\infty} \theta_{M, a}(x) \leq L_{\l, a_{\ell}}(x) \quad \text{for all } x \in [0,R).
\end{equation*}
Thus, the limit in \eqref{vii.3} is finite in $[0,R)$ and \eqref{vii.4} is proved. Moreover, letting $x\ua R$ in this relation shows that
$$
   \lim_{x\ua R} L_{\l,a}(x)=+\infty.
$$
Finally, to show that $L_{\l,a}$ indeed satisfies the differential equation and the boundary condition at $x=0$ in \eqref{eq:1.1}, one can proceed by adapting, \emph{mutatis mutandis}, the proof of Proposition \ref{pr:4.2}, working with $\theta_{M, a}$ instead of $\t_M$.
\par
\medskip
\noindent\emph{Step 2: Minimality of the solution $L_{\l, a}$.} In this step we will show that $L^{\rm min}_{\l,a}=L_{\l, a}$, in the sense that any positive solution $L$ of \eqref{eq:1.1} satisfies
\begin{equation}
\label{7.6}
L_{\l,a}(x) \leq L(x) \quad \text{for all } x \in [0,R).
\end{equation}
Indeed, let $L$ be any positive solution of \eqref{eq:1.1} and consider the solution $\theta_{M, a}$ of \eqref{vii.2} constructed in Step 1. Thanks to Corollary \ref{co:2.3} applied with $c=0$ and $d=R$,
\begin{equation}
\label{7.7}
\theta_{M, a}(x) < L(x)  \quad \hbox{for all } x \in (0,R).
\end{equation}
Moreover,
\begin{equation*}
   \left\{ \begin{array}{ll} \theta_{M, a}(0) < L(0)  & \;\; \hbox{if}\;\; \mc{B}\in \{\mc{N},\mc{R}_\b:\b\in\R\},\\[1ex]
     \theta'_{M, a}(0) < L'(0)& \;\; \hbox{if}\;\; \mc{B} = \mathcal{D}.\end{array}\right.
\end{equation*}
Thus, letting $M \uparrow +\infty$ in \eqref{7.7}  it follows from \eqref{vii.3} that \eqref{7.6} holds.
\par
\medskip
\noindent\emph{Step 3: Existence of the maximal solution.} For convenience, for every $r \in (0, R)$, we will denote by $L^{\rm min}_{\l,a,[0,r)}$ the minimal positive solution of
\begin{equation*}
	\left \{ \begin{array}{l}
		-u''=\lambda u-a(x) |u|^{p-1}u\qquad \hbox{in} \;\; [0,r),\\[1ex]
		\mathcal{B}u(0)=0,\quad     u(r) = +\infty,         \end{array} \right.
\end{equation*}
whose existence can be established adapting the construction carried out in Steps 1 and 2. We claim that
the point-wise limit
\begin{equation}
	\label{7.8}
L^{\rm max}_{\l,a}(x):=\lim_{\e \da 0} L^{\rm min}_{\l,a,[0,R-\e)}(x), \quad x \in [0,R),
\end{equation}
provides us with the maximal positive solution of problem \eqref{eq:1.1}. First, we will show that
$L^{\rm max}_{\l,a}$ is well defined. This relies on the fact that $\e \mapsto L^{\rm min}_{\l,a,[0,R-\e)}$ is increasing. Indeed, suppose that $0 < \e_2 < \e_1$. Then, applying Corollary \ref{co:2.3} with the choices
$$
  (c,d)=(0,R-\e_1), \quad u:=L^{\rm min}_{\l,a,[0,R-\e_1)},\quad v:=L^{\rm min}_{\l,a,[0,R-\e_2)},
$$
it becomes apparent that
\begin{equation}
\label{7.9}
L^{\rm min}_{\l,a,[0,R-\e_2)}(x) < L^{\rm min}_{\l,a,[0,R-\e_1)}(x) \quad \hbox{for all } x \in (0,R- \e_1).
\end{equation}
Moreover, at the left endpoint  of the interval we have that
\begin{equation}
   \left\{ \begin{array}{ll} L^{\rm min}_{\l,a,[0,R-\e_2)}(0) < L^{\rm min}_{\l,a,[0,R-\e_1)}(0)  & \;\; \hbox{if}\;\; \mc{B}\in \{\mc{N},\mc{R}_\b:\b\in\R\},\\[1ex]
     \left(L^{\rm min}_{\l,a,[0,R-\e_2)}\right)'(0) < \left(L^{\rm min}_{\l,a,[0,R-\e_1)}\right)'(0)& \;\; \hbox{if}\;\; \mc{B} = \mathcal{D}.\end{array}\right.
\end{equation}
Thus, the point-wise limit in \eqref{7.8} is well defined. Moreover, such a limit is finite in $[0,R)$. Indeed, by letting $\e_2\da 0$ in \eqref{7.9}, we find that, for every $\e \in (0,R)$,
\begin{equation}
\label{7.10}
L^{\rm max}_{\l,a}(x) \leq L^{\rm min}_{\l,a,[0,R-\e)}(x) \quad \hbox{for all } x \in [0,R- \e).
\end{equation}
Consequently, $L^{\rm max}_{\l,a}$ is well defined and it is uniformly bounded in any compact subinterval of $[0,R)$.
\par
Based on the uniform a priori bounds \eqref{7.10} and on the Ascoli--Arzel\`a Theorem, one can easily adapt the proof of Proposition \ref{pr:4.2} to prove that  $L^{\rm max}_{\l,a}$ is indeed a positive solution of \eqref{eq:1.1}.
\par
Finally, we will show that $L^{\rm max}_{\l,a}$ is maximal, in the sense that any other positive solution $L$ of \eqref{eq:1.1} satisfies
\begin{equation}
\label{7.11}
L(x) \leq L^{\rm max}_{\l,a}(x) \quad \hbox{for all } x \in [0,R).
\end{equation}
Indeed, let $L$ be any positive solution of \eqref{eq:1.1}. Then, for every $\e\in(0,R)$, by applying Corollary \ref{co:2.3} with the choices
$$
   (c,d)=(0,R-\e),\quad u:=L^{\rm min}_{\l,a,[0,R-\e)},\quad v:=L,
$$
we can infer that
\begin{equation}
\label{7.12}
L(x) < L^{\rm min}_{\l,a,[0,R-\e)}(x) \quad \hbox{for all } x \in (0,R-\e).
\end{equation}
Moreover,
\begin{equation}
   \left\{ \begin{array}{ll} L(0) < L^{\rm min}_{\l,a,[0,R-\e)}(0)  & \;\; \hbox{if}\;\; \mc{B}\in \{\mc{N},\mc{R}_\b:\b\in\R\},\\[1ex]
     L'(0) < \left(L^{\rm min}_{\l,a,[0,R-\e)}\right)'(0)& \;\; \hbox{if}\;\; \mc{B} = \mathcal{D}.\end{array}\right.
\end{equation}
Therefore, by letting $\e \da 0$ in \eqref{7.12} and using \eqref{7.8}, \eqref{7.11} holds.
\par
\medskip
\noindent\emph{Step 4: Limiting behavior of $L^{\rm min}_{\l,a}$ as $\l \uparrow +\infty$ and of $L^{\rm max}_{\l,a}$ as $\l \downarrow -\infty$.} According to \eqref{vii.4}, it follows from the analysis already done in Step 2 that
\begin{equation}
\label{7.13}
L_{\l,a_m}(x) \leq L^{\rm min}_{\l,a}(x) \quad \hbox{for all } x \in [0,R).
\end{equation}
Letting $\l \uparrow  +\infty$ in \eqref{7.13}, the property \eqref{i.5} is a direct consequence of Theorem \ref{th:1.1}~(c). Similarly, by considering the problem
\begin{equation*}
	\left \{ \begin{array}{l}
		-u''=\lambda u-a(x) |u|^{p-1}u\qquad \hbox{in} \;\; [0,R-\e),\\[1ex]
		\mathcal{B}u(0)=0,\quad            u(R-\e) = +\infty,         \end{array} \right.
\end{equation*}
and adapting the proof of \eqref{vii.4}, it becomes apparent that
\begin{equation}
\label{7.14}
L^{\rm min}_{\l,a,[0,R-\e)}(x) \leq L_{\l,a_{\ell},[0,R-\e)}(x) \quad \hbox{for all } x \in [0,R-\e),
\end{equation}
where we are denoting by $L_{\l,a_{\ell},[0,R-\e)}$ the unique positive solution of
\begin{equation*}
	\left \{ \begin{array}{l}
		-u''=\lambda u-a_\ell |u|^{p-1}u\qquad \hbox{in} \;\; [0,R-\e),\\[1ex]
		\mathcal{B}u(0)=0,\quad            u(R-\e) = +\infty.         \end{array} \right.
\end{equation*}
Therefore, letting $\e\da 0$ in \eqref{7.14}, it follows from \eqref{7.8} that
\begin{equation}
\label{7.15}
	L^{\rm max}_{\l,a}(x) \leq L_{\l,a_{\ell}}(x) \quad \hbox{for all } x \in [0,R).
\end{equation}
Note that, adapting the argument of Step 3,
$$
 L_{\l,a_{\ell}}^\mathrm{max}:= \lim_{\e\da 0}  L_{\l,a_{\ell},[0,R-\e)}
$$
provides us with the maximal positive solution of \eqref{eq:1.1} with $a\equiv a_\ell$. Thus, by Theorem \ref{th:1.1},
$$
   \lim_{\e\da 0}  L_{\l,a_{\ell},[0,R-\e)} = L_{\l,a_\ell}\quad \hbox{in}\;\; [0,R).
$$
Consequently, \eqref{7.15} holds. Finally, letting $\l \da -\infty$ in \eqref{7.15}, it follows from
Theorem \ref{th:1.1} (d) that \eqref{i.6} holds.
\par
\medskip
\noindent\emph{Step 5: Proof of Theorem \ref{th:1.2} (d).}
Suppose that $a(x)$ is non-increasing and that
$$
   \l \geq 0\;\; \hbox{and}\;\; \mc{B} \in \{ \mc{D}, \mc{N}, \mc{R}_\b : \b <0 \}.
$$
For every $\varepsilon \in (0,R)$, set
\begin{equation}
	\label{7.16}
\rho_{\varepsilon} := \frac{R}{R- \varepsilon} > 1,
\end{equation}
and, for each $x\in[0,R-\varepsilon)$, consider the function
\begin{equation}
	\label{7.17}
\hat{L}_{\e}(x) := \rho_{\e}^{\gamma} L^{\rm min}_{\l,a}\left(\r_{\e} x\right) \quad \hbox{with}\;\;
 \gamma = \frac{2}{p-1}>0.
\end{equation}
Thanks to \eqref{7.16} and \eqref{7.17},  we have that, for every $x \in [0,R)$,
\begin{equation}
\label{7.18}
\lim_{\e \da 0} \hat{L}_\e(x) = L^{\rm min}_{\l,a}(x).
\end{equation}
By the definition of $\hat{L}_\e$, $\hat{L}_{\varepsilon}(R-\varepsilon) = +\infty$ and, for every $x \in [0,R-\varepsilon)$, we have that
\begin{equation}
\label{7.19}
-\hat{L}_{\e}''(x) =  \l \r_\e^2  \hat{L}_{\e}(x) - \rho_{\varepsilon}^{-\gamma p + \gamma +2} a \left( \rho_{\varepsilon}x \right) \hat{L}_{\varepsilon}^p(x).
\end{equation}
Since $\l \geq 0$ and $a$ is non-increasing in $[0,R]$, it follows from \eqref{7.19} that
\begin{equation}
-\hat{L}_{\varepsilon}''(x) \geq \lambda \hat{L}_{\varepsilon}(x) - a(x) \hat{L}_{\varepsilon}^p(x) \quad \hbox{for all } x \in [0,R-\e).
\end{equation}
Moreover, if $\mc{B} \in \{ \mc{D}, \mc{N} \}$, by a direct computation, we find that $\mc{B}\hat{L}_{\e}(0) = 0$. Thus, in such cases, the function $\hat{L}_\e$ satisfies
\begin{equation}
	\left \{ \begin{array}{l}
		-\hat{L}_\e'' \geq \lambda \hat{L}_\e-a\hat{L}_\e^p  \quad \hbox{in} \;\; [0,R-\e),\\[1ex]
		\mathcal{B}\hat{L}_\e(0)= 0,\quad
           \hat{L}_\e(R-\e) = +\infty.         \end{array} \right.
\end{equation}
Hence, applying Corollary \ref{co:2.3} for the choices
$$
   (c,d) = (0, R-\e),\quad u := \hat{L}_\e,\quad v := L^{\rm max}_{\l,a},
$$
we obtain that
\begin{equation}
L^{\rm max}_{\l,a}(x) < \hat{L}_{\e}(x)   \quad \hbox{for all } x \in (0,R-\e).
\end{equation}
Therefore, letting $\varepsilon \downarrow 0$ in this relation, we can conclude from \eqref{7.18} that
\begin{equation}
L^{\rm max}_{\l,a}(x) \leq L^{\rm min}_{\l,a}(x) \quad \hbox{for all } x \in [0,R).
\end{equation}
Combining this estimate with \eqref{i.4}, shows the uniqueness of the positive solution of \eqref{eq:1.1}
under these circumstances, which ends the proof  of Theorem \ref{th:1.2} in these cases.
\par
Finally, when $\mc{B} = \mc{R}_\b$ with $\b \in \R$,
$$
 \mc{R}_\b \hat{L}_{\e}(0) = -\r_\e^{\gamma+1}\left(L^{\rm min}_{\l,a}\right)'(0)+ \b \r_\e^\gamma L^{\rm min}_{\l,a}(0).
$$
Thus, by adding and subtracting the term $\b \r_\e^{\gamma+1} L^{\rm min}_{\l,a}(0)$ in the previous identity,
and taking into account that together $\mc{R}_\b L^{\rm min}_{\l,a}(0) =0$ for all $\b \in \R$, we can infer that, for every $\b \in \R$,
$$
\mc{R}_\b \hat{L}_\e(0) = \b \r_\e^\gamma L^{\rm min}_{\l,a}(0) (1- \r_\e).
$$
Consequently, $\mc{R}_\b \hat{L}_\e(0) >0$ if $\b <0$. Hence, in the case $\b<0$, it turns out that $\hat{L}_{\e}$ satisfies
\begin{equation}
\label{7.20}
	\left \{ \begin{array}{l}
		-\hat{L}_{\e}'' \geq \lambda \hat{L}_{\e}-a (x)\hat{L}_{\e}^p  \quad \hbox{in} \;\; [0,R-\e),\\[1ex]
		\mc{R}_\b \hat{L}_{\e}(0)>0,\quad            \hat{L}_{\e}(R-\e) = +\infty.         \end{array} \right.
\end{equation}
Thus, arguing as above, applying Proposition \ref{pr:2.2} for the choices
$$
   (c,d) = (0, R-\e),\quad u := \hat{L}_\e,\quad v := L^{\rm max}_{\l,a},
$$
we also find that
\begin{equation}
L^{\rm max}_{\l,a}(x) < \hat{L}_{\e}(x)   \quad \hbox{for all } x \in [0,R-\e).
\end{equation}
Therefore, letting $\varepsilon \downarrow 0$, we can conclude from \eqref{7.18} that
\begin{equation}
L^{\rm max}_{\l,a}(x) \leq L^{\rm min}_{\l,a}(x) \quad \hbox{for all } x \in [0,R).
\end{equation}
Combining this estimate with  \eqref{i.4} completes the proof of Theorem \ref{th:1.2}.

\end{document}